\documentclass[oneside,reqno]{amsart}
\usepackage[T1]{fontenc}
\usepackage[latin9]{inputenc}
\usepackage{geometry}
\usepackage{color}
\usepackage{mathrsfs}
\usepackage{mathtools}
\usepackage{enumitem}
\usepackage{bm}
\usepackage{amsbsy}
\usepackage{amstext}
\usepackage{amsthm}
\usepackage{amssymb}
\usepackage{stmaryrd}
\usepackage[unicode=true,pdfusetitle,
 bookmarks=true,bookmarksnumbered=false,bookmarksopen=false,
 breaklinks=false,pdfborder={0 0 1},backref=false,colorlinks=true]
 {hyperref}

\makeatletter
\numberwithin{equation}{section}
\numberwithin{figure}{section}
\theoremstyle{plain}
\newtheorem{thm}{\protect\theoremname}
\theoremstyle{plain}
\newtheorem{conjecture}[thm]{\protect\conjecturename}
\theoremstyle{definition}
\newtheorem{defn}[thm]{\protect\definitionname}
\theoremstyle{definition}
\newtheorem{example}[thm]{\protect\examplename}
\theoremstyle{remark}
\newtheorem{rem}[thm]{\protect\remarkname}
\theoremstyle{plain}
\newtheorem{prop}[thm]{\protect\propositionname}
\theoremstyle{plain}
\newtheorem{lem}[thm]{\protect\lemmaname}

\usepackage{url}
\usepackage{shuffle}
\usepackage{bbm}
\usepackage[geometry]{ifsym}

\makeatother

\providecommand{\conjecturename}{Conjecture}
\providecommand{\definitionname}{Definition}
\providecommand{\examplename}{Example}
\providecommand{\lemmaname}{Lemma}
\providecommand{\propositionname}{Proposition}
\providecommand{\remarkname}{Remark}
\providecommand{\theoremname}{Theorem}

\begin{document}
\address[Minoru Hirose]{Institute for Advanced Research, Nagoya University,  Furo-cho, Chikusa-ku, Nagoya, 464-8602, Japan}
\email{minoru.hirose@math.nagoya-u.ac.jp}
\address[Nobuo Sato]{Department of Mathematics, National Taiwan University, No. 1, Sec. 4, Roosevelt Rd., Taipei 10617, Taiwan (R.O.C.)}
\email{nbsato@ntu.edu.tw}
\title{Block shuffle identities for multiple zeta values}
\author{Minoru Hirose and Nobuo Sato}
\begin{abstract}
In 1998, Borwein, Bradley, Broadhurst and Lison\v{e}k posed two families
of conjectural identities among multiple zeta values, later generalized
by Charlton using his alternating block notation. In this paper, we
prove a new class of identities among multiple zeta values that simultaneously
resolve and generalize these conjectures.
\end{abstract}


\maketitle

\global\long\def\mshuffle{\mathbin{{\mathbin{\raisebox{-.23em}{\SmallDiamondshape}}}}}%
\global\long\def\mdelta{{\widehat{\Delta}}}%
\global\long\def\dshuffle{{\mathbin{\raisebox{-.23em}{\SmallRightDiamond}}}}%

\section{Introduction}

The multiple zeta values are real numbers defined by the multiple
Dirichlet series
\[
\zeta(\boldsymbol{k}):=\sum_{0<m_{1}<\cdots<m_{d}}\frac{1}{m_{1}^{k_{1}}\cdots m_{d}^{k_{d}}}
\]
for a multi-index $\boldsymbol{k}=(k_{1},\ldots,k_{d})$. Here we
assume $k_{1},\ldots,k_{d-1}\geq1$ and $k_{d}>1$ for the convergence.
The sum of the entries $k_{1}+\cdots+k_{d}$ is called the weight
while $d$ is called the depth of the multiple zeta values. In the
case of depth $d=1$, a multiple zeta value is a Riemann zeta value,
and a famous evaluation formula of Euler says that 
\[
\zeta(2)=\frac{\pi^{2}}{6}
\]
which gave a surprising answer to Barsel's problem. Furthermore, Euler
showed
\[
\zeta(4)=\frac{\pi^{4}}{90},\:\zeta(6)=\frac{\pi^{6}}{945},\:\zeta(8)=\frac{\pi^{8}}{9450},\:\zeta(10)=\frac{\pi^{10}}{93555},\:\zeta(12)=\frac{691\pi^{12}}{638512875}
\]
and in general
\[
\zeta(2m)\in\mathbb{Q}\cdot\pi^{2m}
\]
for a positive integer $m$. On the other hand, it is widely believed
that none of the `odd' Riemann zeta values $\zeta(2m+1)$ are rational
multiples of powers of $\pi$, and it is even conjectured that all
of `odd' Riemann zeta values and $\pi$ are algebraically independent.
Let us refer to a multiple zeta value which is a rational multiple
of some power of $\pi$ as an `Eulerian' multiple zeta value. In weight
two, there is only one multiple zeta value $\zeta(2)$ which is Eulerian
by the aforementioned Euler's formula. Then, in weight three, we have
two multiple zeta values $\zeta(3),\zeta(1,2)$ which turn out to
be equal:
\[
\zeta(1,2)=\zeta(3)
\]
as shown also by Euler. Thus, at least conjecturally, neither of them
are Eulerian. In weight four, we have four multiple zeta values $\zeta(4),\zeta(1,3),\zeta(2,2),\zeta(1,1,2)$,
and it is known that
\[
\zeta(4)=\zeta(1,1,2)=\frac{\pi^{4}}{90},\ \zeta(1,3)=\frac{\pi^{4}}{360},\ \zeta(2,2)=\frac{\pi^{4}}{120},
\]
which shows that they are all Eulerian. In higher weights, we have
a infinite series of Eulerian multiple zeta values 
\[
\zeta(\overbrace{2,2\ldots,2}^{n})=\frac{\pi^{2n}}{(2n+1)!}
\]
which is an easy consequence of the infinite product expression of
sine function. Moreover, it is not hard to show that

\begin{equation}
\zeta(\overbrace{2m+2,2m+2,\ldots,2m+2}^{n\text{ times}})\label{eq:Eulerian_kkk}
\end{equation}
is an Eulerian multiple zeta value for $m,n$ non-negative integers,
and by the general `duality' relation of the multiple zeta values,
this also implies that
\begin{equation}
\zeta(\overbrace{\{1\}^{2m},2,\{1\}^{2m},2,\ldots,\{1\}^{2m},2}^{n\text{ times}})\label{eq:Eulerian_kkk_dual}
\end{equation}
is also an Eulerian multiple zeta value for $m,n$ non-negative integers.
These two series of Eulerian multiple zeta values cover $\zeta(2,2)$
and $\zeta(1,1,2)$, but does not cover $\zeta(1,3)$. As a generalization
of $\zeta(1,3)=\frac{\pi^{4}}{360}$, Zagier conjectured a simple
evaluation formula 
\[
\zeta(\{1,3\}^{n})=\frac{2\pi^{4n}}{(4n+2)!}
\]
which was later proved in \cite{BBB97}. A natural question is then
``what are \emph{all} Eulerian multiple zeta values?''. Borwein-Bradley-Broadhurst
performed a numerical investigation on multiple zeta values and discovered
the following remarkable series of conjecturally Eulerian multiple
zeta values:
\begin{conjecture}[{\cite[equation (18)]{BBB97}}]
\label{conj:BBB}For $m,n\geq0$,
\begin{equation}
\zeta(\{\{2\}^{m},1,\{2\}^{m},3\}^{n},\{2\}^{m})\stackrel{?}{=}\frac{1}{2n+1}\frac{\pi^{\mathrm{wt}}}{(\mathrm{wt}+1)!}\quad\left(=\frac{1}{2n+1}\zeta(\{2\}^{{\rm wt}/2})\right)\label{eq:Eulerian_main}
\end{equation}
where $\mathrm{wt}=4(m+1)n+2m$ is the weight of the multiple zeta
values.
\end{conjecture}

Extensive numerical evidence\footnote{The authors have verified this fact up to weight $22$ based on the
data of The Multiple Zeta Values data mine \cite{Datamine}.} shows that it is very likely that (\ref{eq:Eulerian_kkk}), (\ref{eq:Eulerian_kkk_dual})
and (\ref{eq:Eulerian_main}) exhaust all the Eulerian multiple zeta
values. Furthermore, Borwein-Bradley-Broadhurst-Lison\v{e}k \cite{BBBL}
introduced the $Z$-notation
\[
Z(m_{0},\dots,m_{2n}):=\zeta(\{2\}^{m_{0}},1,\{2\}^{m_{1}},3,\{2\}^{m_{2}},\dots,1,\{2\}^{m_{2n-1}},3,\{2\}^{m_{2n}})
\]
and generalized Conjecture \ref{conj:BBB} as follows:
\begin{conjecture}[{Cyclic insertion conjecture, \cite[Conjecture 1]{BBBL}}]
\label{Conj:BBBL1}For $n\geq0$ and $m_{0},\dots,m_{2n}\geq0$,
we have
\[
\sum_{j=0}^{2n}Z(m_{2n-j+1},\dots,m_{2n},m_{0},\dots,m_{2n-j})=\zeta(\{2\}^{2n+\sum_{j=0}^{2n}m_{j}}).
\]
\end{conjecture}

Using the $Z$-notation, they also gave the following conjecture.
\begin{conjecture}[{\cite[Conjecture 2]{BBBL}}]
\label{Conj:BBBL2}For $a_{1},a_{2},a_{3},b_{1},b_{2}\geq0$, we
have
\begin{align*}
 & Z(a_{1},b_{1},a_{2},b_{2},a_{3})+Z(a_{2},b_{1},a_{3},b_{2},a_{1})+Z(a_{3},b_{1},a_{1},b_{2},a_{2})\\
= & Z(a_{1},b_{2},a_{2},b_{1},a_{3})+Z(a_{2},b_{2},a_{3},b_{1},a_{1})+Z(a_{3},b_{2},a_{1},b_{1},a_{2}).
\end{align*}
\end{conjecture}

Conjectures \ref{conj:BBB}, \ref{Conj:BBBL1} and \ref{Conj:BBBL2}
were later generalized by Charlton by his invention of the (alternating)
block notation for the motivic multiple zeta values (\cite{Ch16,Ch21}).
To state his conjectures, let us define some notations. For $a_{1},\dots,a_{k}$
with $a_{1},\dots,a_{k}\in\{0,1\}$, we denote by
\[
I^{\mathfrak{m}}(0,a_{1},\dots,a_{k},1)
\]
the motivic iterated integral of $\frac{dt}{t-a_{1}},\dots,\frac{dt}{t-a_{k}}$
along the straight path from $0$ (with the tangential vector $\frac{d}{dt}$)
to $1$ (with the tangential vector $-\frac{d}{dt}$). We call a (possibly
empty) sequence $(l_{1},\dots,l_{d})$ of positive integers an index,
and say that $(l_{1},\dots,l_{d})$ is admissible if it is empty or
$l_{1}>1$ and $l_{d}>1$. Also, we define the parity (even/odd) of
$(l_{1},\dots,l_{d})$ to be the parity of $\sum_{i=1}^{d}(l_{i}-1)$.
\begin{defn}[{(Alternating) block notation, \cite[Definition 2.2.8]{Ch16}, \cite[Notation 3.6]{Ch21}}]
For an odd index $(l_{1},\dots,l_{d})$, we define
\[
I_{\mathrm{bl}}^{\mathfrak{m}}(l_{1},\ldots,l_{d}):=I^{\mathfrak{m}}(a_{0},\ldots,a_{n+1})
\]
where $n+2=l_{1}+\cdots+l_{d}$ and $a_{i}$'s are defined by $a_{0}=0$
and 
\[
a_{i+1}=\begin{cases}
a_{i} & i\in\{l_{1},l_{1}+l_{2},\ldots,l_{1}+\cdots+l_{d-1}\}\\
1-a_{i} & \mathrm{otherwise}.
\end{cases}
\]
Notice that $a_{n+1}=1$ by the oddness of $(l_{1},\ldots,l_{d})$.
\end{defn}

\begin{example}
For example,
\[
I_{\mathrm{bl}}^{\mathfrak{m}}(3,4,5)=I^{\mathfrak{m}}(\overbrace{0,1,0}^{3},\overbrace{0,1,0,1}^{4},\overbrace{1,0,1,0,1}^{5}).
\]
\end{example}

The block notation can be viewed as a generalization of $Z$-notation
since
\[
Z^{\mathfrak{m}}(m_{0},\dots,m_{2n})=(-1)^{m_{0}+\cdots+m_{2n}}I_{{\rm bl}}^{\mathfrak{m}}(2m_{0}+2,\dots,2m_{2n}+2)
\]
with the motivic version $Z^{\mathfrak{m}}$ of $Z$. With this setup,
Charlton's conjectures are stated as follows.

\begin{conjecture}[{Generalized cyclic insertion conjecture, \cite[Conjecture 2.5.1]{Ch16}\cite[Conjecture 6.1]{Ch21}}]
\label{Conj:Ch1}For $d\geq1$ and $l_{1},\dots,l_{d}\in\mathbb{Z}_{\geq2}$
such that $\sum_{i=1}^{d}(l_{i}-1)$ is odd, we have
\[
\sum_{i=0}^{d-1}I_{\mathrm{bl}}^{\mathfrak{m}}(l_{i+1},\dots,l_{d},l_{1},\dots,l_{i})=\begin{cases}
I_{\mathrm{bl}}^{\mathfrak{m}}(l_{1}+\cdots+l_{d}) & d:\,{\rm odd}\\
0 & d:\,{\rm even}.
\end{cases}
\]
\end{conjecture}

\begin{conjecture}[{Generalization of Conjecture \ref{Conj:BBBL2}, \cite[Conjecture 2.8.2]{Ch16}\cite[Conjecture 8.3]{Ch21}}]
\label{Conj:Ch2}For $n\geq1$ and $a_{1},\dots,a_{n+1},\allowbreak b_{1},\dots,b_{n}\geq2$
such that $\sum_{i=1}^{n+1}a_{i}+\sum_{i=1}^{n}b_{i}$ is even, we
have
\[
\sum_{\sigma\in\mathfrak{S}_{n+1}}{\rm sgn}(\sigma)I_{{\rm bl}}^{\mathfrak{m}}(a_{\sigma(1)},b_{1},a_{\sigma(2)},b_{2},\dots,a_{\sigma(n)},b_{n},a_{\sigma(n+1)})=0.
\]
\end{conjecture}

The purpose of this paper is to prove a new extensive class of identities
among multiple zeta values simultaneously generalizing Conjectures
\ref{Conj:Ch1} and \ref{Conj:Ch2}, which we call \emph{block shuffle}
\emph{identities}. To state the block shuffle identities, we start
with preparing some algebraic settings. Let $\mathfrak{X}\coloneqq\mathbb{Q}\left\langle x_{1},x_{2},x_{3},\dots\right\rangle $
be the free non-commutative algebra generated by the formal variables
$x_{1},x_{2},x_{3},\dots$ over $\mathbb{Q}$. For an index $\Bbbk=(k_{1},\dots,k_{d})$,
we put $x_{\Bbbk}\coloneqq x_{k_{1}}\cdots x_{k_{d}}$. We denote
by $\mathfrak{X}^{+}$ (resp. $\mathfrak{X}_{{\rm ev}}$, $\mathfrak{X}_{{\rm od}}$)
the subspace of $\mathfrak{X}$ spanned by $x_{\Bbbk}$'s for all
non-empty (resp. even, odd) indices $\Bbbk$, and put $\mathfrak{X}_{{\rm ev}}^{+}=\mathfrak{X}^{+}\cap\mathfrak{X}_{{\rm ev}}$.
For $k\geq0$, we define a $\mathbb{Q}$-linear map $s_{k}:\mathfrak{X}\to\mathfrak{X}$
by 
\begin{align*}
s_{k}(x_{i}w) & =x_{i+k}w\quad(w\in\mathfrak{X}),\\
s_{k}(1) & =0.
\end{align*}
We define a $\mathbb{Q}$-bilinear map $\mshuffle:\mathfrak{X}\times\mathfrak{X}\to\mathfrak{X}$
by
\begin{align*}
u\mshuffle1= & 1\mshuffle u=u\quad(u\in\mathfrak{X}),\\
x_{a}u\mshuffle x_{b}v= & x_{a}(u\mshuffle x_{b}v)+x_{b}(x_{a}u\mshuffle v)-s_{a+b}(u\mshuffle v)\quad(u,v\in\mathfrak{X})
\end{align*}
or equivalently
\[
x_{a_{1}}\cdots x_{a_{l}}\mshuffle x_{b_{1}}\cdots x_{b_{m}}=\sum_{n=0}^{l+m}(-1)^{(l+m-n)/2}\sum_{\substack{(f,g)\in\mathscr{D}_{n}}
}x_{c_{1}(f,g)}\cdots x_{c_{n}(f,g)}
\]
where 
\[
\mathscr{D}_{n}\coloneqq\left\{ \left.\substack{\substack{f}
:\{1,\dots,l\}\to\{1,\dots,n\}\\
\substack{g}
:\{1,\dots,m\}\to\{1,\dots,n\}
}
\right|\substack{f,g:\;\text{weakly increasing,}\\
\forall i\in\{1,\dots,n\},\ \#f^{-1}(i)-\#g^{-1}(i)=\pm1
}
\right\} 
\]
and $c_{i}(f,g):=\sum_{j\in f^{-1}(i)}a_{j}+\sum_{j\in f^{-1}(i)}b_{j}$.
For example, $x_{a}x_{b}\mshuffle x_{c}=x_{c}x_{a}x_{b}+x_{a}x_{c}x_{b}+x_{a}x_{b}x_{c}-x_{a+b+c}$.
Let $\mathcal{H}$ be the ring of motivic periods of mixed Tate motives
over $\mathbb{Q}$. We define a $\mathbb{Q}$-linear map $L_{B}^{\shuffle}:\mathfrak{X}_{{\rm od}}\to\mathcal{H}$
by 
\[
L_{B}^{\shuffle}(x_{l_{1}}\cdots x_{l_{d}}):=I_{\mathrm{bl}}^{\mathfrak{m}}(l_{1},\dots,l_{d}).
\]
Also, we put 
\[
\mathfrak{X}_{{\rm ev}}'\coloneqq\bigoplus_{d=1}^{\infty}\bigoplus_{\substack{(l_{1},\dots,l_{d}):\,{\rm even}\\
(l_{1},\dots,l_{d})\neq(1,\dots,1)
}
}\mathbb{Q}x_{l_{1}}\cdots x_{l_{d}}.
\]

\begin{defn}[Block regularization]
Define a $\mathbb{Q}$-linear map $L_{B}:\mathfrak{X}_{{\rm od}}\to\mathcal{H}$
by the generating series
\[
\sum_{\Bbbk:\,{\rm odd}}L_{B}(x_{\Bbbk})X_{\Bbbk}\coloneqq\Gamma_{1}(X_{1})^{-2}\left(\sum_{\Bbbk:\,{\rm odd}}L_{B}^{\shuffle}(x_{\Bbbk})X_{\Bbbk}\right)\Gamma_{1}(-X_{1})^{-2}\ \ \ \ \Big(\in\mathcal{H}\langle\langle X_{1},X_{2},\dots\rangle\rangle\Big)
\]
where $X_{(k_{1},\dots,k_{d})}=X_{k_{1}}\cdots X_{k_{d}}$ and $\Gamma_{1}(t)=\exp(\sum_{n=2}^{\infty}\zeta^{\mathfrak{m}}(n)(-t)^{n}/n)\in\mathcal{H}[[t]]$.
\end{defn}

Note that $L_{B}(x_{\Bbbk})=L_{B}^{\shuffle}(x_{\Bbbk})$ for any
odd admissible index $\Bbbk$ by definition. The following are the
main results of this paper.
\begin{thm}[Block shuffle identity, the first form]
\label{thm:BlockShuffleIdentity}For $u\in\mathfrak{X}_{{\rm ev}}'$
and $v\in\mathfrak{X}_{{\rm od}}$, we have 
\[
L_{B}^{\shuffle}(u\mshuffle v)=0.
\]
\end{thm}

\begin{thm}[Block shuffle identity, the second form]
\label{thm:BlockShuffleIdentity2}For $u\in\mathfrak{X}_{{\rm ev}}^{+}$
and $v\in\mathfrak{X}_{{\rm od}}$, we have 
\[
L_{B}(u\mshuffle v)=0.
\]
\end{thm}

Using some algebraic identities (Propositions \ref{prop:cyclic_in_msh}
and \ref{prop:odd_alt_sum_in_msh}), our main theorem yields all the
aforementioned conjectures:
\begin{thm}
Conjectures \ref{Conj:Ch1} and \ref{Conj:Ch2} (and thus Conjectures
\ref{conj:BBB}, \ref{Conj:BBBL1}, \ref{Conj:BBBL2}) are true.
\end{thm}

\begin{rem}
Charlton formulated a generalization of Conjecture \ref{Conj:Ch1}
to the non-admissible case, which is also a consequence of Theorem
\ref{thm:BlockShuffleIdentity2} (see Proposition \ref{prop:Charlton_full_ver}).
\end{rem}

\begin{rem}
Let us make a few historical remarks. The case $m_{0}+\cdots+m_{2n}=1$
of Conjecture \ref{Conj:BBBL1} is proved in \cite[Theorem 2]{BBBL}
(``dressed with 2'' identity). Later on, Bowman and Bradley showed
a more general result
\[
\sum_{m_{0}+\cdots+m_{2n}=m}Z(m_{0},\dots,m_{2n})=\frac{1}{2n+1}\binom{m+2n}{m}\cdot\frac{\pi^{\mathrm{wt}}}{(\mathrm{wt}+1)!}
\]
which is also a consequence (obtained by summing all the cases with
$m_{0}+\cdots+m_{2n}=m$) of Conjecture \ref{Conj:BBBL1}. By a method
using motivic multiple zeta values, Charlton \cite{Ch15_Eulerinan}
showed that $\zeta(\{\{2\}^{m},1,\{2\}^{m},3\}^{n},\{2\}^{m})$ are
Eulerian, i.e., $\zeta(\{\{2\}^{m},1,\{2\}^{m},3\}^{n},\{2\}^{m})\in\mathbb{Q}\cdot\pi^{{\rm wt}}$
for all $m,n\geq0$. It should be stressed that Charlton's result
only gives the existence of the rational numbers $q_{n,m}\in\mathbb{Q}$
such that
\[
\zeta(\{\{2\}^{m},1,\{2\}^{m},3\}^{n},\{2\}^{m})=q_{n,m}\cdot\pi^{{\rm wt}},
\]
but does not provide a way to determine $q_{n,m}$ itself. Also, Keilthy
showed in \cite[Corollary 2.4.4]{Kei_thesis}, \cite[Corollary 6.4]{Kei_21}
that 
\[
L_{B}^{\shuffle}(u\shuffle v)\equiv0
\]
holds modulo lower block degree and products. To compare his result
with our results, notice that the top degree ($n=l+m$) terms in the
definition of $x_{a_{1}}\cdots x_{a_{l}}\mshuffle x_{b_{1}}\cdots x_{b_{m}}$
agrees with
\[
x_{a_{1}}\cdots x_{a_{l}}\shuffle x_{b_{1}}\cdots x_{b_{m}}.
\]
Thus, Keilthy's result says that the block shuffle identity (Theorems
\ref{thm:BlockShuffleIdentity}, \ref{thm:BlockShuffleIdentity2})
holds modulo lower block degree and products. As a final remark, Hoffman's
conjectural identities
\[
\zeta(3,\{2\}^{n},1,2)=\zeta(\{2\}^{n+3})+2\zeta(3,3,\{2\}^{n})
\]
are another consequence of Conjecture \ref{Conj:Ch1}, which has been
proved by the authors in their previous paper \cite{HoffmanConjProof}.
\end{rem}

Another interesting consequence of our main theorem is:
\begin{prop}
\label{prop:23_shuffle}Let $W_{2,3}$ denote the finite sequences
of $2$ and $3$, and $\mathbb{W}_{2,3}$ the set of formal $\mathbb{Q}$-linear
sums of elements of $W_{2,3}$. We define a binary operation $\star:W_{2,3}\times W_{2,3}\to\mathbb{W}_{2,3}$
recursively by
\begin{align*}
\Bbbk\star\emptyset=\emptyset\star\Bbbk & =\Bbbk\\
(\Bbbk,2)\star\Bbbk'=\Bbbk\star(\Bbbk',2) & =(\Bbbk\star\Bbbk',2)\\
(\Bbbk,3)\star(\Bbbk',3) & =(\Bbbk\star(\Bbbk',3),3)+((\Bbbk,3)\star\Bbbk',3)+(\Bbbk\star\Bbbk',2,2,2).
\end{align*}
Then for $\Bbbk,\Bbbk'\in W_{2,3}$, we have
\[
\zeta^{\mathfrak{m}}(\Bbbk',\Bbbk^{\dagger})=\zeta^{\mathfrak{m}}(\Bbbk\star\Bbbk')
\]
where $\Bbbk^{\dagger}$ is the dual index of $\Bbbk$.
\end{prop}

\begin{rem}
\label{rem:Hoffman}Notice that, the left-hand side $\zeta^{\mathfrak{m}}(\Bbbk',\Bbbk^{\dagger})$
of this propostion is of the form
\[
\zeta^{\mathfrak{m}}(\{2\}^{a_{1}},3,\ldots,\{2\}^{a_{r}},3,\{2\}^{c},1,\{2\}^{b_{s}+1},\ldots,1,\{2\}^{b_{1}+1})\quad\left(a_{1},\ldots,a_{r},b_{1},\ldots,b_{s},c\geq0\right)
\]
which is equal to 
\[
\pm I_{{\rm bl}}^{\mathfrak{m}}(2a_{1}+3,\ldots,2a_{r}+3,2c+2,2b_{s}+3,\ldots,2b_{1}+3)
\]
under the block notation. Given an index $(l_{1},\dots,l_{d})$ with
$l_{1},\dots,l_{d}\geq2$ and $\#\{l_{i}:{\rm even}\}=1$, this proposition
gives a way to express $I_{{\rm bl}}^{\mathfrak{m}}(l_{1},\dots,l_{d})$
as a sum of $\zeta^{\mathfrak{m}}(\Bbbk)$ with $\Bbbk\in W_{2,3}$.
As we know that $\{\zeta^{\mathfrak{m}}(\Bbbk)\}_{\Bbbk\in W_{2,3}}$
(namely, the Hoffman basis) are linearly independent by Brown's theorem,
our theorem gives all $\mathbb{Q}$-linear relations among such $I_{{\rm bl}}^{\mathfrak{m}}(l_{1},\dots,l_{d})$'s.
This proposition also says that $\zeta^{\mathfrak{m}}(\Bbbk',\Bbbk^{\dagger})$
is a $\mathbb{Z}$-linear sum of $\zeta^{\mathfrak{m}}(\Bbbk)$ with
$\Bbbk\in W_{2,3}$ (in light of its significance, we restate this
fact as Theorem \ref{thm:characterizing_BS} in the final section
and discuss related experimental observations). At least numerically,
it seems that these are the only multiple zeta values having such
a property, and no other multiple zeta values are $\mathbb{Z}$-linear
sum of the Hoffman basis. Nevertheless, we also conjecture general
$2$-integrality of the coefficients in the expansion by the Hoffman
basis, that is, all the motivic multiple zeta values are $\mathbb{Z}_{(2)}$-linear
sums of $\zeta^{\mathfrak{m}}(\Bbbk)$ for $\Bbbk\in W_{2,3}$.
\end{rem}

The outline of the proof of Theorems \ref{thm:BlockShuffleIdentity}
and \ref{thm:BlockShuffleIdentity2} is as follows. In a nutshell,
the basic strategy is to construct a lift of the block shuffle identity
to hyperlogarithms and make use of their differential equations. To
be a bit more precise, we will use a refined algebraic version of
the differential equations in order to work at the level of motivic
iterated integrals. The proof is roughly divided into three parts;
the construction of the multivariable block shuffle product, the proof
of the algebraic differential formula of the multivariable block shuffle
product, and the derivation of the block shuffle identity from the
algebraic differential formula. Here, for the final step of the proof,
we use the motivicity result of the confluence relation proved in
another article \cite{ConfMotivicityArxiv} of the authors and additional
arguments concerning the non-admissible generalization.

This paper is organized as follows. To begin with, in Section \ref{sec:conseq},
we discuss various consequences of the block shuffle identities described
above. We will show that all the aforementioned conjectures (Conjectures
\ref{Conj:Ch1}, \ref{Conj:Ch2}) together with Proposition \ref{prop:23_shuffle}
are corollaries of the block shuffle identity (Theorem \ref{thm:BlockShuffleIdentity2})
via certain algebraic identities/congruences. The subsequent sections
(Sections \ref{sec:dshuffle}, \ref{sec:diff_formula} and \ref{sec:proof_of_mshuffle})
are mainly devoted to the proof of the block shuffle identity. Firstly,
in Section \ref{sec:dshuffle}, we will construct a multivariable
lift of the block shuffle product; we give a definition of the multivariable
block shuffle product and prove its compatibility with the original
block shuffle product. Compared to the original block shuffle product,
the definition of the multivariable block shuffle product has a much
more complex form. For this reason, we will describe its combinatorial
nature in a more intuitive way and also give several recurrence formulas
satisfied by the multivariable block shuffle product (Propositions
\ref{prop:recur_dshuffle} and \ref{prop:dsh_sep_as_as1}). Then,
in Section \ref{sec:diff_formula}, we will prove a key differential
formula of the multivariable block shuffle product. Finally, in Section
\ref{sec:proof_of_mshuffle}, we will prove the full versions of the
block shuffle identity (Theorems \ref{thm:BlockShuffleIdentity},
\ref{thm:BlockShuffleIdentity2}). The derivation of the full version
(non-admissible extension) of Charlton's generalized cyclic insertion
conjecture and analog of the block shuffle identity for refined symmetric
multiple zeta values (Conjecture \ref{conj:block_shuffle_for_symmetric_MZV},
Theorem \ref{thm:block_shuffle_for_symmetric_MZV}) are also discussed
at the end of Section \ref{sec:proof_of_mshuffle}. Finally, in Section
\ref{sec:remarks}, we will discuss a certain maximality of the block
shuffle relation and their significance from the perspective of the
size of the coefficients.

\section{consequences of the block shuffle identity\label{sec:conseq}}

\subsection{Derivation of Charlton's conjectures from the block shuffle identity}

Let $\mathfrak{X}^{+}\mshuffle\mathfrak{X}^{+}$ be the $\mathbb{Q}$-vector
subspace of $\mathfrak{X}$ spanned by $\{u\mshuffle v\mid u,v\in\mathfrak{X}^{+}\}$.
The purpose of this section is prove some algebraic identities and
derive Conjectures \ref{Conj:Ch1} and \ref{Conj:Ch2} from Theorem
\ref{thm:BlockShuffleIdentity2}. For convenience, we extend the definition
of $L_{B}$ to $\mathfrak{X}$ by setting $L_{B}(\mathfrak{X}_{{\rm ev}})=\{0\}$.
Since Theorem \ref{thm:BlockShuffleIdentity2} says that
\[
\mathfrak{X}^{+}\mshuffle\mathfrak{X}^{+}\subset\ker L_{B},
\]
it is enough to prove the following propositions.
\begin{prop}
\label{prop:cyclic_in_msh}For $d\geq1$ and $l_{1},\dots,l_{d}\in\mathbb{Z}_{\geq1}$,
we have
\[
\sum_{i=0}^{d-1}x_{l_{i+1}}\cdots x_{l_{d}}x_{l_{1}}\cdots x_{l_{i}}\equiv\begin{cases}
x_{l_{1}+\cdots+l_{d}} & d:{\rm odd}\\
0 & d:{\rm even}
\end{cases}\pmod{\mathfrak{X}^{+}\mshuffle\mathfrak{X}^{+}}.
\]
\end{prop}

\begin{prop}
\label{prop:odd_alt_sum_in_msh}For $n\geq1$ and $a_{1},\dots a_{n+1},b_{1},\dots,b_{n}\geq1$,
we have
\[
\sum_{\sigma\in\mathfrak{S}_{n+1}}{\rm sgn}(\sigma)x_{a_{\sigma(1)}}x_{b_{1}}x_{a_{\sigma(2)}}x_{b_{2}}\cdots x_{a_{\sigma(n)}}x_{b_{n}}x_{a_{\sigma(n+1)}}\in\mathfrak{X}^{+}\mshuffle\mathfrak{X}^{+}.
\]
\end{prop}

To prove these propositions, we need some lemmas. For $c\geq1$, we
define $\mathbb{Q}$-bilinear maps $f_{c}:\mathfrak{X}\times\mathfrak{X}\to\mathfrak{X}^{+}$
by
\begin{align*}
f_{c}(x_{a_{1}}\cdots x_{a_{l}},x_{b_{1}}\cdots x_{b_{m}})= & \sum_{n=1}^{l+m+1}(-1)^{(l+m+1-n)/2}\sum_{\substack{(f,g)\in\widetilde{\mathscr{D}}_{n}}
}x_{\alpha_{1}(f,g)}\cdots x_{\alpha_{n}(f,g)}
\end{align*}
where
\[
\widetilde{\mathscr{D}}_{n}\coloneqq\left\{ \left.\substack{\substack{f}
:\{1,\dots,l\}\to\{1,\dots,n\}\\
\substack{g}
:\{1,\dots,m\}\to\{1,\dots,n\}
}
\right|\substack{f,g:\;\text{weakly increasing,}\\
\#f^{-1}(i)-\#g^{-1}(i)=\left\{ \substack{0\\
\pm1
}
\substack{i=1\quad\\
\;\;i\in\{2,\dots,n\}
}
\right.
}
\right\} 
\]
and
\[
\alpha_{i}(f,g)=c\delta_{1,i}+\sum_{j\in f^{-1}(i)}a_{j}+\sum_{j\in f^{-1}(i)}b_{j}.
\]
By definition, $f_{c}$ satisfies 
\begin{equation}
f_{c}(u,1)=f_{c}(1,u)=x_{c}u\quad(u\in\mathfrak{X})\label{eq:f_initial}
\end{equation}
 and
\begin{align}
f_{c}(x_{a}u,x_{b}v)= & x_{c}f_{a}(u,x_{b}v)+x_{c}f_{b}(x_{a}u,v)-f_{c+a+b}(u,v)\quad(u,v\in\mathfrak{X}),\label{eq:f_recurrence}
\end{align}
which also gives an inductive definition of $f_{c}.$ Note also that
\begin{equation}
x_{a}u\mshuffle x_{b}v=f_{a}(u,x_{b}v)+f_{b}(x_{a}u,v)\quad(u,v\in\mathfrak{X})\label{eq:diamond_f}
\end{equation}
by definition.
\begin{lem}
\label{lem:move_an}For $n\geq1$ and $a_{1},\dots,a_{n}\in\mathbb{Z}_{\geq1}$
and $v\in\mathfrak{X}$, we have
\[
\sum_{i=0}^{n-1}S(x_{a_{1}}\cdots x_{a_{i}})\mshuffle x_{a_{i+1}}\cdots x_{a_{n}}v=f_{a_{n}}(S(x_{a_{1}}\cdots x_{a_{n-1}}),v)
\]
where $S$ is the anti-automorphism of $\mathfrak{X}$ defined by
$S(x_{k})=-x_{k}$. In particular,
\[
x_{a_{1}}\cdots x_{a_{n}}v\equiv f_{a_{n}}(S(x_{a_{1}}\cdots x_{a_{n-1}}),v)\pmod{\mathfrak{X}^{+}\mshuffle\mathfrak{X}^{+}}.
\]
\end{lem}

\begin{proof}
Put 
\[
g_{i}=f_{a_{i}}(S(x_{a_{1}}\cdots x_{a_{i-1}}),x_{a_{i+1}}\cdots x_{a_{n}}v)
\]
for $i=1,\dots,n$ and $g_{0}=0$. Then, by (\ref{eq:diamond_f}),
we have
\[
S(x_{a_{1}}\cdots x_{a_{i}})\mshuffle x_{a_{i+1}}\cdots x_{a_{n}}v=g_{i+1}-g_{i}
\]
for $i=0,\dots,n-1$. Thus, by telescoping the sum, we get
\[
\sum_{i=0}^{n-1}S(x_{a_{1}}\cdots x_{a_{i}})\mshuffle x_{a_{i+1}}\cdots x_{a_{n}}v=g_{n}-g_{0}=f_{a_{n}}(S(x_{a_{1}}\cdots x_{a_{n-1}}),v)
\]
which proves the claim.
\end{proof}
\begin{lem}
\label{lem:flem_for_cyclic}For $n\geq0$ and $a_{1},\dots,a_{n}\geq1$,
we have
\[
\sum_{i=0}^{n}f_{c}(S(x_{a_{i+1}}\cdots x_{a_{n}}),x_{a_{1}}\cdots x_{a_{i}})=\begin{cases}
x_{c+a_{1}+\cdots+a_{n}} & n:{\rm even}\\
0 & n:{\rm odd}.
\end{cases}
\]
\end{lem}

\begin{proof}
We prove the lemma by the induction on $n$. The cases $n=0,1$ are
trivial by definition. Suppose that $n\geq2$. Then, the left-hand
side can be computed as
\begin{eqnarray*}
 &  & \sum_{i=0}^{n}f_{c}(S(x_{a_{1}}\cdots x_{a_{i}}),x_{a_{i+1}}\cdots x_{a_{n}})\\
 & \stackrel{\text{(\ref{eq:f_initial})}}{=} & \overbrace{x_{c}x_{a_{1}}\cdots x_{a_{n}}}^{(i=0\text{ term})}+\overbrace{x_{c}S(x_{a_{1}}\cdots x_{a_{n}})}^{(i=n\text{ term})}+\sum_{i=1}^{n-1}f_{c}(S(x_{a_{i+1}}\cdots x_{a_{n}}),x_{a_{1}}\cdots x_{a_{i}})\\
 & \stackrel{\text{(\ref{eq:f_recurrence})}}{=} & x_{c}x_{a_{1}}\cdots x_{a_{n}}+x_{c}S(x_{a_{1}}\cdots x_{a_{n}})+\sum_{i=1}^{n-1}f_{c+a_{1}+a_{n}}(S(x_{a_{i+1}}\cdots x_{a_{n-1}}),x_{a_{2}}\cdots x_{a_{i}})\\
 &  & -\sum_{i=1}^{n-1}x_{c}f_{a_{n}}(S(x_{a_{i+1}}\cdots x_{a_{n-1}}),x_{a_{1}}\cdots x_{a_{i}})+\sum_{i=1}^{n-1}x_{c}f_{a_{1}}(S(x_{a_{i+1}}\cdots x_{a_{n}}),x_{a_{2}}\cdots x_{a_{i}})\\
 & \stackrel{\text{(\ref{eq:f_initial})}}{=} & \sum_{i=1}^{n-1}f_{c+a_{1}+a_{n}}(S(x_{a_{i+1}}\cdots x_{a_{n-1}}),x_{a_{2}}\cdots x_{a_{i}})\\
 &  & -\sum_{i=0}^{n-1}x_{c}f_{a_{n}}(S(x_{a_{i+1}}\cdots x_{a_{n-1}}),x_{a_{1}}\cdots x_{a_{i}})+\sum_{i=1}^{n}x_{c}f_{a_{1}}(S(x_{a_{i+1}}\cdots x_{a_{n}}),x_{a_{2}}\cdots x_{a_{i}})\\
 & \stackrel{\text{ind.\ hyp.}}{=} & \begin{cases}
x_{c+a_{1}+a_{n}+(a_{2}+\cdots+a_{n-1})} & n:{\rm even}\\
0 & n:{\rm odd},
\end{cases}
\end{eqnarray*}
by which we complete the proof.
\end{proof}
\begin{lem}
\label{lem:flem_for_alt}For $n\geq1$ and $a_{1},\dots a_{n+1},b_{1},\dots,b_{n}\geq1$,
we have
\[
\sum_{\sigma\in\mathfrak{S}_{n+1}}\mathrm{sgn}(\sigma)f_{b_{1}}(x_{a_{\sigma(1)}},\,x_{a_{\sigma(2)}}x_{b_{2}}\cdots x_{a_{\sigma(n)}}x_{b_{n}}x_{a_{\sigma(n+1)}})=0.
\]
\end{lem}

\begin{proof}
We prove the claim by induction on $n$. The case $n=1$ is trivial
by definition, so let $n\geq2$. For $\sigma\in\mathfrak{S}_{n+1}$,
put
\begin{align*}
u_{\sigma}= & x_{a_{\sigma(3)}}x_{b_{3}}\cdots x_{a_{\sigma(n)}}x_{b_{n}}x_{a_{\sigma(n+1)}}.
\end{align*}
Then the left hand side of lemma can be written as 
\[
\sum_{\sigma\in\mathfrak{S}_{n+1}}\mathrm{sgn}(\sigma)f_{b_{1}}(x_{a_{\sigma(1)}},\,x_{a_{\sigma(2)}}x_{b_{2}}u_{\sigma}).
\]
By applying (\ref{eq:f_recurrence}) twice together with (\ref{eq:f_initial}),
we get
\begin{align}
f_{b_{1}}(x_{a_{\sigma(1)}},\,x_{a_{\sigma(2)}}x_{b_{2}}u_{\sigma})= & -x_{b_{1}+a_{\sigma(1)}+a_{\sigma(2)}}x_{b_{2}}u_{\sigma}-x_{b_{1}}x_{a_{\sigma(1)}+a_{\sigma(2)}+b_{2}}u_{\sigma}\label{eq:twice_expanded}\\
 & +x_{b_{1}}x_{a_{\sigma(1)}}x_{a_{\sigma(2)}}x_{b_{2}}u_{\sigma}+x_{b_{1}}x_{a_{\sigma(2)}}x_{a_{\sigma(1)}}x_{b_{2}}u_{\sigma}\nonumber \\
 & +x_{b_{1}}x_{a_{\sigma(2)}}f_{b_{2}}(x_{a_{\sigma(1)}},u_{\sigma}).\nonumber 
\end{align}
Noting the trivial identities
\begin{align*}
\sum_{\sigma\in\mathfrak{S}_{n+1}}\mathrm{sgn}(\sigma)x_{b_{1}+a_{\sigma(1)}+a_{\sigma(2)}}x_{b_{2}}u_{\sigma}= & 0,\\
\sum_{\sigma\in\mathfrak{S}_{n+1}}\mathrm{sgn}(\sigma)x_{b_{1}}x_{a_{\sigma(1)}+a_{\sigma(2)}+b_{2}}u_{\sigma}= & 0,\\
\sum_{\sigma\in\mathfrak{S}_{n+1}}\mathrm{sgn}(\sigma)x_{b_{1}}\left(x_{a_{\sigma(1)}}x_{a_{\sigma(2)}}+x_{a_{\sigma(2)}}x_{a_{\sigma(1)}}\right)x_{b_{2}}u_{\sigma}= & 0,
\end{align*}
the alternating sum of (\ref{eq:twice_expanded}) is computed as
\[
\sum_{\sigma\in\mathfrak{S}_{n+1}}\mathrm{sgn}(\sigma)f_{b_{1}}(x_{a_{\sigma(1)}},\,x_{a_{\sigma(2)}}x_{b_{2}}u_{\sigma})=\sum_{\sigma\in\mathfrak{S}_{n+1}}\mathrm{sgn}(\sigma)x_{b_{1}}x_{a_{\sigma(2)}}f_{b_{2}}(x_{a_{\sigma(1)}},u_{\sigma})=0,
\]
where for the last equality we used the induction hypothesis. This
proves the lemma.
\end{proof}
Now, let us prove Propositions \ref{prop:cyclic_in_msh} and \ref{prop:odd_alt_sum_in_msh}.

\begin{proof}[Proof of Proposition \ref{prop:cyclic_in_msh}]
Let $d\geq1$ and $l_{1},\dots,l_{d}\in\mathbb{Z}_{\geq1}$. By Lemma
\ref{lem:move_an}, we have
\begin{align*}
x_{l_{i+1}}\cdots x_{l_{d}}x_{l_{1}}\cdots x_{l_{i}}\equiv & f_{l_{d}}(S(x_{l_{i+1}}\cdots x_{l_{d-1}}),x_{l_{1}}\cdots x_{l_{i}})
\end{align*}
modulo $\mathfrak{X}^{+}\mshuffle\mathfrak{X}^{+}$ for $i=0,\dots,d-1$.
Summing this for $0\leq i\leq d-1$, and using Lemma \ref{lem:flem_for_cyclic},
we have
\[
\sum_{i=0}^{d-1}x_{l_{i+1}}\cdots x_{l_{d}}x_{l_{1}}\cdots x_{l_{i}}\equiv\sum_{i=0}^{d-1}f_{l_{d}}(S(x_{l_{i+1}}\cdots x_{l_{d-1}}),x_{l_{1}}\cdots x_{l_{i}})=\begin{cases}
x_{l_{1}+\cdots+l_{d}} & d:{\rm odd}\\
0 & d:{\rm even},
\end{cases}
\]
which proves the claim.
\end{proof}

\begin{proof}[Proof of Proposition \ref{prop:odd_alt_sum_in_msh}]
Let $n\geq1$ and $a_{1},\dots a_{n+1},b_{1},\dots,b_{n}\geq1$.
Again, by Lemma \ref{lem:move_an}, we have
\[
x_{a_{1}}x_{b_{1}}x_{a_{2}}x_{b_{2}}\cdots x_{a_{n}}x_{b_{n}}x_{a_{n+1}}\equiv-f_{b_{1}}(x_{a_{1}},\,x_{a_{2}}x_{b_{2}}\cdots x_{a_{n}}x_{b_{n}}x_{a_{n+1}}).
\]
modulo $\mathfrak{X}^{+}\mshuffle\mathfrak{X}^{+}$. Taking the alternating
sum over all permutations of $a_{i}$'s, and using Lemma \ref{lem:flem_for_alt},
\begin{align*}
\sum_{\sigma\in\mathfrak{S}_{n+1}}{\rm sgn}(\sigma)x_{a_{\sigma(1)}}x_{b_{1}}x_{a_{\sigma(2)}}x_{b_{2}}\cdots x_{a_{\sigma(n)}}x_{b_{n}}x_{a_{\sigma(n+1)}} & \equiv-\sum_{\sigma\in\mathfrak{S}_{n+1}}\mathrm{sgn}(\sigma)f_{b_{1}}(x_{a_{\sigma(1)}},\,x_{a_{\sigma(2)}}x_{b_{2}}\cdots x_{a_{\sigma(n)}}x_{b_{n}}x_{a_{\sigma(n+1)}})\\
 & =0,
\end{align*}
which proves the claim.
\end{proof}

\subsection{Derivation of Proposition \ref{prop:23_shuffle} from the block shuffle
identity}

In this section, we derive Proposition \ref{prop:23_shuffle} from
Theorem \ref{thm:BlockShuffleIdentity2}. For this purpose, we first
reformulate Proposition \ref{prop:23_shuffle} in terms of $L_{B}$.
Let $W_{2,3}'$ be a subset of $W_{2,3}$ consisting of indices which
do not end with $2$. Define 
\[
\iota:W_{2,3}'\to\mathfrak{X},\ \theta:W_{2,3}\to\mathfrak{X}
\]
by
\begin{align*}
\iota(\{2\}^{k_{d}},3,\dots\{2\}^{k_{1}},3) & =(-1)^{k_{1}+\cdots+k_{d}}x_{2k_{1}+3}\cdots x_{2k_{d}+3},\\
\theta(\{2\}^{k_{d}},3,\dots\{2\}^{k_{1}},3,\{2\}^{k_{0}}) & =(-1)^{k_{0}+\cdots+k_{d}}x_{2k_{0}+2}x_{2k_{1}+3}\cdots x_{2k_{d}+3}.
\end{align*}
By definition, they satisfy the recursion
\begin{align*}
\iota(\Bbbk,\{2\}^{c},3) & =(-1)^{c}x_{2c+3}\cdot\iota(\Bbbk)\quad(\Bbbk\in W_{2,3}'),\\
\theta(\Bbbk,3,\{2\}^{c}) & =(-1)^{c}x_{2c+2}\cdot s_{1}(\theta(\Bbbk))\quad(\Bbbk\in W_{2,3}).
\end{align*}
By definition,
\begin{align*}
 & \zeta^{\mathfrak{m}}(\{2\}^{k_{d}},3,\dots\{2\}^{k_{1}},3,\{2\}^{k_{0}})\\
 & =(-1)^{k_{0}+\cdots+k_{d}+d}I^{\mathfrak{m}}(\overbrace{0,1,\dots,0,1,0}^{2k_{d}+3},\overbrace{0,1,\dots,0,1,0}^{2k_{d-1}+3},\dots,\overbrace{0,1,\dots,0,1,0}^{2k_{1}+3},\overbrace{0,1,\dots,0,1}^{2k_{0}+2})\\
 & =(-1)^{k_{0}+\cdots+k_{d}}L_{B}(x_{2k_{0}+2}x_{2k_{1}+3}\cdots x_{2k_{d}+3}),
\end{align*}
and thus
\[
\zeta^{\mathfrak{m}}(\Bbbk)=L_{B}(\theta(\Bbbk))\ \ \ (\Bbbk\in W_{2,3}).
\]
Now write $\Bbbk,\Bbbk'\in W_{2,3}$ in the forms
\begin{align*}
\Bbbk & =(\bm{k},\{2\}^{a})\quad(\bm{k}\in W_{2,3}'),\\
\Bbbk' & =(\bm{l},\{2\}^{b})\quad(\bm{l}\in W_{2,3}').
\end{align*}
Then $\zeta^{\mathfrak{m}}(\Bbbk\star\Bbbk')$ (the right-hand side
of Proposition \ref{prop:23_shuffle}) is equal to 
\[
L_{B}(\theta(\bm{k}\star\bm{l},\{2\}^{a+b}))
\]
as checked above, and by a similar calculation $\zeta^{\mathfrak{m}}(\Bbbk',\Bbbk^{\dagger})$
(the left-hand side of Proposition \ref{prop:23_shuffle}) is equal
to 
\[
(-1)^{a+b}L_{B}(S(\iota(\bm{l}))x_{2+2a+2b}\iota(\bm{k})).
\]
Hence, Proposition \ref{prop:23_shuffle} is equivalent to
\[
S(\iota(\bm{l}))x_{2+2c}\iota(\bm{k})-(-1)^{c}\theta(\bm{k}\star\bm{l},\{2\}^{c})\in\ker L_{B}
\]
with $c=a+b$.

Now let us derive Proposition \ref{prop:23_shuffle} from Theorem
\ref{thm:BlockShuffleIdentity2}. For this purpose, it is sufficient
to show
\[
S(\iota(\bm{l}))x_{2+2c}\iota(\bm{k})-(-1)^{c}\theta(\bm{k}\star\bm{l},\{2\}^{c})\in\mathfrak{X}^{+}\mshuffle\mathfrak{X}^{+}.
\]
Noting that the first term is congruent to
\[
f_{2+2c}(\iota(\bm{l}),\iota(\bm{k})),
\]
modulo $\mathfrak{X}^{+}\mshuffle\mathfrak{X}^{+}$ by Lemma \ref{lem:move_an},
it is also enough to prove
\[
f_{2+2c}(\iota(\bm{l}),\iota(\bm{k}))\equiv(-1)^{c}\theta(\bm{k}\star\bm{l},\{2\}^{c})\pmod{\mathfrak{X}^{+}\mshuffle\mathfrak{X}^{+}}.
\]
In fact, the following stronger equality holds:
\begin{lem}
For $c\geq0$ and $\bm{k},\bm{l}\in W_{2,3}'$, we have
\[
f_{2+2c}(\iota(\bm{l}),\iota(\bm{k}))=(-1)^{c}\theta(\bm{k}\star\bm{l},\{2\}^{c}).
\]
\end{lem}

\begin{proof}
We prove the claim by induction on $\bm{k}$ and $\bm{l}$. If $\bm{k}=\emptyset$
(or $\bm{l}=\emptyset$), then the equality holds since
\begin{align*}
f_{2+2c}(\iota(\bm{l}),\iota(\emptyset)) & =x_{2+2c}\iota(\bm{l})\\
 & =(-1)^{c}\theta(\bm{l},\{2\}^{c}).
\end{align*}
Now, let $\bm{k}=(\bm{k}',\{2\}^{s},3)$ and $\bm{l}=(\bm{l}',\{2\}^{t},3)$.
Then,
\begin{align*}
f_{2+2c}(\iota(\bm{l}),\iota(\bm{k})) & =(-1)^{s+t}f_{2+2c}(x_{2t+3}\iota(\bm{l}'),x_{2s+3}\iota(\bm{k}'))\\
 & =(-1)^{s}x_{2+2c}f_{2s+3}(\iota(\bm{l}),\iota(\bm{k}'))+(-1)^{t}x_{2+2c}f_{2t+3}(\iota(\bm{l}'),\iota(\bm{k}))\\
 & \ -(-1)^{s+t}f_{2+2c+2s+2t+6}(\iota(\bm{l}'),\iota(\bm{k}'))\\
 & =(-1)^{s}x_{2+2c}s_{1}(f_{2s+2}(\iota(\bm{l}),\iota(\bm{k}')))+(-1)^{t}x_{2+2c}s_{1}(f_{2t+2}(\iota(\bm{l}'),\iota(\bm{k})))\\
 & \ +(-1)^{s+t+3}f_{2+2c+2s+2t+6}(\iota(\bm{l}'),\iota(\bm{k}')).
\end{align*}
Using induction hypothesis, the last expression is equal to
\begin{align*}
 & x_{2+2c}s_{1}(\theta(\bm{k}'\star\bm{l},\{2\}^{s}))+x_{2+2c}s_{1}(\theta(\bm{k}\star\bm{l}',\{2\}^{t}))+(-1)^{c}\theta(\bm{k}'\star\bm{l}',\{2\}^{c+s+t+3})\\
 & =(-1)^{c}\left(\theta(\bm{k}'\star\bm{l},\{2\}^{s},3,\{2\}^{c})+\theta(\bm{k}\star\bm{l}',\{2\}^{t},3,\{2\}^{c})+\theta(\bm{k}'\star\bm{l}',\{2\}^{s+t+3+c})\right)\\
 & =(-1)^{c}\theta(\bm{k}\star\bm{l},\{2\}^{c}),
\end{align*}
which completes the proof.
\end{proof}
\begin{proof}[Proof of Theorem \ref{thm:characterizing_BS}]
It suffices to show
\[
\ker L_{B}\cap\text{\ensuremath{\mathfrak{X}}}'\subset\mathfrak{X}^{+}\mshuffle\mathfrak{X}^{+}.
\]
Then, by Brown's theorem \cite{Brown12}, $\ker L_{B}\cap\theta(W_{2,3})=0$.
Notice that $\left\{ S(\iota(\bm{l}))x_{2c+2}\iota(\bm{k})\bigr|\bm{l},\bm{k}\in W_{2,3}',c\geq0\right\} $
forms a basis of $\mathfrak{X}'$ and 
\[
\xi_{\bm{k},\bm{l},c}\coloneqq S(\iota(\bm{l}))x_{2c+2}\iota(\bm{k})-(-1)^{c}\theta(\bm{k}\star\bm{l},\{2\}^{c})\in\ker L_{B}
\]
for $\bm{l},\bm{k}\in W_{2,3}'$ and $c\geq0$ by the proof above.
Thus $\left\{ \xi_{\bm{k},\bm{l},c}\bigr|\bm{l},\bm{k}\in W_{2,3}',c\geq0\right\} $
spans $\ker L_{B}\cap\mathfrak{X}'$. Since $\xi_{\bm{k},\bm{l},c}\in\mathfrak{X}^{+}\mshuffle\mathfrak{X}^{+}$,
this proves the desired claim.
\end{proof}

\section{Multivariable generalization of the block shuffle product\label{sec:dshuffle}}

Recall that the block shuffle product is a product defined on $\mathfrak{X}$.
Charlton's alternating block notation (see Definition \ref{def:Block_notation}
for the precise definition) gives a natural way to identify the non-constant
part of $\mathfrak{X}$ as a subspace of $\mathbb{Q}\left\langle e_{0},e_{1}\right\rangle $.
The goal of this section is to extend the block shuffle product to
a multi-variable setting. To give a little more detail, we let $\mathcal{A}\coloneqq\mathbb{Q}\left\langle \left.e_{x}\right|x\in\mathbb{C}\right\rangle $
and $\mathcal{A}^{+}\coloneqq\bigoplus_{\substack{n>0\\
x_{1},\ldots,x_{n}\in\mathbb{C}
}
}\mathbb{Q}\,e_{x_{1}}\cdots e_{x_{n}}$ be its subspace without the constant term. We shall construct a $\mathbb{Q}$-bilinear
map
\[
\dshuffle:\mathfrak{X}\times\mathcal{A}^{+}\to\mathcal{A}^{+}
\]
that extends the block shuffle product
\[
\mshuffle:\mathfrak{X}\times\mathfrak{X}^{+}\to\mathfrak{X}^{+}.
\]
This extension will play a crucial role of our proof of the block
shuffle identity in later section, and the construction of this multivariable
lift is indeed the most non-trivial part throughout the entire proof
of the block shuffle identity.

\subsection{\label{subsec:Def_of_diaprod}Definition of the multivariable block
shuffle product}

In this section, we define a multivariable extension of the block
shuffle product. To begin with, let $\iota$ be an involution on $\mathbb{C}$
defined by $\iota(z)=1-z$ and $\bm{\iota}$ its induced involution
on $\mathcal{A}$ i.e., $\bm{\iota}(e_{z})=e_{1-z}$. We write $\iota^{k}$
(resp. $\bm{\iota}^{k}$) for the $k$-times composition of $\iota$
(resp. $\bm{\iota}^{k}$). For $b\in\mathbb{Z}_{\geq0}$, we put
\[
U(b):=\overbrace{e_{0}e_{1}e_{0}e_{1}\cdots e_{\iota^{b}(1)}}^{b}+\overbrace{e_{1}e_{0}e_{1}e_{0}\cdots e_{\iota^{b}(0)}}^{b}\in\mathcal{A}.
\]

\begin{defn}[Multi-variable block shuffle product]
We define a $\mathbb{Q}$-bilinear map $\dshuffle:\mathfrak{X}\times\mathcal{A}^{+}\to\mathcal{A}^{+}$
by
\[
x_{b_{1}}\cdots x_{b_{m}}\dshuffle e_{a_{1}}\cdots e_{a_{n}}=\sum_{d=n}^{n+m}(-1)^{n+m-d}\sum_{(i_{1},\ldots,i_{d},f)\in\mathscr{B}_{d}}e_{a_{i_{1}}}V_{1}\bm{\iota}^{s_{1}}(e_{a_{i_{2}}})\cdots V_{d-1}\bm{\iota}^{s_{d-1}}(e_{a_{i_{d}}})
\]
where
\[
\mathscr{B}_{d}\coloneqq\left\{ \left.\substack{1=i_{1}\leq\cdots\leq i_{d}=n\\
f:\,\{1,\dots,m\}\to\{1,\dots,d-1\}
}
\right|\substack{f:\text{ weakly increasing}\\
\delta_{\emptyset,f^{-1}(\{k\})}\leq i_{k+1}-i_{k}\leq1.
}
\right\} ,
\]
$s_{1},\dots,s_{d-1}\in\mathbb{Z}$ and $V_{1},\dots,V_{d-1}\in\mathcal{A}$
(both of which depend on $(i_{1},\ldots,i_{d},f)\in\mathscr{B}_{d}$)
are given by 
\[
s_{k}=\sum_{f(j)\leq k}(b_{j}-1)
\]
and
\[
V_{k}=\begin{cases}
1 & f^{-1}(\{k\})=\emptyset\\
U(i_{k+1}-i_{k}-1+\sum_{f(j)=k}b_{j}) & f^{-1}(\{k\})\neq\emptyset.
\end{cases}
\]
\end{defn}

\begin{example}
For example, $x_{k}x_{l}\dshuffle e_{a}e_{b}$ is given by
\begin{align*}
 & +e_{a}U(k+l)e_{b} &  & (d=2,(i_{1},i_{2})=(1,2),f(1)=1,f(2)=1)\\
 & -e_{a}U(k-1)e_{a}U(l)e_{b} &  & (d=3,(i_{1},i_{2},i_{3})=(1,1,2),f(1)=1,f(2)=2)\\
 & -e_{a}U(k+l-1)e_{a}e_{b} &  & (d=3,(i_{1},i_{2},i_{3})=(1,1,2),f(1)=1,f(2)=1)\\
 & -e_{a}U(k)e_{b}U(l-1)e_{b} &  & (d=3,(i_{1},i_{2},i_{3})=(1,2,2),f(1)=1,f(2)=2)\\
 & -e_{a}e_{b}U(k+l-1)e_{b} &  & (d=3,(i_{1},i_{2},i_{3})=(1,2,2),f(1)=2,f(2)=2)\\
 & +e_{a}U(k-1)e_{a}U(l-1)e_{a}e_{b} &  & (d=4,(i_{1},i_{2},i_{3},i_{4})=(1,1,1,2),f(1)=1,f(2)=2)\\
 & +e_{a}U(k-1)e_{a}e_{b}U(l-1)e_{b} &  & (d=4,(i_{1},i_{2},i_{3},i_{4})=(1,1,2,2),f(1)=1,f(2)=3)\\
 & +e_{a}e_{b}U(k-1)e_{b}U(l-1)e_{b} &  & (d=4,(i_{1},i_{2},i_{3},i_{4})=(1,2,2,2),f(1)=2,f(2)=3).
\end{align*}
\end{example}

\begin{rem}
We conjecture that $\dshuffle$ is associative i.e., $u\dshuffle(v\dshuffle w)=(u\mshuffle v)\dshuffle w$
for $u,v\in\mathfrak{X}$ and $w\in\mathcal{A}^{+}$.
\end{rem}

\begin{rem}
To have a better idea of what terms appear in the multivariable block
shuffle product, we give another more intuitive description of $x_{b_{1}}\cdots x_{b_{m}}\dshuffle e_{a_{1}}\cdots e_{a_{n}}$
here. Let $S_{n,m}$ be the set of words in letters $A_{1},\ldots,A_{n}$
and $B_{1},\ldots,B_{m}$ satisfying the following conditions.
\begin{itemize}
\item $A_{1},\ldots,A_{n}$ appear in this order and each $A_{i}$ appear
at least once.
\item $B_{1},\ldots,B_{m}$ appear in this order and each $B_{i}$ appear
exactly once.
\item For each $i$, $A_{i}$ does not appear consecutively.
\item The first letter isomorphism $A_{1}$ while the last letter is $A_{n}$.
\end{itemize}
For example, 
\[
S_{2,2}=\left\{ \begin{array}{c}
A_{1}B_{1}B_{2}A_{2},\,A_{1}B_{1}A_{1}B_{2}A_{2},\,A_{1}B_{1}B_{2}A_{1}A_{2},\,A_{1}B_{1}A_{2}B_{2}A_{2},\,A_{1}A_{2}B_{1}B_{2}A_{2},\\
A_{1}B_{1}A_{1}B_{2}A_{1}A_{2},\,A_{1}B_{1}A_{1}A_{2}B_{2}A_{2},\,A_{1}A_{2}B_{1}A_{2}B_{2}A_{2}
\end{array}\right\} .
\]
(In general, the cardinality of $S_{n,m}$ is $2^{m-1}\left\{ \binom{n+m-1}{m}+\binom{n+m-2}{m}\right\} $).
Then the previous definition is equivalent to
\[
x_{b_{1}}\cdots x_{b_{m}}\dshuffle e_{a_{1}}\cdots e_{a_{n}}\coloneqq\sum_{W\in S_{n,m}}h(W)
\]
where $h$ is a map defined recursively by
\begin{align*}
h(A_{p}B_{i}\cdots B_{j}A_{p}W) & =(-1)^{j-i}e_{a_{p}}U(\sum_{k=i}^{j}b_{k}-1)\cdot\bm{\iota}^{\sum_{k=i}^{j}(b_{k}-1)}\left(h(A_{p}W)\right),\\
h(A_{p}B_{i}\cdots B_{j}A_{p+1}W) & =(-1)^{j-i+1}e_{a_{p}}U(\sum_{k=i}^{j}b_{k})\cdot\bm{\iota}^{\sum_{k=i}^{j}(b_{k}-1)}\left(h(A_{p+1}W)\right),\\
h(A_{p}A_{p+1}W) & =e_{a_{p}}\left(h(A_{p+1}W)\right)
\end{align*}
and $h(A_{n})=e_{a_{n}}$.
\end{rem}

By definition of the multivariable block shuffle product, it satisfies
the following recurrence relation.
\begin{prop}
\label{prop:recur_dshuffle}For $m\geq0,n\geq1$ with $(m,n)\neq(0,1)$,
the above-defined $\dshuffle$ satisfies the recursion formulas
\begin{align*}
x_{b_{1}}\cdots x_{b_{m}}\dshuffle e_{a_{1}}\cdots e_{a_{n}} & =\sum_{k=1}^{m}(-1)^{k-1}e_{a_{1}}U(b_{1}+\cdots+b_{k}-1)\cdot\bm{\iota}^{b_{1}+\cdots+b_{k}-k}(x_{b_{k+1}}\cdots x_{b_{m}}\dshuffle e_{a_{1}}\cdots e_{a_{n}})\\
 & \quad+\sum_{k=1}^{m}(-1)^{k}e_{a_{1}}U(b_{1}+\cdots+b_{k})\cdot\bm{\iota}^{b_{1}+\cdots+b_{k}-k}(x_{b_{k+1}}\cdots x_{b_{m}}\dshuffle e_{a_{2}}\cdots e_{a_{n}})\\
 & \quad+e_{a_{1}}(x_{b_{1}}\cdots x_{b_{m}}\dshuffle e_{a_{2}}\cdots e_{a_{n}})
\end{align*}
and
\begin{align*}
x_{b_{1}}\cdots x_{b_{m}}\dshuffle e_{a_{1}}\cdots e_{a_{n}} & =\sum_{l=0}^{m-1}(-1)^{m-1-l}(x_{b_{1}}\cdots x_{b_{l}}\dshuffle e_{a_{1}}\cdots e_{a_{n}})U(b_{l+1}+\cdots+b_{m}-1)\bm{\iota}^{b_{1}+\cdots+b_{m}-m}(e_{a_{n}})\\
 & \quad+\sum_{l=0}^{m-1}(-1)^{m-l}(x_{b_{1}}\cdots x_{b_{l}}\dshuffle e_{a_{1}}\cdots e_{a_{n-1}})U(b_{l+1}+\cdots+b_{m})\bm{\iota}^{b_{1}+\cdots+b_{m}-m}(e_{a_{n}})\\
 & \quad+(x_{b_{1}}\cdots x_{b_{m}}\dshuffle e_{a_{1}}\cdots e_{a_{n-1}})\bm{\iota}^{b_{1}+\cdots+b_{m}-m}(e_{a_{n}}),
\end{align*}
where we understand that the terms of the form $u\dshuffle1$ with
$u\in\mathfrak{X}$ on the right-hand side are zero.
\end{prop}

The following identity is also useful.
\begin{prop}
\label{prop:dsh_sep_as_as1}Let $m\geq0$ and $n\geq2$. For $1\leq s\leq n-1$,
we have
\begin{align*}
x_{b_{1}}\cdots x_{b_{m}}\dshuffle e_{a_{1}}\cdots e_{a_{n}} & =\sum_{k=0}^{m}(x_{b_{1}}\cdots x_{b_{k}}\dshuffle e_{a_{1}}\cdots e_{a_{s}})\cdot\bm{\iota}^{b_{1}+\cdots+b_{k}-k}(x_{b_{k+1}}\cdots x_{b_{m}}\dshuffle e_{a_{s+1}}\cdots e_{a_{n}})\\
 & \quad+\sum_{0\leq k<l\leq m}(-1)^{l-k}(x_{b_{1}}\cdots x_{b_{k}}\dshuffle e_{a_{1}}\cdots e_{a_{s}})\cdot U(b_{k+1}+\cdots+b_{l})\\
 & \quad\quad\quad\quad\quad\cdot\bm{\iota}^{b_{1}+\cdots+b_{l}-l}(x_{b_{l+1}}\cdots x_{b_{m}}\dshuffle e_{a_{s+1}}\cdots e_{a_{n}}).
\end{align*}
\end{prop}

\begin{proof}
We can verify the claim by separating each term of $x_{b_{1}}\cdots x_{b_{m}}\dshuffle e_{a_{1}}\cdots e_{a_{n}}$
at the right of the last (= the rightmost) $e_{a_{s}}$ (or $e_{\iota(a_{s})}$)
and at the left of the first (= the leftmost) $e_{a_{s+1}}$ (or $e_{\iota(a_{s+1})}$).
\end{proof}

\subsection{\label{subsec:Compatibility_diamonds}Compatibility with the usual
block shuffle product}
\begin{defn}
\label{def:Block_notation}We define a map $B:\mathfrak{X}^{+}\to\mathbb{Q}\left\langle e_{0},e_{1}\right\rangle $
by
\[
B(x_{b_{1}}\cdots x_{b_{m}})=e_{a_{1}}\cdots e_{a_{b_{1}+\cdots+b_{m}}}
\]
where $a_{i}$'s are defined by $a_{1}=0$ and 
\[
a_{i+1}=\begin{cases}
a_{i} & i\in\{b_{1},b_{1}+b_{2},\ldots,b_{1}+\cdots+b_{m-1}\}\\
1-a_{i} & \mathrm{otherwise}.
\end{cases}
\]
\end{defn}

\begin{thm}
\label{thm:block_diamond_compatibility}For $u,v\in\mathfrak{X}$,
we have
\[
u\dshuffle B(v)=B(u\mshuffle v).
\]
\end{thm}

\begin{proof}
By linearity, it is enough to consider the case when $u$ and $v$
are monomials. Let
\[
u=x_{a_{1}}\cdots x_{a_{m}},\ v=x_{b_{1}}\cdots x_{b_{n}}.
\]
We shall prove the claim by induction on $n+m$. Put $d_{k}=a_{1}+\cdots+a_{k}$
and $u_{k}=x_{a_{k+1}}x_{a_{k+2}}\cdots x_{a_{m}}$ for $k=0,\dots,m$.
We write $B(v)$ as
\[
e_{0}\bm{\iota}^{t}(B(v'))
\]
with $t\in\{0,1\}$ and $v'\in\mathfrak{X}$ . More explicitly, 
\[
(t,v')=\begin{cases}
(1,x_{b_{1}-1}x_{b_{2}}\cdots x_{b_{n}}) & b_{1}>1\\
(0,x_{b_{2}}\cdots x_{b_{n}}) & b_{1}=1.
\end{cases}
\]
By the first formula of Proposition \ref{prop:recur_dshuffle}, we
have
\begin{align*}
u\dshuffle B(v) & =\sum_{k=1}^{m}(-1)^{k-1}e_{0}U(d_{k}-1)\bm{\iota}^{d_{k}-k}\left(u_{k}\dshuffle B(v)\right)\\
 & \ \ +\sum_{k=1}^{m}(-1)^{k}e_{0}U(d_{k})\bm{\iota}^{d_{k}-k+t}\left(u_{k}\dshuffle B(v')\right)\\
 & \ \ +e_{0}\bm{\iota}^{t}(u\dshuffle B(v')).
\end{align*}
By the induction hypothesis, the first term $\sum_{k=1}^{m}(-1)^{k-1}e_{0}U(d_{k}-1)\bm{\iota}^{d_{k}-k}\left(u_{k}\dshuffle B(v)\right)$
is equal to
\begin{align*}
 & \sum_{k=1}^{m}(-1)^{k-1}e_{0}U(d_{k}-1)\bm{\iota}^{d_{k}-k}\left(B\left(u_{k}\mshuffle v\right)\right)\\
 & =-\sum_{\substack{k=1\\
k:\,{\rm even}
}
}^{m}B(s_{d_{k}}(u_{k}\mshuffle v)+x_{1}x_{d_{k}-1}(u_{k}\mshuffle v))+\sum_{\substack{k=1\\
k:\,{\rm odd}
}
}^{m}B(x_{d_{k}}(u_{k}\mshuffle v)+x_{1}s_{d_{k}-1}(u_{k}\mshuffle v)),
\end{align*}
by noting $e_{0}U(d_{k}-1)=\prod_{j=0}^{d_{k-1}}e_{\iota^{j}(0)}+e_{0}\prod_{j=0}^{d_{k-2}}e_{\iota^{j}(0)},$
$B(u_{k}\mshuffle v)=e_{0}\cdots$ and $\bm{\iota}^{d_{k}-k}(e_{0}\cdots)=e_{\iota^{d_{k}-k}(0)}\cdots$.
Similarly, the second term $\sum_{k=1}^{m}(-1)^{k}e_{0}U(d_{k})\bm{\iota}^{d_{k}-k+t}\left(u_{k}\dshuffle B(v')\right)$
is rewritten as
\begin{align*}
 & \sum_{k=1}^{m}(-1)^{k}e_{0}U(d_{k})\bm{\iota}^{d_{k}-k+t}\left(B(u_{k}\mshuffle v')\right)\\
 & =\sum_{\substack{k=1\\
k+t:\,{\rm odd}
}
}^{m}(-1)^{k}B(s_{d_{k}+1}(u_{k}\mshuffle v')+x_{1}x_{d_{k}}(u_{k}\mshuffle v'))+\sum_{\substack{k=1\\
k+t:\,{\rm even}
}
}^{m}(-1)^{k}B(x_{d_{k}+1}(u_{k}\mshuffle v')+x_{1}s_{d_{k}}(u_{k}\mshuffle v'))\\
 & =(-1)^{t}\left(-\sum_{\substack{k=1\\
k+t:\,{\rm odd}
}
}^{m}B(s_{d_{k}+1}(u_{k}\mshuffle v')+x_{1}x_{d_{k}}(u_{k}\mshuffle v'))+\sum_{\substack{k=1\\
k+t:\,{\rm even}
}
}^{m}B(x_{d_{k}+1}(u_{k}\mshuffle v')+x_{1}s_{d_{k}}(u_{k}\mshuffle v'))\right),
\end{align*}
and the last term $e_{0}\bm{\iota}^{t}(u\dshuffle B(v'))$ is equal
to
\[
e_{0}\bm{\iota}^{t}(B(u\mshuffle v'))=\begin{cases}
B(x_{1}(u\mshuffle v')) & t=0\\
B(s_{1}(u\mshuffle v')) & t=1.
\end{cases}
\]
Thus if $t=1$ then $u\dshuffle B(v)$ is equal to the image of
\begin{align*}
 & -\left[\sum_{\substack{k=0\\
k:\,{\rm even}
}
}^{m}s_{d_{k}}\left(u_{k}\mshuffle v-s_{1}(u_{k}\mshuffle v')\right)-\sum_{\substack{k=1\\
k:\,{\rm odd}
}
}^{m}\left(x_{d_{k}}(u_{k}\mshuffle v)-x_{d_{k}+1}(u_{k}\mshuffle v')\right)\right]\\
 & +\left[\sum_{\substack{k=1\\
k:\,{\rm odd}
}
}^{m}x_{1}s_{d_{k}-1}\left(u_{k}\mshuffle v-s_{1}(u_{k}\mshuffle v')\right)-\sum_{\substack{k=1\\
k:\,{\rm even}
}
}^{m}\left(x_{1}x_{d_{k}-1}(u_{k}\mshuffle v)-x_{1}x_{d_{k}}(u_{k}\mshuffle v')\right)\right]\\
 & +u\mshuffle v\\
 & =0+0+u\mshuffle v=u\mshuffle v
\end{align*}
under $B$. Here, we used
\[
u_{k}\mshuffle v-s_{1}(u_{k}\mshuffle v')=\begin{cases}
x_{a_{k+1}}(u_{k+1}\mshuffle v)-x_{a_{k+1}+1}(u_{k+1}\mshuffle v') & k<m\\
0 & k=m,
\end{cases}
\]
which follows from $v=s_{1}(v')$. Similarly, if $t=0$ then $u\dshuffle B(v)$
is equal to the image of
\begin{align*}
 & -\left[\sum_{\substack{k=0\\
k:\,{\rm even}
}
}^{m}s_{d_{k}}\left(u_{k}\mshuffle v-x_{1}(u_{k}\mshuffle v')\right)-\sum_{\substack{k=1\\
k:\,{\rm odd}
}
}^{m}\left(x_{d_{k}}(u_{k}\mshuffle v)-s_{d_{k}+1}(u_{k}\mshuffle v')\right)\right]\\
 & +\left[\sum_{\substack{k=1\\
k:\,{\rm odd}
}
}^{m}x_{1}s_{d_{k}-1}\left(u_{k}\mshuffle v-x_{1}(u_{k}\mshuffle v')\right)-\sum_{\substack{k=1\\
k:\,{\rm even}
}
}^{m}\left(x_{1}x_{d_{k}-1}(u_{k}\mshuffle v)-x_{1}s_{d_{k}}(u_{k}\mshuffle v')\right)\right]\\
 & +u\mshuffle v.\\
 & =0+0+u\mshuffle v=u\mshuffle v
\end{align*}
under $B$. Here, we use
\[
u_{k}\mshuffle v-x_{1}(u_{k}\mshuffle v')=\begin{cases}
x_{a_{k+1}}(u_{k+1}\mshuffle v)-s_{a_{k+1}+1}(u_{k+1}\mshuffle v') & k<m\\
0 & k=m
\end{cases}
\]
which follows from $v=x_{1}v'$. Thus in either case, we have $u\dshuffle B(v)=B(u\mshuffle v)$
as desired.
\end{proof}

\section{Algebraic differential formula for the multivariable block shuffle
product\label{sec:diff_formula}}

In this section, we prove the block shuffle identity by constructing
its multivariable version. We prove a certain compatibility with an
algebraic differential operator. Then, using the compatibility, we
shall prove the block shuffle identities as well as its multi-variable
generalization. 

For a subset $S$ of $\mathbb{C}$, we denote by $\mathcal{W}(S)$
(resp. $\mathcal{A}(S)$) the free monoid (resp. the free algebra
over $\mathbb{Q}$) generated by the formal symbols $\{e_{z}\mid z\in S\}$.
By definition, $\mathcal{A}(S)$ can be regard as the free $\mathbb{Q}$-vector
space generated by $\mathcal{W}(S)$. For $s,t\in\mathbb{C}$, define
the subspace $\mathcal{A}_{s,t}^{0}(S)\subset\mathcal{A}(S)$ by
\[
\mathcal{A}_{s,t}^{0}(S)=\mathbb{Q}\oplus\bigoplus_{k=1}^{\infty}\bigoplus_{\substack{a_{1},\dots,a_{k}\in S\\
a_{1}\neq s,a_{k}\neq t
}
}\mathbb{Q}e_{a_{1}}\cdots e_{a_{k}}.
\]
For $X\subset\mathbb{C}$ and tangential basepoints $\hat{s}=\overrightarrow{u}_{s}$
, $\hat{t}=\overrightarrow{v}_{t}$ at $s,t\in\mathbb{C}$ with tangential
vectors $u,v$, we denote by $\pi_{1}(X,\hat{s},\hat{t})$ the set
of homotopy classes of piece-wisely smooth paths $\gamma$ on $X$
from $\hat{s}$ to $\hat{t}$ (i.e., $\gamma:[0,1]\to\mathbb{C}$
with $\gamma((0,1))\subset X$, $\gamma(0)=s$, $\gamma(1)=t$, $\gamma'(0)=u$
and $\gamma'(1)=-v$). Moreover, we define $\Pi(X;s,t)$ to be the
union of $\pi_{1}(X,\hat{s},\hat{t})$ over all tangential basepoints
$\hat{s}$ and $\hat{t}$ based at $s,t\in\mathbb{C}$.
\begin{defn}
Let $S$ be a subset of $\mathbb{C}$, and $\hat{s}$, $\hat{t}$
tangential basepoints at $s,t\in\mathbb{C}$. For $\gamma\in\pi_{1}(\mathbb{C}\setminus S,\hat{s},\hat{t})$,
we define a $\mathbb{Q}$-linear map $L_{\gamma}:e_{s}\mathcal{A}_{s,t}^{0}(S)e_{t}\to\mathbb{C}$
by
\[
L_{\gamma}(e_{s}e_{z_{1}}\cdots e_{z_{n}}e_{t})\coloneqq I_{\gamma}(\hat{s};z_{1},\dots,z_{n};\hat{t})
\]
where $I_{\gamma}(\hat{s};z_{1},\dots,z_{n};\hat{t})$ is the iterated
integral symbol introduced by Goncharov (see, for example, \cite{Gon05GalSym}\footnote{The definition differs by scaling factor of $(2\pi i)^{n}$.}).
In particular, when $s\neq z_{1}$ and $t\neq z_{2}$, $I_{\gamma}(\hat{s};z_{1},\dots,z_{n};\hat{t})$
is expressed as
\[
\int_{0<t_{1}<\cdots<t_{n}<1}\prod_{j=1}^{n}\frac{d\gamma(t_{j})}{\gamma(t_{j})-z_{j}}.
\]
\end{defn}

\subsection{\label{subsec:Algebraic-differential-formula}Algebraic differential
formula for the multivariable block shuffle product}

We fix complex numbers $a_{1},\dots,a_{n}$ such that $a_{1},\dots,a_{n},1-a_{1},\dots,1-a_{n},0,1$
are distinct. Put 
\[
S=\{a_{1},\dots,a_{n},1-a_{1},\dots,1-a_{n},0,1\}.
\]
Let us fix positive integers $b_{1},\dots,b_{m}$. Note that $x_{b_{1}}\cdots x_{b_{m}}\dshuffle e_{a_{1}}\cdots e_{a_{n}}\in\mathcal{A}(S)$
by definition. For $\alpha,\beta\in S$, we define a $\mathbb{Q}$-linear
map $\partial_{\alpha,\beta}:\mathcal{A}(S)\to\mathcal{A}(S)$ by
$\partial_{\alpha,\beta}(1)=\partial_{\alpha,\beta}(e_{z})=0$ and
\begin{align*}
\partial_{\alpha,\beta}(e_{z_{0}}\cdots e_{z_{k+1}}) & =\sum_{i=1}^{k}\left(\Delta_{z_{i},z_{i+1}}^{\alpha,\beta}-\Delta_{z_{i},z_{i-1}}^{\alpha,\beta}\right)e_{z_{0}}\cdots\widehat{e_{z_{i}}}\cdots e_{z_{k+1}}
\end{align*}
for $k\geq0$, where $\Delta_{x,y}^{\alpha,\beta}\in\{0,1\}$ is defined
by
\[
\Delta_{x,y}^{\alpha,\beta}=\begin{cases}
1 & \{\alpha,\beta\}\in\{\{x,y\},\{1-x,1-y\}\}\\
0 & \{\alpha,\beta\}\notin\{\{x,y\},\{1-x,1-y\}\}.
\end{cases}
\]
Note that for $u\in\mathcal{A}(S)$, $L_{\gamma}(u)$ can be considered
as a complex function of $a_{1},\dots,a_{n}$, and its total differential
is given by
\begin{equation}
dL_{\gamma}(u)=\sum_{1\leq i<j\leq n}L_{\gamma}(\partial_{a_{i},a_{j}}u)\cdot d\log(a_{i}-a_{j})+\sum_{1\leq i\leq j\leq n}L_{\gamma}(\partial_{1-a_{i},a_{j}}u)\cdot d\log(1-a_{i}-a_{j})\label{eq:diff_eq_S}
\end{equation}
using Goncharov's differential formula
\[
dI(z_{0};z_{1},\dots,z_{n};z_{n+1})=\sum_{i=1}^{n}d\log\left(\frac{z_{i+1}-z_{i}}{z_{i}-z_{i-1}}\right)I(z_{0};z_{1},\dots,\widehat{z_{i}},\dots,z_{n};z_{n+1}).
\]
This is the central motivation to calculate $\partial_{\alpha,\beta}(x_{b_{1}}\cdots x_{b_{m}}\dshuffle e_{a_{1}}\cdots e_{a_{n}}\in\mathcal{A}(S))$.
The purpose of this section is to obtain the following formula for
$\partial_{\alpha,\beta}(x_{b_{1}}\cdots x_{b_{m}}\dshuffle e_{a_{1}}\cdots e_{a_{n}})$
(Theorem \ref{thm:AlgDerFormula}):
\begin{align*}
\partial_{\alpha,\beta}(x_{b_{1}}\cdots x_{b_{m}}\dshuffle e_{a_{1}}\cdots e_{a_{n}}) & =x_{b_{1}}\cdots x_{b_{m}}\dshuffle(\partial_{\alpha,\beta}(e_{a_{1}}\cdots e_{a_{n}}))\\
 & \ \ \ +\delta_{b_{1},1}(\Delta_{a_{1},0}^{\alpha,\beta}+\Delta_{a_{1},1}^{\alpha,\beta}-\Delta_{a_{1},a_{1}}^{\alpha,\beta})x_{b_{2}}\cdots x_{b_{m}}\dshuffle e_{a_{1}}\cdots e_{a_{n}}\\
 & \ \ \ -\delta_{b_{m},1}(\Delta_{a_{n},0}^{\alpha,\beta}+\Delta_{a_{n},1}^{\alpha,\beta}-\Delta_{a_{n},a_{n}}^{\alpha,\beta})x_{b_{1}}\cdots x_{b_{m-1}}\dshuffle e_{a_{1}}\cdots e_{a_{n}}.
\end{align*}

\begin{defn}
Define 
\[
D_{\alpha,\beta}(w_{1},a,w_{2})=\sum_{s\in S}\Delta_{a,s}^{\alpha,\beta}(w_{1}e_{s}w_{2,s}-w_{1,s}e_{s}w_{2})
\]
where $w_{1,s},w_{2,s}\in\mathcal{A}$ are defined by
\[
w_{1}=q_{1}+\sum_{s\in S}w_{1,s}e_{s},\ w_{2}=q_{2}+\sum_{s\in S}e_{s}w_{2,s}\ \ \ \ (q_{1},q_{2}\in\mathbb{Q}).
\]
\end{defn}

We put $\mathcal{A}(S)^{+}\coloneqq\mathcal{A}(S)\cap\mathcal{A}^{+}$. 
\begin{lem}
\label{lem:derivation_split}For $w_{1},w_{2}\in\mathcal{A}(S)^{+}$
and $\alpha,\beta,a\in S$, we have
\[
\partial_{\alpha,\beta}(w_{1}e_{a}w_{2})=\partial_{\alpha,\beta}(w_{1}e_{a})w_{2}+w_{1}\partial_{\alpha,\beta}(e_{a}w_{2})+D_{\alpha,\beta}(w_{1},a,w_{2}).
\]
\end{lem}

This lemma is an immediate consequence of the definition of $\partial_{\alpha,\beta}$
and $D_{\alpha,\beta}$. For $t\in\{0,1\}$ and $b\in\mathbb{Z}$,
we put
\[
W(t;b)=\begin{cases}
e_{\iota^{0}(t)}e_{\iota^{1}(t)}\cdots e_{\iota^{b-1}(t)} & b\geq0\\
0 & b<0.
\end{cases}
\]
For convenience, hereafter we regard $U(b)=0$ for $b<0$. Notice
that $W(0;b)+W(1;b)=U(b)$ definition.
\begin{lem}
\label{lem:Dsum}Fix $\alpha,\beta\in S$. For $a\in S\setminus\{0,1\}$,
$k\geq2$ and $b_{1},\dots,b_{k}\in\mathbb{Z}_{\geq1}$, $s,s'\in\{0,1\}$
we have
\begin{align*}
 & \sum_{j=1}^{k-1}D_{\alpha,\beta}\left(U(b_{1}+\cdots+b_{j}-s),\iota^{b_{1}+\cdots+b_{j}-j}(a),U(b_{j+1}+\cdots+b_{k}-s')\right)\\
 & =-(-1)^{s}\sum_{t=0}^{1}(\Delta_{a,t}^{\alpha,\beta}-\Delta_{a,\iota^{k+1}(t)}^{\alpha,\beta})W(t;b_{1}+\cdots+b_{k}-s-s')\\
 & \ \ +(\Delta_{a,0}^{\alpha,\beta}+\Delta_{a,1}^{\alpha,\beta})\delta_{b_{1},1}\delta_{s,1}U(b_{1}+\cdots+b_{k}-s-s')\\
 & \ \ -(\Delta_{a,0}^{\alpha,\beta}+\Delta_{a,1}^{\alpha,\beta})\delta_{b_{k},1}\delta_{s',1}U(b_{1}+\cdots+b_{k}-s-s').
\end{align*}
\end{lem}

\begin{proof}
Note that for $b,b'\in\mathbb{Z}_{\geq0}$ and $c\in\{a,1-a\}$, we
have
\begin{align*}
D_{\alpha,\beta}\left(U(b),c,U(b')\right) & =\sum_{t=0}^{1}(1-\delta_{b',0})\Delta_{c,t}^{\alpha,\beta}U(b)W(t;b')-\sum_{t=0}^{1}(1-\delta_{b,0})\Delta_{c,t}^{\alpha,\beta}W(\iota^{b+1}(t);b)U(b')\\
 & =\sum_{t=0}^{1}\Delta_{c,t}^{\alpha,\beta}\left(U(b)W(t;b')-W(\iota^{b+1}(t);b)U(b')\right)-\sum_{t=0}^{1}\delta_{b',0}\Delta_{c,t}^{\alpha,\beta}U(b)+\sum_{t=0}^{1}\delta_{b,0}\Delta_{c,t}^{\alpha,\beta}U(b')\\
 & =\sum_{t=0}^{1}\Delta_{c,t}^{\alpha,\beta}\left(\left(W(\iota^{b+1}(t);b)+W(\iota^{b}(t);b)\right)W(t;b')-W(\iota^{b+1}(t);b)\left(W(t;b')+W(\iota(t);b')\right)\right)\\
 & \ \ \ +\left(\delta_{b,0}U(b')-\delta_{b',0}U(b)\right)\sum_{t=0}^{1}\Delta_{c,t}^{\alpha,\beta}.\\
 & =\sum_{t=0}^{1}\Delta_{c,t}^{\alpha,\beta}\left(W(\iota^{b}(t);b+b')-W(\iota^{b+1}(t);b+b')\right)+\left(\delta_{b,0}-\delta_{b',0}\right)(\Delta_{c,0}^{\alpha,\beta}+\Delta_{c,1}^{\alpha,\beta})\cdot U(b+b')\\
 & =(-1)^{b}(\Delta_{c,0}^{\alpha,\beta}-\Delta_{c,1}^{\alpha,\beta})\cdot U^{-}(b+b')+\left(\delta_{b,0}-\delta_{b',0}\right)(\Delta_{c,0}^{\alpha,\beta}+\Delta_{c,1}^{\alpha,\beta})\cdot U(b+b')
\end{align*}
where
\[
U^{-}(b)\coloneqq W(0;b)-W(1;b).
\]
Therefore we have
\begin{align*}
 & \sum_{j=1}^{k-1}D_{\alpha,\beta}\left(U(b_{1}+\cdots+b_{j}-s),\iota^{b_{1}+\cdots+b_{j}-j}(a),U(b_{j+1}+\cdots+b_{k}-s')\right)\\
 & =\sum_{j=1}^{k-1}(-1)^{b_{1}+\cdots+b_{j}-s}\left(\Delta_{\iota^{b_{1}+\cdots+b_{j}-j}(a),0}^{\alpha,\beta}-\Delta_{\iota^{b_{1}+\cdots+b_{j}-j}(a),1}^{\alpha,\beta}\right)U^{-}(b_{1}+\cdots+b_{k}-s-s')\\
 & \ \ \ +\sum_{j=1}^{k-1}\left(\delta_{b_{1}+\cdots+b_{j}-s,0}-\delta_{b_{j+1}+\cdots+b_{k}-s',0}\right)\left(\Delta_{\iota^{b_{1}+\cdots+b_{j}-j}(a),0}^{\alpha,\beta}+\Delta_{\iota^{b_{1}+\cdots+b_{j}-j}(a),1}^{\alpha,\beta}\right)U(b_{1}+\cdots+b_{k}-s-s')\\
 & =\sum_{j=1}^{k-1}(-1)^{j-s}\left(\Delta_{a,0}^{\alpha,\beta}-\Delta_{a,1}^{\alpha,\beta}\right)U^{-}(b_{1}+\cdots+b_{k}-s-s')\\
 & \ \ \ +\sum_{j=1}^{k-1}\left(\delta_{j,1}\delta_{b_{1},1}\delta_{s,1}-\delta_{j,k-1}\delta_{b_{k},1}\delta_{s',1}\right)\left(\Delta_{a,0}^{\alpha,\beta}+\Delta_{a,1}^{\alpha,\beta}\right)U(b_{1}+\cdots+b_{k}-s-s')\\
 & =\left(\Delta_{a,0}^{\alpha,\beta}-\Delta_{a,1}^{\alpha,\beta}\right)U^{-}(b_{1}+\cdots+b_{k}-s-s')\sum_{j=1}^{k-1}(-1)^{j-s}\\
 & \ \ \ +\left(\delta_{b_{1},1}\delta_{s,1}-\delta_{b_{k},1}\delta_{s',1}\right)\left(\Delta_{a,0}^{\alpha,\beta}+\Delta_{a,1}^{\alpha,\beta}\right)U(b_{1}+\cdots+b_{k}-s-s')\\
 & =\begin{cases}
(-1)^{s+1}\left(\Delta_{a,0}^{\alpha,\beta}-\Delta_{a,1}^{\alpha,\beta}\right)U^{-}(b_{1}+\cdots+b_{k}-s-s') & k:\,{\rm even}\\
0 & k:\,{\rm odd}
\end{cases}\\
 & \ \ \ +\left(\delta_{b_{1},1}\delta_{s,1}-\delta_{b_{k},1}\delta_{s',1}\right)\left(\Delta_{a,0}^{\alpha,\beta}+\Delta_{a,1}^{\alpha,\beta}\right)U(b_{1}+\cdots+b_{k}-s-s')\\
 & =(-1)^{s+1}\left(\Delta_{a,0}^{\alpha,\beta}-\Delta_{a,\iota^{k}(1)}^{\alpha,\beta}\right)U^{-}(b_{1}+\cdots+b_{k}-s-s')\\
 & \ \ \ +\left(\delta_{b_{1},1}\delta_{s,1}-\delta_{b_{k},1}\delta_{s',1}\right)\left(\Delta_{a,0}^{\alpha,\beta}+\Delta_{a,1}^{\alpha,\beta}\right)U(b_{1}+\cdots+b_{k}-s-s').
\end{align*}
This completes the proof.
\end{proof}

\begin{prop}
\label{prop:AlgDerFormula_a0}For $\alpha\in S\setminus\{0,1\}$,
$\beta\in\{0,1\}$, $m>0$ and $b_{1},\dots,b_{m}\in\mathbb{Z}_{>0}$,
we have
\begin{align*}
\partial_{\alpha,\beta}(x_{b_{1}}\cdots x_{b_{m}}\dshuffle e_{a_{1}}\cdots e_{a_{n}}) & =(\Delta_{a_{1},0}^{\alpha,\beta}+\Delta_{a_{1},1}^{\alpha,\beta})\delta_{b_{1},1}x_{b_{2}}\cdots x_{b_{m}}\dshuffle e_{a_{1}}\cdots e_{a_{n}}\\
 & \ \ -(\Delta_{a_{n},0}^{\alpha,\beta}+\Delta_{a_{n},1}^{\alpha,\beta})\delta_{b_{m},1}x_{b_{1}}\cdots x_{b_{m-1}}\dshuffle e_{a_{1}}\cdots e_{a_{n}}.
\end{align*}
\end{prop}

\begin{proof}
We prove the claim by induction on $n+m.$ The case $m=n=1$ is easily
checked from the definition. Suppose that $(m,n)\neq(1,1)$. Decomposing
$x_{b_{1}}\cdots x_{b_{m}}\dshuffle e_{a_{1}}\cdots e_{a_{n}}$ as
\[
x_{b_{1}}\cdots x_{b_{m}}\dshuffle e_{a_{1}}\cdots e_{a_{n}}=\sum_{l}w_{l}e_{p_{l}}w'_{l}\quad\left(w_{l}\in e_{a_{1}}\mathcal{A}(\{0,1\}),w_{l}'\in\mathcal{A},p_{l}\in\{a_{1},a_{2}\}\right),
\]
and applying Lemma \ref{lem:derivation_split}, we obtain
\[
\partial_{\alpha,\beta}(x_{b_{1}}\cdots x_{b_{m}}\dshuffle e_{a_{1}}\cdots e_{a_{n}})=F_{1}+F_{2}+F_{3}
\]
where
\begin{align*}
F_{1} & =\sum_{l}\partial_{\alpha,\beta}(w_{l}e_{p_{l}})w'_{l}\\
F_{2} & =\sum_{l}w_{l}\partial_{\alpha,\beta}(e_{p_{l}}w'_{l})\\
F_{3} & =\sum_{l}D_{\alpha,\beta}(w_{l},p_{l},w'_{l}).
\end{align*}
To compute these $F_{i}$'s further, let us recall the recurrence
relation of $\dshuffle$-product (the first formula of Proposition
\ref{prop:recur_dshuffle}). We put $c_{k}=\sum_{j=1}^{k}(b_{j}-1)$,
$c_{k}'=\sum_{j=m-k+1}^{m}(b_{j}-1)$,
\begin{align*}
g_{k,i} & =\begin{cases}
(-1)^{k+i}\boldsymbol{\iota}^{c_{k}}(x_{b_{k+1}}\cdots x_{b_{m}}\dshuffle e_{a_{i+1}}\cdots e_{a_{n}}) & 0\leq i<n\\
0 & i\geq n
\end{cases},
\end{align*}
\[
G_{k,i}=\begin{cases}
(-1)^{k+i}\boldsymbol{\iota}^{c_{k}}(x_{b_{k+1}}\cdots x_{b_{m-1}}\dshuffle e_{a_{i+1}}\cdots e_{a_{n}}) & 0\leq i<n\\
0 & i\geq n
\end{cases},
\]
and
\[
Y_{k,t}(s)=e_{a_{1}}W(t;c_{k}+k+s-2)g_{k,s}.
\]
Then the recurrence relation is expressed as
\begin{equation}
x_{b_{1}}\cdots x_{b_{m}}\dshuffle e_{a_{1}}\cdots e_{a_{n}}=-\sum_{s=0}^{1}\sum_{k=1}^{m}e_{a_{1}}U(c_{k}+k+s-1)g_{k,s}-e_{a_{1}}g_{0,1}.\label{eq:ds_recur_byg}
\end{equation}
We put $a_{j}=a_{1}$ for $j\notin\{1,\dots,n\}$ to avoid undefined
notations. For the calculation of $F_{1}$, notice that $g_{k,s}\in e_{\iota^{c_{k}}(a_{s+1})}\mathcal{A}(S)$
and
\[
\partial_{\alpha,\beta}(e_{a_{1}}U(c_{k}+k+s-1)e_{\iota^{c_{k}}(a_{s+1})})=\sum_{t=0}^{1}\bigg(\Delta_{\iota^{c_{k}}(a_{s+1}),\iota^{c_{k}+k+s-2}(t)}^{\alpha,\beta}-\sum_{t=0}^{1}\Delta_{a_{1},\iota(t)}^{\alpha,\beta}\bigg)e_{a_{1}}W(t;c_{k}+k+s-2)e_{\iota^{c_{k}}(a_{s+1})}
\]
since if $c_{k}+k+s-1\geq1$
\begin{align*}
\partial_{\alpha,\beta}(e_{a_{1}}U(c_{k}+k+s-1)e_{\iota^{c_{k}}(a_{s+1})}) & =\partial_{\alpha,\beta}\bigg(\sum_{t=0}^{1}e_{a_{1}}\overbrace{e_{\iota(t)}e_{\iota^{2}(t)}\cdots e_{\iota^{c_{k}+k+s-1}(t)}}^{c_{k}+k+s-1}e_{\iota^{c_{k}}(a_{s+1})}\bigg)\\
 & =\sum_{t=0}^{1}\Delta_{\iota^{c_{k}}(a_{s+1}),\iota^{c_{k}+k+s-1}(t)}^{\alpha,\beta}e_{a_{1}}\overbrace{e_{\iota(t)}e_{\iota^{2}(t)}\cdots e_{\iota^{c_{k}+k+s-2}(t)}}^{c_{k}+k+s-2}e_{\iota^{c_{k}}(a_{s+1})}\\
 & \quad-\sum_{t=0}^{1}\Delta_{a_{1},\iota(t)}^{\alpha,\beta}e_{a_{1}}\overbrace{e_{\iota^{2}(t)}\cdots e_{\iota^{c_{k}+k+s-1}(t)}}^{c_{k}+k+s-2}e_{\iota^{c_{k}}(a_{s+1})}\\
 & =\sum_{t=0}^{1}\bigg(\Delta_{\iota^{c_{k}}(a_{s+1}),\iota^{c_{k}+k+s}(t)}^{\alpha,\beta}-\sum_{t=0}^{1}\Delta_{a_{1},\iota(t)}^{\alpha,\beta}\bigg)e_{a_{1}}W(t;c_{k}+k+s-2)e_{\iota^{c_{k}}(a_{s+1})}
\end{align*}
and the case $c_{k}+k+s-1=0$ is trivial. Thus
\begin{align*}
F_{1} & =\sum_{s=0}^{1}\sum_{k=1}^{m}\sum_{t=0}^{1}\left(\Delta_{a_{1},\iota(t)}^{\alpha,\beta}-\Delta_{\iota^{c_{k}}(a_{s+1}),\iota^{c_{k}+k+s}(t)}^{\alpha,\beta}\right)e_{a_{1}}W(t;c_{k}+k+s-2)g_{k,s}\\
 & =\sum_{s=0}^{1}\sum_{k=1}^{m}\sum_{t=0}^{1}\left(\Delta_{a_{1},\iota(t)}^{\alpha,\beta}-\Delta_{a_{s+1},\iota^{k+s}(t)}^{\alpha,\beta}\right)Y_{k,t}(s).
\end{align*}
Next, we have $F_{2}=F_{2}^{(1)}+F_{2}^{(2)}$ where
\begin{align*}
F_{2}^{(1)} & =-\sum_{s=0}^{1}\sum_{k=1}^{m}e_{a_{1}}U(c_{k}+k+s-1)\partial_{\alpha,\beta}(g_{k,s}),\\
F_{2}^{(2)} & =-e_{a_{1}}\partial_{\alpha,\beta}(g_{0,1}).
\end{align*}
By the induction hypothesis and the fact $\partial_{\alpha,\beta}(g_{m,s})=0$,
$F_{2}^{(1)}$ is further decomposed as $F_{2}^{(1)}=F_{2}^{(1,1)}+F_{2}^{(1,2)}$
where
\begin{align*}
F_{2}^{(1,1)} & =\sum_{s=0}^{1}\sum_{k=1}^{m-1}(\Delta_{a_{s+1},0}^{\alpha,\beta}+\Delta_{a_{s+1},1}^{\alpha,\beta})\delta_{b_{k+1},1}\cdot e_{a_{1}}U(c_{k+1}+k+s-1)g_{k+1,s}\\
 & =\sum_{s=0}^{1}\sum_{k=2}^{m}(\Delta_{a_{s+1},0}^{\alpha,\beta}+\Delta_{a_{s+1},1}^{\alpha,\beta})\delta_{b_{k},1}\cdot e_{a_{1}}U(c_{k}+k+s-2)g_{k,s}\\
 & =\sum_{s=0}^{1}\sum_{k=1}^{m}(\Delta_{a_{s+1},0}^{\alpha,\beta}+\Delta_{a_{s+1},1}^{\alpha,\beta})\delta_{b_{k},1}\cdot e_{a_{1}}U(c_{k}+k+s-2)g_{k,s}-\sum_{s=0}^{1}(\Delta_{a_{s+1},0}^{\alpha,\beta}+\Delta_{a_{s+1},1}^{\alpha,\beta})\delta_{b_{1},1}\cdot e_{a_{1}}U(c_{1}+s-1)g_{1,s}\\
 & =\sum_{s=0}^{1}\sum_{k=1}^{m}\sum_{t=0}^{1}(\Delta_{a_{s+1},0}^{\alpha,\beta}+\Delta_{a_{s+1},1}^{\alpha,\beta})\delta_{b_{k},1}\cdot Y_{k,t}(s)-2(\Delta_{a_{2},0}^{\alpha,\beta}+\Delta_{a_{2},1}^{\alpha,\beta})\delta_{b_{1},1}e_{a_{1}}g_{1,1}.
\end{align*}
and
\begin{align*}
F_{2}^{(1,2)} & =\sum_{s=0}^{1}\sum_{k=1}^{m-1}(\Delta_{a_{n},0}^{\alpha,\beta}+\Delta_{a_{n},1}^{\alpha,\beta})\delta_{b_{m},1}\cdot e_{a_{1}}U(c_{k}+k+s-1)G_{k,s}\\
 & =-(\Delta_{a_{n},0}^{\alpha,\beta}+\Delta_{a_{n},1}^{\alpha,\beta})\delta_{b_{m},1}\cdot\left(x_{b_{1}}\cdots x_{b_{m-1}}\dshuffle e_{a_{1}}\cdots e_{a_{n}}+e_{a_{1}}G{}_{0,1}\right),
\end{align*}
where in the last equality we have used the recurrence relation for
$x_{b_{1}}\cdots x_{b_{m-1}}\dshuffle e_{a_{1}}\cdots e_{a_{n}}$.
By the induction hypothesis, for $n\geq2$, we have 
\begin{align*}
F_{2}^{(2)} & =e_{a_{1}}\partial_{\alpha,\beta}(x_{b_{1}}\cdots x_{b_{m}}\dshuffle e_{a_{2}}\cdots e_{a_{n}})\\
 & =(\Delta_{a_{2},0}^{\alpha,\beta}+\Delta_{a_{2},1}^{\alpha,\beta})\delta_{b_{1},1}\cdot e_{a_{1}}g_{1,1}+(\Delta_{a_{n},0}^{\alpha,\beta}+\Delta_{a_{n},1}^{\alpha,\beta})\delta_{b_{m},1}\cdot e_{a_{1}}G_{0,1}.
\end{align*}
Now, let us compute $F_{3}$. Since
\[
x_{b_{1}}\cdots x_{b_{m}}\dshuffle e_{a_{1}}\cdots e_{a_{n}}=-\sum_{s=0}^{1}\sum_{j=1}^{m}e_{a_{1}}U(c_{j}+j+s-1)g_{j,s}-e_{a_{1}}g_{0,1}
\]
and
\begin{align*}
g_{j,s} & =-\sum_{s'=0}^{1}\sum_{j<k\leq m}\bm{\iota}^{c_{j}}(e_{a_{s+1}})U(b_{j+1}+\cdots+b_{k}+s'-1)g_{k,s+s'}-\bm{\iota}^{c_{j}}(e_{a_{s+1}})g_{j,s+1},\\
g_{0,1} & =-\sum_{s=0}^{1}\sum_{k=1}^{m}e_{a_{2}}U(c_{k}+k+s-1)g_{k,s+1}-e_{a_{2}}g_{0,2},
\end{align*}
we have
\begin{align*}
x_{b_{1}}\cdots x_{b_{m}}\dshuffle e_{a_{1}}\cdots e_{a_{n}} & =\sum_{s=0}^{1}\sum_{s'=0}^{1}\sum_{1\leq j<k\leq m}e_{a_{1}}U(c_{j}+j+s-1)\bm{\iota}^{c_{j}}(e_{a_{s+1}})U(b_{j+1}+\cdots+b_{k}+s'-1)g_{k,s+s'}\\
 & \ \ +\sum_{\substack{s=0}
}^{1}\sum_{k=1}^{m}e_{a_{1}}U(c_{k}+k+s-1)\bm{\iota}^{c_{k}}(e_{a_{s+1}})g_{k,s+1}\\
 & \ \ +\sum_{\substack{s=0}
}^{1}\sum_{k=1}^{m}e_{a_{1}}e_{a_{2}}U(c_{k}+k+s-1)g_{k,s+1}\\
 & \ \ +e_{a_{1}}e_{a_{2}}g_{0,2}.
\end{align*}
Therefore, we have
\begin{align*}
F_{3} & =\sum_{s=0}^{1}\sum_{s'=0}^{1}\sum_{1\leq j<k\leq m}e_{a_{1}}D_{\alpha,\beta}(U(b_{1}+\cdots+b_{j}+s-1),\iota^{c_{j}}(a_{s+1}),U(b_{j+1}+\cdots+b_{k}+s'-1))g_{k,s+s'}\\
 & \quad-\sum_{s=0}^{1}\sum_{k=1}^{m}\sum_{t=0}^{1}(1-\delta_{b_{1},1}\delta_{k,1}\delta_{s,0})(\Delta_{a_{s+1},\iota^{k+s}(t)}^{\alpha,\beta})e_{a_{1}}W(t,c_{k}+k+s-1)g_{k,s+1}\\
 & \quad+\sum_{s=0}^{1}\sum_{k=1}^{m}\sum_{t=0}^{1}(1-\delta_{b_{1},1}\delta_{k,1}\delta_{s,0})(\Delta_{a_{2},t}^{\alpha,\beta})e_{a_{1}}W(t,c_{k}+k+s-1)g_{k,s+1}\\
 & =\sum_{s=0}^{1}\sum_{s'=0}^{1}\sum_{1\leq j<k\leq m}e_{a_{1}}D_{\alpha,\beta}(U(b_{1}+\cdots+b_{j}+s-1),\iota^{c_{j}}(a_{s+1}),U(b_{j+1}+\cdots+b_{k}+s'-1))g_{k,s+s'}\\
 & \ \ +\sum_{s=0}^{1}\sum_{k=1}^{m}\sum_{t=0}^{1}(1-\delta_{b_{1},1}\delta_{k,1}\delta_{s,0})(\Delta_{a_{2},t}^{\alpha,\beta}-\Delta_{a_{s+1},\iota^{k-s}(t)}^{\alpha,\beta})Y_{k,t}(s+1)\\
 & =F_{3}^{(1)}+F_{3}^{(2)}+F_{3}^{(3)},
\end{align*}
where
\[
F_{3}^{(1)}=\sum_{s=0}^{1}\sum_{s'=0}^{1}\sum_{1\leq j<k\leq m}e_{a_{1}}D_{\alpha,\beta}(U(b_{1}+\cdots+b_{j}+s-1),\iota^{c_{j}}(a_{s+1}),U(b_{j+1}+\cdots+b_{k}+s'-1))g_{k,s+s'},
\]
\begin{align*}
F_{3}^{(2)} & =\sum_{s=0}^{1}\sum_{k=1}^{m}\sum_{t=0}^{1}(\Delta_{a_{2},t}^{\alpha,\beta}-\Delta_{a_{s+1},\iota^{k-s}(t)}^{\alpha,\beta})Y_{k,t}(s+1),
\end{align*}
and
\begin{align*}
F_{3}^{(3)} & =-\sum_{t=0}^{1}(\Delta_{a_{2},t}^{\alpha,\beta}-\Delta_{a_{1},1-t}^{\alpha,\beta})\delta_{b_{1},1}Y_{1,t}(1)\\
 & =-\sum_{t=0}^{1}(\Delta_{a_{2},t}^{\alpha,\beta}-\Delta_{a_{1},1-t}^{\alpha,\beta})\delta_{b_{1},1}e_{a_{1}}W(t;b_{1}-1)g_{1,1}\\
 & =-\sum_{t=0}^{1}(\Delta_{a_{2},t}^{\alpha,\beta}-\Delta_{a_{1},1-t}^{\alpha,\beta})\delta_{b_{1},1}e_{a_{1}}g_{1,1}\\
 & =(\Delta_{a_{1},0}^{\alpha,\beta}+\Delta_{a_{1},1}^{\alpha,\beta}-\Delta_{a_{2},0}^{\alpha,\beta}-\Delta_{a_{2},1}^{\alpha,\beta})\delta_{b_{1},1}e_{a_{1}}g_{1,1}.
\end{align*}
The term $F_{3}^{(1)}$ can be computed using Lemma \ref{lem:Dsum}.
In fact, by setting $a\mapsto a_{s+1}$, $s\mapsto1-s$, $s'\mapsto1-s'$
in Lemma \ref{lem:Dsum}, we find
\begin{align*}
 & \sum_{1\leq j<k}D_{\alpha,\beta}\left(U(b_{1}+\cdots+b_{j}+s-1),\iota^{c_{j}}(a_{s+1}),U(b_{j+1}+\cdots+b_{k}+s'-1)\right)\\
 & =(-1)^{s}\sum_{t=0}^{1}(\Delta_{a_{s+1},t}^{\alpha,\beta}-\Delta_{a_{s+1},\iota^{k}(1-t)}^{\alpha,\beta})W(t;b_{1}+\cdots+b_{k}+s+s'-2)\\
 & \ \ -(\Delta_{a_{s+1},0}^{\alpha,\beta}+\Delta_{a_{s+1},1}^{\alpha,\beta})\delta_{b_{k},1}\delta_{s',0}U(b_{1}+\cdots+b_{k}+s+s'-2)\\
 & \ \ +(\Delta_{a_{s+1},0}^{\alpha,\beta}+\Delta_{a_{s+1},1}^{\alpha,\beta})\delta_{b_{1},1}\delta_{s,0}U(b_{1}+\cdots+b_{k}+s+s'-2)
\end{align*}
for $k\geq2$, and thus we have
\[
F_{3}^{(1)}=F_{3}^{(1,1)}+F_{3}^{(1,2)}+F_{3}^{(1,3)}
\]
where
\begin{align*}
F_{3}^{(1,1)} & =\sum_{s=0}^{1}\sum_{s'=0}^{1}\sum_{k=2}^{m}e_{a_{1}}\Bigg((-1)^{s}\sum_{t=0}^{1}(\Delta_{a_{s+1},t}^{\alpha,\beta}-\Delta_{a_{s+1},\iota^{k+1}(t)}^{\alpha,\beta})W(t;b_{1}+\cdots+b_{k}+s+s'-2)\Bigg)g_{k,s+s'}\\
 & =\sum_{s=0}^{1}\sum_{s'=0}^{1}\sum_{k=2}^{m}\sum_{t=0}^{1}(-1)^{s}(\Delta_{a_{s+1},t}^{\alpha,\beta}-\Delta_{a_{s+1},\iota^{k+1}(t)}^{\alpha,\beta})Y_{k,t}(s+s')\\
 & =\sum_{s=0}^{1}\sum_{s'=0}^{1}\sum_{k=1}^{m}\sum_{t=0}^{1}(-1)^{s}(\Delta_{a_{s+1},t}^{\alpha,\beta}-\Delta_{a_{s+1},\iota^{k+1}(t)}^{\alpha,\beta})Y_{k,t}(s+s'),
\end{align*}
\begin{align*}
F_{3}^{(1,2)} & =\sum_{s=0}^{1}\sum_{s'=0}^{1}\sum_{k=2}^{m}e_{a_{1}}\Bigg((\Delta_{a_{s+1},0}^{\alpha,\beta}+\Delta_{a_{s+1},1}^{\alpha,\beta})\delta_{b_{1},1}\delta_{s,0}U(b_{1}+\cdots+b_{k}+s+s'-2)\Bigg)g_{k,s+s'}\\
 & =\sum_{s'=0}^{1}\sum_{k=2}^{m}e_{a_{1}}(\Delta_{a_{1},0}^{\alpha,\beta}+\Delta_{a_{1},1}^{\alpha,\beta})\delta_{b_{1},1}U(b_{1}+\cdots+b_{k}+s'-2)g_{k,s'}\\
 & =(\Delta_{a_{1},0}^{\alpha,\beta}+\Delta_{a_{1},1}^{\alpha,\beta})\delta_{b_{1},1}\sum_{s'=0}^{1}\sum_{k=2}^{m}e_{a_{1}}U(b_{1}+\cdots+b_{k}+s'-2)g_{k,s'}\\
 & =(\Delta_{a_{1},0}^{\alpha,\beta}+\Delta_{a_{1},1}^{\alpha,\beta})\delta_{b_{1},1}\sum_{s'=0}^{1}\sum_{k=2}^{m}e_{a_{1}}U(b_{2}+\cdots+b_{k}+s'-1)g_{k,s'}\\
 & =(\Delta_{a_{1},0}^{\alpha,\beta}+\Delta_{a_{1},1}^{\alpha,\beta})\delta_{b_{1},1}\bigg(x_{b_{2}}\cdots x_{b_{m}}\dshuffle e_{a_{1}}\cdots e_{a_{n}}-e_{a_{1}}g_{1,1}\bigg),
\end{align*}
and
\begin{align*}
F_{3}^{(1,3)} & =\sum_{s=0}^{1}\sum_{s'=0}^{1}\sum_{k=2}^{m}e_{a_{1}}\Bigg(-(\Delta_{a_{s+1},0}^{\alpha,\beta}+\Delta_{a_{s+1},1}^{\alpha,\beta})\delta_{b_{k},1}\delta_{s',0}U(b_{1}+\cdots+b_{k}+s+s'-2)\Bigg)g_{k,s+s'}\\
 & =-\sum_{s=0}^{1}\sum_{k=2}^{m}(\Delta_{a_{s+1},0}^{\alpha,\beta}+\Delta_{a_{s+1},1}^{\alpha,\beta})\delta_{b_{k},1}e_{a_{1}}U(b_{1}+\cdots+b_{k}+s-2)g_{k,s}\\
 & =-\sum_{s=0}^{1}\sum_{k=2}^{m}\sum_{t=0}^{1}(\Delta_{a_{s+1},0}^{\alpha,\beta}+\Delta_{a_{s+1},1}^{\alpha,\beta})\delta_{b_{k},1}e_{a_{1}}W(t;b_{1}+\cdots+b_{k}+s-2)g_{k,s}\\
 & =-\sum_{s=0}^{1}\sum_{k=2}^{m}\sum_{t=0}^{1}(\Delta_{a_{s+1},0}^{\alpha,\beta}+\Delta_{a_{s+1},1}^{\alpha,\beta})\delta_{b_{k},1}Y_{k,t}(s)\\
 & =-\sum_{s=0}^{1}\sum_{k=1}^{m}\sum_{t=0}^{1}(\Delta_{a_{s+1},0}^{\alpha,\beta}+\Delta_{a_{s+1},1}^{\alpha,\beta})\delta_{b_{k},1}Y_{k,t}(s)+\sum_{s=0}^{1}\sum_{t=0}^{1}(\Delta_{a_{s+1},0}^{\alpha,\beta}+\Delta_{a_{s+1},1}^{\alpha,\beta})\delta_{b_{1},1}Y_{1,t}(s)\\
 & =-\sum_{s=0}^{1}\sum_{k=1}^{m}\sum_{t=0}^{1}(\Delta_{a_{s+1},0}^{\alpha,\beta}+\Delta_{a_{s+1},1}^{\alpha,\beta})\delta_{b_{k},1}Y_{k,t}(s)+\sum_{s=0}^{1}\sum_{t=0}^{1}(\Delta_{a_{s+1},0}^{\alpha,\beta}+\Delta_{a_{s+1},1}^{\alpha,\beta})\delta_{b_{1},1}e_{a_{1}}W(t;c_{1}+1+s-2)g_{1,s}\\
 & =-\sum_{s=0}^{1}\sum_{k=1}^{m}\sum_{t=0}^{1}(\Delta_{a_{s+1},0}^{\alpha,\beta}+\Delta_{a_{s+1},1}^{\alpha,\beta})\delta_{b_{k},1}Y_{k,t}(s)+2(\Delta_{a_{2},0}^{\alpha,\beta}+\Delta_{a_{2},1}^{\alpha,\beta})\delta_{b_{1},1}e_{a_{1}}g_{1,1}.
\end{align*}
Thus, finally we have
\begin{align*}
 & \partial_{\alpha,\beta}(x_{b_{1}}\cdots x_{b_{m}}\dshuffle e_{a_{1}}\cdots e_{a_{n}})\\
 & =F_{1}+F_{2}^{(1,1)}+F_{2}^{(1,2)}+F_{2}^{(2)}+F_{3}^{(1,1)}+F_{3}^{(1,2)}+F_{3}^{(1,3)}+F_{3}^{(2)}+F_{3}^{(3)}\\
 & =\delta_{b_{1},1}(\Delta_{a_{1},0}^{\alpha,\beta}+\Delta_{a_{1},1}^{\alpha,\beta})\left(x_{b_{2}}\cdots x_{b_{m}}\dshuffle e_{a_{1}}\cdots e_{a_{n}}\right)-\delta_{b_{m},1}(\Delta_{a_{n},0}^{\alpha,\beta}+\Delta_{a_{n},1}^{\alpha,\beta})\left(x_{b_{1}}\cdots x_{b_{m-1}}\dshuffle e_{a_{1}}\cdots e_{a_{n}}\right)\\
 & \quad+\sum_{k=1}^{m}\sum_{t=0}^{1}\left(\Delta_{a_{1},\iota(t)}^{\alpha,\beta}+\Delta_{a_{1},t}^{\alpha,\beta}-\Delta_{a_{1},\iota^{k}(1-t)}^{\alpha,\beta}-\Delta_{a_{1},\iota^{k}(t)}^{\alpha,\beta}\right)\bigg(Y_{k,t}(0)+Y_{k,t}(1)\bigg)\\
 & =\delta_{b_{1},1}(\Delta_{a_{1},0}^{\alpha,\beta}+\Delta_{a_{1},1}^{\alpha,\beta})\left(x_{b_{2}}\cdots x_{b_{m}}\dshuffle e_{a_{1}}\cdots e_{a_{n}}\right)-\delta_{b_{m},1}(\Delta_{a_{n},0}^{\alpha,\beta}+\Delta_{a_{n},1}^{\alpha,\beta})\left(x_{b_{1}}\cdots x_{b_{m-1}}\dshuffle e_{a_{1}}\cdots e_{a_{n}}\right),
\end{align*}
which completes the proof.
\end{proof}
\begin{prop}
\label{lem:AlgDerFormula_aa}Assume that $n>1$. For $s\in\{1,\dots,n-1\}$,
$m>0$ and $b_{1},\dots,b_{m}\in\mathbb{Z}_{>0}$, we have
\begin{align*}
\partial_{a_{s},a_{s+1}}(x_{b_{1}}\cdots x_{b_{m}}\dshuffle e_{a_{1}}\cdots e_{a_{n}}) & =x_{b_{1}}\cdots x_{b_{m}}\dshuffle\partial_{a_{s},a_{s+1}}(e_{a_{1}}\cdots e_{a_{n}}).
\end{align*}
\end{prop}

\begin{proof}
We put
\begin{align*}
c_{k} & \coloneqq\sum_{j=1}^{k}(b_{j}-1),\\
d_{k} & \coloneqq\sum_{j=1}^{k}b_{j}=c_{k}+k,
\end{align*}
and
\begin{align*}
X_{i,j} & =X_{i,j}'\bm{\iota}^{c_{i}}(e_{a_{j}})=x_{b_{1}}\cdots x_{b_{i}}\dshuffle e_{a_{1}}\cdots e_{a_{j}}\\
Y_{i,j} & =e_{a_{j+1}}Y_{i,j}'=x_{b_{i+1}}\cdots x_{b_{m}}\dshuffle e_{a_{j+1}}\cdots e_{a_{n}}.
\end{align*}
By Proposition \ref{prop:dsh_sep_as_as1}
\[
x_{b_{1}}\cdots x_{b_{m}}\dshuffle e_{a_{1}}\cdots e_{a_{n}}=\sum_{k=0}^{m}X_{k,s}\cdot\bm{\iota}^{c_{k}}(Y_{k,s})+\sum_{0\leq k<l\leq m}(-1)^{l-k}X_{k,s}\cdot U(d_{l}-d_{k})\cdot\bm{\iota}^{c_{l}}(Y_{l,s}).
\]
Since
\[
\partial_{a_{s},a_{s+1}}(X_{k,s}\bm{\iota}^{c_{k}}(Y_{k,s}))=\begin{cases}
X_{k,s}'\iota_{\mathcal{A}}^{c_{k}}(Y_{k,s}) & (k,s)\neq(0,1)\\
0 & (k,s)=(0,1)
\end{cases}-\begin{cases}
X_{k,s}\iota_{\mathcal{A}}^{c_{k}}(Y_{k,s}') & (k,s)\neq(m,n-1)\\
0 & (k,s)=(m,n-1)
\end{cases}
\]
and
\[
\partial_{a_{s},a_{s+1}}(X_{k,s}U(d_{l}-d_{k})\bm{\iota}^{c_{l}}(Y_{l,s}))=0,
\]
we have
\[
\partial_{a_{s},a_{s+1}}(x_{b_{1}}\cdots x_{b_{m}}\dshuffle e_{a_{1}}\cdots e_{a_{n}})=\sum_{k=0}^{m}X_{k,s}'\bm{\iota}^{c_{k}}(Y_{k,s})-\sum_{k=0}^{m}X_{k,s}\bm{\iota}^{c_{k}}(Y_{k,s}')-\delta_{s,1}Y_{0,1}+\delta_{s,n-1}X_{m,n-1}.
\]
By the second formula of Proposition \ref{prop:recur_dshuffle}, 
\begin{align*}
X_{k,s}' & =-\sum_{l=0}^{k-1}(-1)^{k-l}X_{l,s}U(d_{k}-d_{l}-1)\\
 & \ \ +\begin{cases}
X_{k,s-1}+\sum_{l=0}^{k-1}(-1)^{k-l}X_{l,s-1}U(d_{k}-d_{l}) & s>1\\
\delta_{k,0} & s=1,
\end{cases}
\end{align*}
and by the first formula of Proposition \ref{prop:recur_dshuffle},
\begin{align*}
\bm{\iota}^{c_{k}}(Y_{k,s}') & =-\sum_{l=k+1}^{m}(-1)^{l-k}U(d_{l}-d_{k}-1)\bm{\iota}^{c_{l}}(Y_{l,s})\\
 & \ \ +\begin{cases}
\bm{\iota}^{c_{k}}(Y_{k,s+1})+\sum_{l=k+1}^{m}(-1)^{l-k}U(d_{l}-d_{k})\bm{\iota}^{c_{l}}(Y_{l,s+1}) & s+1<n\\
\delta_{k,m} & s+1=n.
\end{cases}
\end{align*}
Thus we have
\[
\sum_{k=0}^{m}X_{k,s}'\bm{\iota}^{c_{k}}(Y_{k,s})=G+F_{1}
\]
and
\[
\sum_{k=0}^{m}X_{k,s}\iota_{\mathcal{A}}^{c_{k}}(Y_{k,s}')=G+F_{2}
\]
where
\[
G=-\sum_{0\leq p<q\leq m}(-1)^{q-p}X_{p,s}U(d_{q}-d_{p}-1)\bm{\iota}^{c_{q}}(Y_{q,s}),
\]
\[
F_{1}=\begin{cases}
\sum_{k=0}^{m}X_{k,s-1}\bm{\iota}^{c_{k}}(Y_{k,s})+\sum_{k=0}^{m}\bigg(\sum_{l=0}^{k-1}(-1)^{k-l}X_{l,s-1}U(d_{k}-d_{l})\bigg)\bm{\iota}^{c_{k}}(Y_{k,s}) & s>1\\
\sum_{k=0}^{m}\delta_{k,0}\bm{\iota}^{c_{k}}(Y_{k,s}) & s=1
\end{cases}
\]
and
\[
F_{2}=\begin{cases}
\sum_{k=0}^{m}X_{k,s}\bm{\iota}^{c_{k}}(Y_{k,s+1})+\sum_{k=0}^{m}X_{k,s}\bigg(\sum_{l=k+1}^{m}(-1)^{l-k}U(d_{l}-d_{k})\bm{\iota}^{c_{l}}(Y_{l,s+1})\bigg) & s+1<n\\
\sum_{k=0}^{m}X_{k,s}\delta_{k,m} & s+1=n.
\end{cases}
\]
Furthermore, by Proposition \ref{prop:dsh_sep_as_as1}, we find
\[
F_{1}=x_{b_{1}}\cdots x_{b_{m}}\dshuffle e_{a_{1}}\cdots\widehat{e_{a_{s}}}\cdots e_{a_{n}}
\]
and
\[
F_{2}=x_{b_{1}}\cdots x_{b_{m}}\dshuffle e_{a_{1}}\cdots\widehat{e_{a_{s+1}}}\cdots e_{a_{n}}.
\]
It follows that
\begin{align*}
 & \partial_{a_{s},a_{s+1}}(x_{b_{1}}\cdots x_{b_{m}}\dshuffle e_{a_{1}}\cdots e_{a_{n}})\\
 & =(G+F_{1})-(G+F_{2})-\delta_{s,1}Y_{0,1}+\delta_{s,n-1}X_{m,n-1}\\
 & =(F_{1}-\delta_{s,1}Y_{0,1})-(F_{2}-\delta_{s,n-1}X_{m,n-1})\\
 & =x_{b_{1}}\cdots x_{b_{m}}\dshuffle\left(\begin{cases}
e_{a_{1}}\cdots\widehat{e_{a_{s}}}\cdots e_{a_{n}} & s>1\\
0 & s=1
\end{cases}\quad-\begin{cases}
e_{a_{1}}\cdots\widehat{e_{a_{s+1}}}\cdots e_{a_{n}} & s<n-1\\
0 & s=n-1
\end{cases}\right)\\
 & =x_{b_{1}}\cdots x_{b_{m}}\dshuffle\left(\partial_{a_{s},a_{s+1}}\left(e_{a_{1}}\cdots e_{a_{n}}\right)\right).
\end{align*}
This completes the proof.
\end{proof}
\begin{lem}
\label{lem:AlgDerFormula_eq}For $\alpha\in S$, we have
\[
\partial_{\alpha,\alpha}(x_{b_{1}}\cdots x_{b_{m}}\dshuffle e_{a_{1}}\cdots e_{a_{n}})=\bigg(-\delta_{b_{1},1}\Delta_{a_{1},a_{1}}^{\alpha,\alpha}x_{b_{2}}\cdots x_{b_{m}}+\delta_{b_{m},1}\Delta_{a_{n},a_{n}}^{\alpha,\alpha}x_{b_{1}}\cdots x_{b_{m-1}}\bigg)\dshuffle e_{a_{1}}\cdots e_{a_{n}}.
\]
\end{lem}

\begin{proof}
Put $c=(b_{1}-1)+\cdots+(b_{m}-1)$ and 
\[
w=x_{b_{1}}\cdots x_{b_{m}}\dshuffle e_{a_{1}}\cdots e_{a_{n}}\in e_{a_{1}}\mathcal{A}\cap\mathcal{A}e_{\iota^{c}(a_{n})}.
\]
By the definition of $\partial_{\alpha,\alpha}$, (since all intermediate
terms cancel out) we have
\[
\partial_{\alpha,\alpha}(w)=-\Delta_{a_{1},a_{1}}^{\alpha,\alpha}w_{1}+\Delta_{a_{n},a_{n}}^{\alpha,\alpha}w_{2}
\]
where $e_{a_{1}}w_{1}$ is the sum of the monomials in $w$ starting
from $e_{a_{1}}e_{a_{1}}$ and $w_{2}e_{\iota^{c}(a_{n})}$ is the
sum of monomials in $w$ ending with $e_{\iota^{c}(a_{n})}e_{\iota^{c}(a_{n})}$.
By Proposition \ref{prop:recur_dshuffle}, such $w_{1}$ and $w_{2}$
are exactly
\begin{align*}
w_{1} & =\delta_{b_{1},1}x_{b_{2}}\cdots x_{b_{m}}\dshuffle e_{a_{1}}\dots e_{a_{n}},\\
w_{2} & =\delta_{b_{m},1}x_{b_{1}}\cdots x_{b_{m-1}}\dshuffle e_{a_{1}}\cdots e_{a_{n}}.
\end{align*}
This proves the claim.
\end{proof}
\begin{lem}
\label{lem:AlgDerFormula_01}We have
\[
\partial_{0,1}(x_{b_{1}}\cdots x_{b_{m}}\dshuffle e_{a_{1}}\cdots e_{a_{n}})=0.
\]
\end{lem}

\begin{proof}
Let $S'\coloneqq S\setminus\{0,1\}=\{a_{1},\dots,a_{n},1-a_{1},\dots,1-a_{n}\}$.
Then, $x_{b_{1}}\cdots x_{b_{m}}\dshuffle e_{a_{1}}\cdots e_{a_{n}}$
is a linear combination of the elements of the form
\[
w_{1}U(t_{1})\cdots w_{k}U(t_{k})w_{k+1}\qquad(k\geq0,w_{1},\dots,w_{k+1}\in\mathcal{A}^{+}(S'),t_{1},\dots,t_{k}\in\mathbb{Z}_{\geq1}).
\]
Now, by definition
\[
\partial_{0,1}(w_{1}U(t_{1})\cdots w_{k}U(t_{k})w_{k+1})=\sum_{i=1}^{k}w_{1}U(t_{1})\cdots w_{i}V(t_{i})w_{i+1}\cdots U(t_{k})w_{k+1}
\]
where
\[
V(t_{i})=\widehat{e_{0}}\overbrace{e_{1}e_{0}e_{1}\cdots}^{t_{i}-1}+\widehat{e_{1}}\overbrace{e_{0}e_{1}e_{0}\cdots}^{t_{i}-1}-\overbrace{\cdots e_{1}e_{0}e_{1}}^{t_{i}-1}\widehat{e_{0}}-\overbrace{\cdots e_{0}e_{1}e_{0}}^{t_{i}-1}\widehat{e_{1}}.
\]
Since $V(t_{i})=0$, it follows that
\[
\partial_{0,1}(x_{b_{1}}\cdots x_{b_{m}}\dshuffle e_{a_{1}}\cdots e_{a_{n}})=0.\qedhere
\]
\end{proof}
\begin{thm}[Algebraic differential formula for $\dshuffle$]
\label{thm:AlgDerFormula}For $\alpha,\beta\in S$, $m\geq0$ and
$b_{1},\dots,b_{m}>0$, we have
\begin{align*}
\partial_{\alpha,\beta}(x_{b_{1}}\cdots x_{b_{m}}\dshuffle e_{a_{1}}\cdots e_{a_{n}}) & =x_{b_{1}}\cdots x_{b_{m}}\dshuffle\partial_{\alpha,\beta}(e_{a_{1}}\cdots e_{a_{n}})\\
 & \ \ \ +(\Delta_{a_{1},0}^{\alpha,\beta}+\Delta_{a_{1},1}^{\alpha,\beta}-\Delta_{a_{1},a_{1}}^{\alpha,\beta})\delta_{b_{1},1}\,x_{b_{2}}\cdots x_{b_{m}}\dshuffle e_{a_{1}}\cdots e_{a_{n}}\\
 & \ \ \ -(\Delta_{a_{n},0}^{\alpha,\beta}+\Delta_{a_{n},1}^{\alpha,\beta}-\Delta_{a_{n},a_{n}}^{\alpha,\beta})\delta_{b_{m},1}\,x_{b_{1}}\cdots x_{b_{m-1}}\dshuffle e_{a_{1}}\cdots e_{a_{n}}.
\end{align*}
In particular, 
\[
\partial_{\alpha,\beta}(x_{\Bbbk}\dshuffle v)=x_{\Bbbk}\dshuffle\partial_{\alpha,\beta}(v)
\]
 for an admissible index $\Bbbk$ and $v\in\mathcal{A}^{+}$.
\end{thm}

\begin{proof}
Putting $S'=\{a_{1},\dots,a_{n},1-a_{1},\dots,1-a_{n}\}$ and $S''=\{0,1\}$,
the claim of the theorem can be divided into the following four cases:
\begin{enumerate}
\item $\alpha\neq\beta$ and
\begin{enumerate}
\item $(\alpha,\beta)\in S'\times S''\cup S''\times S',$
\item $(\alpha,\beta)\in S'\times S',$
\item $(\alpha,\beta)\in S''\times S''$.
\end{enumerate}
\item $\alpha=\beta$.
\end{enumerate}
The cases (1-a), (1-b), (1-c) and (2) follow from Lemmas \ref{prop:AlgDerFormula_a0},
\ref{lem:AlgDerFormula_aa}, \ref{lem:AlgDerFormula_01} and \ref{lem:AlgDerFormula_eq},
respectively.
\end{proof}
\begin{thm}
\label{thm:total_dif_of_dshuffle}Assume $(m,n)\neq(0,1)$. Let $x_{b_{1}}\cdots x_{b_{m}}\in\mathfrak{X}$
and $a_{1},\dots,a_{n}$ be complex variables such that $a_{1},\dots,a_{n},1-a_{1},\dots,1-a_{n},0,1$
are all distinct. Fix $\gamma\in\Pi(\mathbb{C}\setminus S;a_{1},\iota^{b_{1}+\cdots+b_{m}-m}(a_{n}))$.
Then as a function of complex variables $a_{1},\dots,a_{n}$, the
total differential of $L_{\gamma}(x_{b_{1}}\cdots x_{b_{m}}\dshuffle e_{a_{1}}\cdots e_{a_{n}})$
is given by
\begin{align*}
dL_{\gamma}(x_{b_{1}}\cdots x_{b_{m}}\dshuffle e_{a_{1}}\cdots e_{a_{n}}) & =\sum_{i=2}^{n-1}L_{\gamma}(x_{b_{1}}\cdots x_{b_{m}}\dshuffle e_{a_{1}}\cdots\widehat{e_{a_{i}}}\cdots e_{a_{n}})d\log(\frac{a_{i}-a_{i+1}}{a_{i}-a_{i-1}})\\
 & \ \ +\delta_{b_{1},1}L_{\gamma}(x_{b_{2}}\cdots x_{b_{m}}\dshuffle e_{a_{1}}\cdots e_{a_{n}})d\log(a_{1}(1-a_{1}))\\
 & \ \ -\delta_{b_{m},1}L_{\gamma}(x_{b_{1}}\cdots x_{b_{m-1}}\dshuffle e_{a_{1}}\cdots e_{a_{n}})d\log(a_{n}(1-a_{n})).
\end{align*}
Here, we regard the first term of the right-hand side as $0$ for
$n=1,2$. Notice that, when $m=n=1$ and $b_{1}=1$, the second and
the third terms of the right-hand side becomes
\[
L_{\gamma}(e_{a_{1}})d\log(a_{1}(1-a_{1}))-L_{\gamma}(e_{a_{1}})d\log(a_{1}(1-a_{1}))
\]
which we consider as $0$, although it contains an undefined expression
$L_{\gamma}(e_{a_{1}})$.
\end{thm}

\begin{proof}
It follows from \ref{eq:diff_eq_S} and Theorem \ref{thm:AlgDerFormula}.
\end{proof}

\subsection{Analytic multivariable block shuffle identity}
\begin{thm}
\label{thmdiamond_prod_identity_pre}Let $b_{1},\dots,b_{m}$ be positive
integers such that $(b_{1},\dots,b_{m})\neq(1,\dots,1)$, and $a_{1},\dots,a_{n}$
complex numbers such that $a_{1},\dots,a_{n},1-a_{1},\dots,1-a_{n},0,1$
are all distinct. Put $S=\{a_{1},\dots,a_{n},1-a_{1},\dots,1-a_{n},0,1\}$.
Then, for $\gamma\in\Pi(\mathbb{C}\setminus\left((-\infty,0]\cup(1,\infty]\cup S\right);a_{1},\iota^{b_{1}+\cdots+b_{m}-m}(a_{n}))$,
we have
\[
L_{\gamma}(x_{b_{1}}\cdots x_{b_{m}}\dshuffle e_{a_{1}}\cdots e_{a_{n}})=0.
\]
\end{thm}

\begin{proof}
For $0<\alpha\leq1$, put $g_{\alpha}(z)=\frac{1}{2}+\alpha(z-\frac{1}{2})$.
For $0<\alpha\leq1$, consider the function
\[
f(\alpha)=L_{\gamma_{\alpha}}(x_{b_{1}}\cdots x_{b_{m}}\dshuffle e_{g_{\alpha}(a_{1})}\cdots e_{g_{\alpha}(a_{n})})
\]
where $\gamma_{\alpha}$ is a path defined by $\gamma_{\alpha}(t)=g_{\alpha}(\gamma(c(t)))$
where $c(t)$ is $\min(\alpha^{-1}t,1/2)$ is for $t\leq1/2$ and
$\max(1-\alpha^{-1}(1-t),1/2)$ for $t>1/2$, so that $\gamma_{\alpha}'(0)=\gamma'(0)$
and $\gamma_{\alpha}'(1)=\gamma'(1)$. Notice that $\lim_{\alpha\to0}L_{\gamma_{\alpha}}(w_{\alpha})=0$
for any word $w_{\alpha}$ in $e_{0}$, $e_{1}$ and $e_{g_{\alpha}(a)}$'s
with $a\in S\setminus\{0,1\}$, containing at least one $e_{0}$ or
$e_{1}$. By the assumption $(b_{1},\dots,b_{m})\neq(1,\dots,1)$,
every term that appears in $x_{b_{1}}\cdots x_{b_{m}}\dshuffle e_{g_{\alpha}(a_{1})}\cdots e_{g_{\alpha}(a_{n})}$
satisfy this condition, and thus 
\[
\lim_{\alpha\to0}f(\alpha)=0.
\]
Furthermore, we have
\begin{align*}
\frac{d}{d\alpha}f(\alpha) & =\delta_{b_{1},1}L_{\gamma}(x_{b_{2}}\cdots x_{b_{m}}\dshuffle e_{g_{\alpha}(a_{1})}\cdots e_{g_{\alpha}(a_{n})})\frac{d}{d\alpha}\log\left(g_{\alpha}(a_{1})(1-g_{\alpha}(a_{1}))\right)\\
 & \ \ -\delta_{b_{m},1}L_{\gamma}(x_{b_{1}}\cdots x_{b_{m-1}}\dshuffle e_{g_{\alpha}(a_{1})}\cdots e_{g_{\alpha}(a_{n})})\frac{d}{d\alpha}\log\left(g_{\alpha}(a_{n})(1-g_{\alpha}(a_{1}))\right)
\end{align*}
by Theorem \ref{thm:total_dif_of_dshuffle}. Thus by the induction
hypothesis, we have
\[
\frac{d}{d\alpha}f(\alpha)=0.
\]
Therefore $f(\alpha)$ is a constant function, and the constant is
$0$ by letting $\alpha\to0$. Therefore, we conclude that
\[
L_{\gamma}(x_{b_{1}}\cdots x_{b_{m}}\dshuffle e_{a_{1}}\cdots e_{a_{n}})=f(1)=0.\qedhere
\]
\end{proof}
Of course, we can extend the scope of the theorem by the identity
theorem of complex functions. Fix an integer $n\geq1$ and tangential
vectors $u,v\in\mathbb{C}^{\times}$. For $\bullet\in\{{\rm ev},{\rm od}\}$,
let $C_{n}^{\bullet}$ be the certain modified configuration space
consisting of the tuples $(a_{1},\dots,a_{n};\gamma)$ where $a_{1},\dots,a_{n}$
are complex numbers such that $S\coloneqq\{a_{1},\dots,a_{n},1-a_{1},\dots,1-a_{n},0,1\}$
has $2n-2$ distinct elements, and $\gamma$ is an element of $\pi_{1}(\mathbb{C}\setminus S,s,t)$
where $s$ and $t$ are the tangential base points
\[
s=\overrightarrow{u}_{a_{1}},\ t=\begin{cases}
\overrightarrow{v}_{a_{n}} & \bullet={\rm ev}\\
\overrightarrow{v}_{\iota(a_{n})} & \bullet={\rm od}.
\end{cases}
\]
Then we can regard $C_{n}^{\bullet}$ as a complex manifold, and for
any $x_{b_{1}}\cdots x_{b_{m}}\in\mathfrak{X}_{\bullet}$, we can
view $L_{\gamma}(x_{b_{1}}\cdots x_{b_{m}}\dshuffle e_{a_{1}}\cdots e_{a_{n}})$
as a holomorphic function on $C_{n}^{\bullet}$. We define a subspace
$C_{n}^{\bullet,{\rm good}}\subset C_{n}^{\bullet}$ consisting of
$(a_{1}',\dots,a_{n}';\gamma')$ such that there exists $(a_{1},\dots,a_{n};\gamma)$
belonging to the same connected component as $(a_{1}',\dots,a_{n}';\gamma)$
and satisfying the same conditions as in Theorem \ref{thmdiamond_prod_identity_pre}.
\begin{example}
Let $(a_{1},\dots,a_{n};\gamma)\in C_{n}^{\bullet}$. If $n\geq2$,
$a_{1},\dots,a_{n}\in\mathbb{R}$, and $\gamma$ is a path on $\{z\mid{\rm Im}(z)\geq0\}$
then $(a_{1},\dots,a_{n};\gamma)\in C_{n}^{\bullet,{\rm good}}$ by
the following reason. Notice that $\gamma$ starts from $a_{1}$ and
ends at $a_{n}$ or $\iota(a_{n})$ according to the parity. In either
case, we can slightly move the entire path slightly upward so that
the resulting path satisfies the condition of Theorem \ref{thmdiamond_prod_identity_pre}.
\end{example}

Now we have:
\begin{thm}[Analytic multivariable block shuffle identity]
\label{thm:Diamond}Let $\bullet\in\{{\rm od},{\rm ev}\}$ and $x_{b_{1}}\cdots x_{b_{m}}\in\mathfrak{X}_{\bullet}$
with $(b_{1},\dots,b_{m})\neq(1,\dots,1)$. Then for $(a_{1},\dots,a_{n};\gamma)\in C_{n}^{\bullet,{\rm good}}$,
we have
\[
L_{\gamma}(x_{b_{1}}\cdots x_{b_{m}}\dshuffle e_{a_{1}}\cdots e_{a_{n}})=0.
\]
\end{thm}

\begin{rem}
From the perspective of Grothendieck's period conjecture, the motivic
version of the theorem with $a_{1},\ldots,a_{n}\in\overline{\mathbb{Q}}$
is presumably true. Also, we expect that the case with $a_{1},\ldots,a_{n}\in\mathbb{Q}$
should follow from the theory of motivic confluence relations discussed
in another article \cite{ConfMotivicityArxiv} by the authors. However,
the main theme of this article is to prove the original block shuffle
relation, and so we will not discuss the detail of the multivariable
block shuffle relation in the motivic setting.
\end{rem}

\section{Proof of the block shuffle identity\label{sec:proof_of_mshuffle}}

\subsection{Proof of the first form of the block shuffle identity}

The purpose of this section is to prove the first form of block shuffle
identity
\[
L_{B}^{\shuffle}(u\mshuffle v)=0
\]
(Theorem \ref{thm:BlockShuffleIdentity}). The real-valued version
of the block shuffle identity can also be regarded as a special instance
of the (analytic) multivariable block shuffle identity stated in the
previous section (Theorem \ref{thm:Diamond}), since $L_{B}^{\shuffle}(u\mshuffle v)=L_{\mathrm{dch}}(u\dshuffle B(v))$
with $\mathrm{dch}(t)=t$ by Theorem \ref{thm:block_diamond_compatibility}
and $L_{\gamma}(u\dshuffle w)=0$ by Theorem \ref{thm:Diamond}. Strictly
speaking, the case $\gamma=\mathrm{dch}$ and $w=B(v)$ itself is
not contained in Theorem \ref{thm:Diamond} since it does not satisfy
the assumption of the theorem, but such cases can be derived from
the generic statement by a suitable limiting argument. In this section,
however, we specify ourselves to the motivic setting and will not
discuss the real-valued version. Our proof below is based on the motivicity
of general confluence relations proved in \cite{ConfMotivicityArxiv}.

Let $\mathcal{H}$ be the ring of motivic periods of mixed Tate motives
over $\mathbb{Q}$. Fix a real number $z$ between $0$ and $1/2$.
We put $\mathcal{A}_{z}=e_{z}\mathcal{A}(\{0,1,z,1-z\})e_{1-z}$.
For $c\in\{0,1,1/2\}$, define a $\mathbb{Q}$-linear map $\partial_{c}:\mathcal{A}_{z}\to\mathcal{A}_{z}$
by
\[
\partial_{c}(e_{a_{0}}e_{a_{1}}\cdots e_{a_{k}}e_{a_{k+1}})=\sum_{r\in\{\pm1\}}r\sum_{i=1}^{k}{\rm ord}_{z-c}(a_{i}-a_{i+r})e_{a_{0}}\cdots\widehat{e_{a_{i}}}\cdots e_{a_{k+1}}.
\]
Here we understand ${\rm ord}_{z-c}(0)$ as $0$. The operators $\partial_{0}$,
$\partial_{1}$ and $\partial_{1/2}$ are essentially equal to $\partial_{z,0}$,
$\partial_{z,1}$ and $\partial_{z,1-z}$ defined in Section \ref{subsec:Algebraic-differential-formula}.
Now we define ${\rm ev}_{0}^{\mathfrak{m}},{\rm ev}_{1/2}^{\mathfrak{m}}:\mathcal{A}(\{0,1,z,1-z\})\to\mathcal{H}$.
We remark that only the properties (1), (2), (3) of Definition \ref{def:ev_m_0}
and (1), (3) of Definition \ref{def:ev_m_0} will be used in our proof.
Hereafter, we mean by $I^{\mathfrak{m}}(a_{0},a_{1},\dots a_{k},a_{k+1})$
the motivic iterated integral along the straight line path from $a_{0}$
to $a_{k+1}$.
\begin{defn}
\label{def:ev_m_0}Define ${\rm ev}_{0}^{\mathfrak{m}}:\mathcal{A}(\{0,1,z,1-z\})\to\mathcal{H}$
as the unique $\mathbb{Q}$-linear map satisfying the following conditions:
\begin{enumerate}
\item ${\rm ev}_{0}^{\mathfrak{m}}(1)=1.$
\item For $u,v\in\mathcal{A}(\{0,1,z,1-z\})$, we have
\[
{\rm ev}_{0}^{\mathfrak{m}}(u\shuffle v)={\rm ev}_{0}^{\mathfrak{m}}(u)\cdot{\rm ev}_{0}^{\mathfrak{m}}(v).
\]
\item Assume $k>0$. If $a_{1}\notin\{0,z\}$ and $a_{k}\notin\{1,1-z\}$
then
\[
{\rm ev}_{0}^{\mathfrak{m}}(e_{a_{1}}\cdots e_{a_{k}})=I^{\mathfrak{m}}(0,a_{1}',\dots a_{k}',1)
\]
where $a_{j}'=\left.a_{j}\right|_{z=0}$ for $j=1,\dots,k$.
\item Assume $k>0$. If $a_{1},\dots,a_{k}\in\{0,z\}$ and $a_{1}\neq z$.
Then
\[
{\rm ev}_{0}^{\mathfrak{m}}(e_{a_{1}}\cdots e_{a_{k}})=I^{\mathfrak{m}}(1,a_{1}',\dots,a_{k}',\infty')
\]
where $a_{j}'=\delta_{a_{j},z}$ for $j=1,\dots,k$ and $\infty'$
is the tangential base point at $+\infty$ such that $I^{\mathfrak{m}}(1,0,\infty')=0$.
\item ${\rm ev}_{0}^{\mathfrak{m}}(e_{z})=0$.
\item Assume $k>0$. If $a_{1},\dots,a_{k}\in\{1,1-z\}$ and $a_{k}\neq1-z$.
Then
\[
{\rm ev}_{0}^{\mathfrak{m}}(e_{a_{1}}\cdots e_{a_{k}})=I^{\mathfrak{m}}(\infty',a_{1}',\dots,a_{k}',1)
\]
where $a_{j}'=\delta_{a_{j},1-z}$ for $j=1,\dots,k$ and $\infty'$
is the tangential base point at $+\infty$ such that $I^{\mathfrak{m}}(1,0,\infty')=0$.
\item ${\rm ev}_{0}^{\mathfrak{m}}(e_{1-z})=0$.
\end{enumerate}
\end{defn}

\begin{defn}
Define ${\rm ev}_{1/2}^{\mathfrak{m}}:\mathcal{A}(\{0,1,z,1-z\})\to\mathcal{H}$
as the unique $\mathbb{Q}$-linear map satisfying the following conditions:
\begin{enumerate}
\item ${\rm ev}_{1/2}^{\mathfrak{m}}(1)=1.$
\item For $u,v\in\mathcal{A}(\{0,1,z,1-z\})$, we have
\[
{\rm ev}_{1/2}^{\mathfrak{m}}(u\shuffle v)={\rm ev}_{1/2}^{\mathfrak{m}}(u)\cdot{\rm ev}_{1/2}^{\mathfrak{m}}(v).
\]
\item If $0\in\{a_{1},\dots,a_{k}\}$ or $1\in\{a_{1},\dots,a_{k}\}$ then
\[
{\rm ev}_{1/2}^{\mathfrak{m}}(e_{a_{1}}\cdots e_{a_{k}})=0.
\]
\item Assume $k>0$. If $a_{1}\neq z$, $a_{k}\neq1-z$, and $\{a_{1},\dots,a_{k}\}\subset\{z,1-z\}$,
then
\[
{\rm ev}_{1/2}^{\mathfrak{m}}(e_{a_{1}}\cdots e_{a_{k}})=I^{\mathfrak{m}}(0,a_{1}',\dots,a_{k}',1)
\]
where $a_{j}'=\left.a_{j}\right|_{z=0}$ for $j=1,\dots,k$.
\item ${\rm ev}_{1/2}^{\mathfrak{m}}(e_{z})=\log^{\mathfrak{m}}(2)\coloneqq I^{\mathfrak{m}}(1,0,2)$.
\item ${\rm ev}_{1/2}^{\mathfrak{m}}(e_{1-z})=-\log^{\mathfrak{m}}(2)$.
\end{enumerate}
\end{defn}

\begin{rem}
The definitions of ${\rm ev}_{x}^{\mathfrak{m}}$ ($x\in\{0,1/2\}$)
are compatible with the definition of ${\rm ev}_{x}^{\mathfrak{m}}$
in \cite{ConfMotivicityArxiv} in the sense that
\[
{\rm ev}_{x}^{\mathfrak{m}}(e_{a_{1}}\cdots e_{a_{k}})=\widetilde{{\rm ev}_{x}^{\mathfrak{m}}}\left(\mathbb{I}_{{\rm dch}}(\overrightarrow{1}_{z};a_{1},\dots,a_{k};-\overrightarrow{1}_{1-z})\right),
\]
where ${\rm dch}$ is a linear path from $\overrightarrow{1}_{z}$
to $-\overrightarrow{1}_{1-z}$ (here we write $\widetilde{{\rm ev}_{x}^{\mathfrak{m}}}$
for ${\rm ev}_{x}^{\mathfrak{m}}$ in \cite{ConfMotivicityArxiv}
to avoid confusion), which guarantees the well-definedness of ${\rm ev}_{x}^{\mathfrak{m}}$.
\end{rem}

The value of ${\rm ev}_{0}^{\mathfrak{m}}(w)$ for a general $w\in\mathcal{A}(\{0,1,z,1-z\})$
is decomposed as follows. We put
\[
\mathcal{A}^{0}\coloneqq\mathbb{Q}\oplus\bigoplus_{\substack{a\in\{1,1-z\}\\
b\in\{0,z\}
}
}e_{a}\mathcal{A}(\{0,z,1,1-z\})e_{b}.
\]
Let ${\rm reg}$ be the composite map 
\[
\mathcal{A}(\{0,1,z,1-z\})\to\mathcal{A}(\{0,z\})\otimes\mathcal{A}^{0}\otimes\mathcal{A}(\{1,1-z\})\xrightarrow{c\otimes{\rm id}\otimes c}\mathcal{A}^{0}
\]
where the first map is inverse of the isomorphism
\[
\mathcal{A}(\{0,z\})\otimes\mathcal{A}^{0}\otimes\mathcal{A}(\{1,1-z\})\to\mathcal{A}(\{0,1,z,1-z\})\ ;\ w_{1}\otimes w_{2}\otimes w_{3}\mapsto w_{1}\shuffle w_{2}\shuffle w_{3},
\]
and $c$ is a $\mathbb{Q}$-linear map defined by $c(e_{a_{1}}\cdots e_{a_{n}})=\delta_{n,0}$.
Then by \cite[(3.2.20)]{PanzerArxiv} we have an explicit decomposition
\begin{align*}
w & =\sum_{n,m=0}^{\infty}\sum_{a_{1},\dots,a_{n}\in\{0,z\}}\sum_{b_{1},\dots,b_{m}\in\{1,1-z\}}e_{a_{1}}\cdots e_{a_{n}}\shuffle e_{b_{1}}\cdots e_{b_{m}}\shuffle{\rm reg}({\rm Cut}_{(a_{1},\dots,a_{n};b_{1},\dots,b_{m})}(w))
\end{align*}
where $w\mapsto{\rm Cut}_{(a_{1},\dots,a_{n};b_{1},\dots,b_{m})}(w)$
is a $\mathbb{Q}$-linear map from $\mathcal{A}(\{0,1,z,1-z\})$ to
$\mathcal{A}(\{0,1,z,1-z\})$ defined by
\[
{\rm Cut}_{(a_{1},\dots,a_{n};b_{1},\dots,b_{m})}(e_{s_{1}}\cdots e_{s_{t}})=\begin{cases}
w' & e_{s_{1}}\cdots e_{s_{t}}\ \text{can be written as }e_{a_{1}}\cdots e_{a_{n}}w'e_{b_{1}}\cdots e_{b_{m}}\\
{\rm otherwise} & 0.
\end{cases}
\]
Thus we have
\begin{align}
{\rm ev}_{0}^{\mathfrak{m}}(w) & =\sum_{n,m=0}^{\infty}\sum_{a_{1},\dots,a_{n}\in\{0,z\}}\sum_{b_{1},\dots,b_{m}\in\{1,1-z\}}{\rm ev}_{0}^{\mathfrak{m}}(e_{a_{1}}\cdots e_{a_{n}}){\rm ev}_{0}^{\mathfrak{m}}(e_{b_{1}}\cdots e_{b_{m}}){\rm ev}_{0}^{\mathfrak{m}}({\rm reg}({\rm Cut}_{(a_{1},\dots,a_{n};b_{1},\dots,b_{m})}(w)))\nonumber \\
 & =\sum_{n,m=0}^{\infty}\sum_{a_{1},\dots,a_{n}\in\{0,z\}}\sum_{b_{1},\dots,b_{m}\in\{1,1-z\}}{\rm ev}_{0}^{\mathfrak{m}}(e_{a_{1}}\cdots e_{a_{n}}){\rm ev}_{0}^{\mathfrak{m}}(e_{b_{1}}\cdots e_{b_{m}})\tilde{L}({\rm reg}({\rm Cut}_{(a_{1},\dots,a_{n};b_{1},\dots,b_{m})}(w))\mid_{z\to0})\nonumber \\
 & =\sum_{n,m=0}^{\infty}\sum_{a_{1},\dots,a_{n}\in\{0,z\}}\sum_{b_{1},\dots,b_{m}\in\{1,1-z\}}{\rm ev}_{0}^{\mathfrak{m}}(e_{a_{1}}\cdots e_{a_{n}}){\rm ev}_{0}^{\mathfrak{m}}(e_{b_{1}}\cdots e_{b_{m}})\tilde{L}({\rm Cut}_{(a_{1},\dots,a_{n};b_{1},\dots,b_{m})}(w)\mid_{z\to0})\label{eq:evm0_dec}
\end{align}
where $\tilde{L}:\mathcal{A}(\{0,1\})\to\mathcal{H}$ is a linear
map defined by $\tilde{L}(e_{a_{1}}\cdots e_{a_{n}})\coloneqq I^{\mathfrak{m}}(0,a_{1},\dots,a_{n},1)$
and $u\mapsto u\mid_{z\to0}$ is a ring homomorphism from $\mathcal{A}(\{0,1,z,1-z\})$
to $\mathcal{A}(\{0,1\})$ defined by $e_{0},e_{z}\mapsto e_{0}$
and $e_{1},e_{1-z}\mapsto e_{1}$. We define a $\mathbb{Q}$-linear
map $\llbracket-\rrbracket:\mathcal{A}(\{0,1,z,1-z\})\to\mathcal{A}(\{0,1,z,1-z\})$
by $\llbracket e_{a_{1}}\cdots e_{a_{k}}\rrbracket\coloneqq e_{a_{2}}\cdots e_{a_{k-1}}$
for $k\geq2$ and $0$ for $k=0,1$.
\begin{lem}
\label{lem:ev0(u_dsh_v)=00003D0}For $u\in\mathfrak{X}_{{\rm ev}}'$
and $v\in\mathcal{A}_{z}$, we have
\[
{\rm ev}_{0}^{\mathfrak{m}}\Big(\llbracket u\dshuffle v\rrbracket\Big)=0.
\]
\end{lem}

\begin{proof}
By the motivicity of generalized confluence relation (\cite[Theorem 6.21]{ConfMotivicityArxiv}),
we have, 
\begin{equation}
{\rm ev}_{0}^{\mathfrak{m}}(\llbracket w\rrbracket)=\sum_{l=0}^{\infty}\sum_{c_{1},\dots,c_{l}\in\{0,1,1/2\}}{\rm ev}_{1/2}^{\mathfrak{m}}(\llbracket\partial_{c_{1}}\cdots\partial_{c_{l}}w\rrbracket)\cdot I^{\mathfrak{m}}(1/2;c_{1},\dots,c_{l};0)\label{eq:motivc_conf}
\end{equation}
for $w\in\mathcal{A}_{z}$. By applying (\ref{eq:motivc_conf}) to
the case $w=u\dshuffle v$, we have
\[
{\rm ev}_{0}^{\mathfrak{m}}(\llbracket u\dshuffle v\rrbracket)=\sum_{l=0}^{\infty}\sum_{c_{1},\dots,c_{l}\in\{0,1,1/2\}}{\rm ev}_{1/2}^{\mathfrak{m}}(\llbracket\partial_{c_{1}}\cdots\partial_{c_{l}}(u\dshuffle v)\rrbracket)\cdot I^{\mathfrak{m}}(1/2;c_{1},\dots,c_{l};0).
\]
By the algebraic differential formula (Theorem \ref{thm:AlgDerFormula})
, we have
\[
\partial_{c}(\mathfrak{X}_{{\rm ev}}'\dshuffle\mathcal{A}_{z})\subset\mathfrak{X}_{{\rm ev}}'\dshuffle\mathcal{A}_{z}
\]
for $c\in\{0,1,1/2\}$, and especially we have
\[
\partial_{c_{1}}\cdots\partial_{c_{l}}(u\dshuffle v)\in\mathfrak{X}_{{\rm ev}}'\dshuffle\mathcal{A}_{z},
\]
and thus
\[
{\rm ev}_{1/2}^{\mathfrak{m}}(\llbracket\partial_{c_{1}}\cdots\partial_{c_{l}}(u\dshuffle v)\rrbracket)=0
\]
for all $l\geq0$, $c_{1},\dots,c_{l}\in\{0,1,1/2\}$. Hence,
\[
{\rm ev}_{0}^{\mathfrak{m}}(\llbracket u\dshuffle v\rrbracket)=0.\qedhere
\]
\end{proof}
We put $\mathscr{X}\coloneqq\mathbb{Q}\langle\langle X_{1},X_{2},\dots\rangle\rangle$
and
\begin{align*}
\psi & \coloneqq\sum_{m,n=1}^{\infty}\sum_{\substack{b_{1},\dots,b_{m},c_{1},\dots,c_{n}\geq1\\
b_{1}+\cdots+b_{m}-m:\,\mathrm{even}_{>0}\\
c_{1}+\cdots+c_{n}-n:\,\mathrm{odd}_{>0}
}
}X_{b_{1}}\cdots X_{b_{m}}\otimes X_{c_{1}}\cdots X_{c_{n}}\otimes\llbracket x_{b_{1}}\cdots x_{b_{m}}\dshuffle B_{z}(x_{c_{1}}\cdots x_{c_{n}})\rrbracket\\
 & \quad\in\mathscr{X}\hat{\otimes}\mathscr{X}\hat{\otimes}\mathcal{A}(\{0,z,1,1-z\})
\end{align*}
where $B_{z}(x_{c_{1}}\cdots x_{c_{n}})$ is the unique element of
$\mathcal{A}(\{z,1-z\})$ such that $B_{z}(x_{c_{1}}\cdots x_{c_{n}})\mid_{z\to0}=B(x_{c_{1}}\cdots x_{c_{m}})$.
Given a $\mathbb{Q}$-algebra $R$ and $\mathbb{Q}$-linear map $f:\mathcal{A}(\{0,z,1,1-z\})\rightarrow R$,
we also denote by $f$ its natural extension to $\mathscr{X}\hat{\otimes}\mathscr{X}\hat{\otimes}\mathcal{A}(\{0,z,1,1-z\})\to\mathscr{X}\hat{\otimes}\mathscr{X}\hat{\otimes}R$
i.e., $f(W\otimes W'\otimes w)\coloneqq W\otimes W'\otimes f(w)$
by an abuse of notation. Let $u=X_{1}\otimes1\otimes1$ and $v=1\otimes X_{1}\otimes1$
be two elements of $\mathscr{X}\hat{\otimes}\mathscr{X}\hat{\otimes}\mathcal{A}(\{0,z,1,1-z\})$.
\begin{lem}
\label{lem:eqs_CutPsi}We have
\begin{align}
{\rm Cut}_{(0,0;\emptyset)}(\psi) & =0,\label{eq:CutPsi(00;)}\\
{\rm Cut}_{(\emptyset;1,1)}(\psi) & =0,\label{eq:CutPsi(;11)}\\
{\rm Cut}_{(z;\emptyset)}(\psi) & =(2u+v)\psi,\label{eq:CutPsi(z;)}\\
{\rm Cut}_{(\emptyset;1-z)}(\psi) & =\psi(2u+v),\label{eq:CutPsi(;1-z)}\\
{\rm Cut}_{(0,z;\emptyset)}(\psi) & =-uv\psi,\label{eq:CutPsi(0z;)}\\
{\rm Cut}_{(\emptyset;1-z,1)}(\psi) & =-\psi uv,\label{eq:CutPsi(;1-z,1)}\\
{\rm Cut}_{(0;\emptyset)}(\psi)\mid_{z\to0} & =-u\cdot(\psi\mid_{z\to0}),\label{eq:CutPsi(0;)}\\
{\rm Cut}_{(\emptyset;1)}(\psi)\mid_{z\to0} & =-(\psi\mid_{z\to0})\cdot u,\label{eq:CutPsi(;1)}\\
{\rm Cut}_{(0;1)}(\psi)\mid_{z\to0} & =u\cdot(\psi\mid_{z\to0})\cdot u.\label{eq:CutPsi(0;1)}
\end{align}
\end{lem}

\begin{proof}
Notice that the terms that appear in $x_{b_{1}}\cdots x_{b_{m}}\dshuffle B(x_{c_{1}}\cdots x_{c_{n}})$
are words consisting of $e_{z},e_{1-z}$ and alternating sequences
of $\{e_{0},e_{1}\}$ being separated by $e_{z}$ or $e_{1-z}$, thus
contain no subsequences of the forms $e_{0}e_{0}$ or $e_{1}e_{1}$.
Hence (\ref{eq:CutPsi(00;)}) and (\ref{eq:CutPsi(;11)}) holds. Writing
$B_{z}(x_{c_{1}}\cdots x_{c_{n}})$ as $e_{a_{1}}\cdots e_{a_{N}}$
($a_{1}=z,\,a_{N}=1-z$), the first formula of Proposition \ref{prop:recur_dshuffle}
says that the terms of $x_{b_{1}}\cdots x_{b_{m}}\dshuffle e_{a_{1}}\cdots e_{a_{N}}$
starting as $e_{z}e_{z}\cdots$ only occur as 
\[
(-1)^{k-1}e_{a_{1}}U(b_{1}+\cdots+b_{k}-1)\cdot\bm{\iota}^{b_{1}+\cdots+b_{k}-k}(x_{b_{k+1}}\cdots x_{b_{m}}\dshuffle e_{a_{1}}\cdots e_{a_{N}})\quad\left(k=b_{1}=1\right)
\]
or as 
\[
e_{a_{1}}(x_{b_{1}}\cdots x_{b_{m}}\dshuffle e_{a_{2}}\cdots e_{a_{N}})\quad\left(a_{2}=z\right).
\]
From this observation, we find (\ref{eq:CutPsi(z;)}). Similarly,
the first formula of Proposition \ref{prop:recur_dshuffle} also says
that the terms of $x_{b_{1}}\cdots x_{b_{m}}\dshuffle e_{a_{1}}\cdots e_{a_{N}}$
starting as $e_{z}e_{0}e_{z}\cdots$ only occur as 
\[
\sum_{k=1}^{m}(-1)^{k}e_{a_{1}}W(0;b_{1}+\cdots+b_{k})\cdot\bm{\iota}^{b_{1}+\cdots+b_{k}-k}(x_{b_{k+1}}\cdots x_{b_{m}}\dshuffle e_{a_{2}}\cdots e_{a_{n}})\quad\left(k=b_{1}=1,\ a_{2}=z\right),
\]
from which we obtain (\ref{eq:CutPsi(0z;)}). By the same argument
for the tail of $x_{b_{1}}\cdots x_{b_{m}}\dshuffle e_{a_{1}}\cdots e_{a_{N}}$
using the second formula of Proposition \ref{prop:recur_dshuffle},
we get (\ref{eq:CutPsi(;1-z)}) and (\ref{eq:CutPsi(;1-z,1)}). Note
that 
\[
\psi\mid_{z\to0}=\sum_{m,n=1}^{\infty}\sum_{\substack{b_{1},\dots,b_{m},c_{1},\dots,c_{n}\geq1\\
b_{1}+\cdots+b_{m}-m:\,\mathrm{even}_{>0}\\
c_{1}+\cdots+c_{n}-n:\,\mathrm{odd}_{>0}
}
}X_{b_{1}}\cdots X_{b_{m}}\otimes X_{c_{1}}\cdots X_{c_{n}}\otimes\left\llbracket B(x_{b_{1}}\cdots x_{b_{m}}\mshuffle x_{c_{1}}\cdots x_{c_{n}})\right\rrbracket 
\]
by \ref{thm:block_diamond_compatibility}. Since
\[
x_{b_{1}}\cdots x_{b_{m}}\mshuffle x_{c_{1}}\cdots x_{c_{n}}=x_{b_{1}}(x_{b_{2}}\cdots x_{b_{m}}\mshuffle x_{c_{1}}\cdots x_{c_{n}})+x_{c_{1}}(x_{b_{1}}\cdots x_{b_{m}}\mshuffle x_{c_{2}}\cdots x_{c_{n}})-s_{b_{1}+c_{1}}(x_{b_{2}}\cdots x_{b_{m}}\mshuffle x_{c_{2}}\cdots x_{c_{n}}),
\]
the terms of $B(x_{b_{1}}\cdots x_{b_{m}}\mshuffle x_{c_{1}}\cdots x_{c_{n}})$
which start as $e_{0}e_{0}\cdots$ only occur as
\[
B(x_{b_{1}}(x_{b_{2}}\cdots x_{b_{m}}\mshuffle x_{c_{1}}\cdots x_{c_{n}}))\quad(b_{1}=1)
\]
or
\[
B(x_{c_{1}}(x_{b_{1}}\cdots x_{b_{m}}\mshuffle x_{c_{2}}\cdots x_{c_{n}}))\quad(c_{1}=1),
\]
and thus
\begin{equation}
{\rm Cut}_{(0;\emptyset)}(\psi\mid_{z\to0})=(u+v)\cdot(\psi\mid_{z\to0}).\label{eq:Cut_Psi0(0;)}
\end{equation}
Similarly, we also have
\begin{equation}
{\rm Cut}_{(\emptyset;1)}(\psi\mid_{z\to0})=(\psi\mid_{z\to0})\cdot(u+v).\label{eq:Cut_Psi0(;1)}
\end{equation}
Then (\ref{eq:CutPsi(0;)}) follows from the calculation
\begin{align*}
{\rm Cut}_{(0;\emptyset)}(\psi)\mid_{z\to0} & ={\rm Cut}_{(0;\emptyset)}(\psi\mid_{z\to0})-{\rm Cut}_{(z;\emptyset)}(\psi)\mid_{z\to0}\\
 & =(u+v)\cdot(\psi\mid_{z\to0})-(2u+v)\cdot(\psi\mid_{z\to0})\quad\quad(\text{using (\ref{eq:Cut_Psi0(0;)}) and (\ref{eq:CutPsi(z;)})})\\
 & =-u\cdot(\psi\mid_{z\to0}),
\end{align*}
and (\ref{eq:CutPsi(;1)}) also follows from the similar calculation.
Finally, (\ref{eq:CutPsi(0;1)}) follows from the following calculation:
\begin{align*}
{\rm Cut}_{(0;1)}(\psi)\mid_{z\to0} & ={\rm Cut}_{(0;1)}(\psi\mid_{z\to0})-{\rm Cut}_{(z;1-z)}(\psi)\mid_{z\to0}-{\rm Cut}_{(z;1)}(\psi)\mid_{z\to0}-{\rm Cut}_{(0;1-z)}(\psi)\mid_{z\to0}\\
 & =(u+v)\cdot(\psi\mid_{z\to0})\cdot(u+v)-(2u+v)\cdot(\psi\mid_{z\to0})\cdot(2u+v)\\
 & \quad\quad-{\rm Cut}_{(\emptyset;1)}((2u+v)\psi)\mid_{z\to0}-{\rm Cut}_{(0;\emptyset)}(\psi(2u+v))\mid_{z\to0}\\
 & \quad(\text{using (\ref{eq:Cut_Psi0(0;)}), (\ref{eq:Cut_Psi0(;1)}), (\ref{eq:CutPsi(z;)}), and (\ref{eq:CutPsi(;1-z)})})\\
 & =(u+v)\cdot(\psi\mid_{z\to0})\cdot(u+v)-(2u+v)\cdot(\psi\mid_{z\to0})\cdot(2u+v)\\
 & \quad\quad+(2u+v)\cdot(\psi\mid_{z\to0})\cdot u+u\cdot(\psi\mid_{z\to0})\cdot(2u+v)\\
 & \quad(\text{using (\ref{eq:CutPsi(;1)}), and (\ref{eq:CutPsi(0;)})})\\
 & =u\cdot(\psi\mid_{z\to0})\cdot u.\qedhere
\end{align*}
\end{proof}
Now we are ready to prove the first form of the block shuffle relation.
\begin{proof}[Proof of Theorem \ref{thm:BlockShuffleIdentity}]
By Theorem \ref{thm:block_diamond_compatibility}, it is enough to
show $\tilde{L}\left(\psi\mid_{z\to0}\right)=0$. By (\ref{eq:evm0_dec}),
we have
\begin{align*}
{\rm ev}_{0}^{\mathfrak{m}}(\psi) & =\sum_{n,m=0}^{\infty}\sum_{\substack{a_{1},\dots,a_{n}\in\{0,z\}\\
b_{1},\dots,b_{m}\in\{1,1-z\}
}
}\big(1\otimes1\otimes{\rm ev}_{0}^{\mathfrak{m}}(e_{a_{1}}\cdots e_{a_{n}}){\rm ev}_{0}^{\mathfrak{m}}(e_{b_{1}}\cdots e_{b_{m}})\big)\cdot\tilde{L}\big({\rm Cut}_{(a_{1},\dots,a_{n};b_{1},\dots,b_{m})}(\psi)\mid_{z\to0}\big)
\end{align*}
For sequence ${\bf a}=(a_{1},\dots,a_{n})\in\{0,z\}^{n}$, define
$p_{{\bf a}}\in\mathbb{Q}[u,v]$ by
\[
p_{{\bf a}}=\begin{cases}
1 & {\bf a}=\emptyset\\
-u & {\bf a}=(0)\\
-uv\cdot p_{{\bf a}'} & {\bf a}=(0,z,{\bf a}')\\
0 & {\bf a}=(0,0,{\bf a}')\\
(2u+v)\cdot p_{{\bf a}'} & {\bf a}=(z,{\bf a}')
\end{cases}
\]
and for sequence ${\bf b}=(b_{1},\dots,b_{m})\in\{1,1-z\}^{m}$, define
$q_{{\bf b}}$ by
\[
q_{{\bf b}}=\begin{cases}
1 & {\bf b}=\emptyset\\
-u & {\bf b}=(1)\\
-uv\cdot q_{{\bf b}'} & {\bf b}=({\bf b}',1-z,1)\\
0 & {\bf b}=({\bf b}',1,1)\\
(2u+v)\cdot q_{{\bf b}'} & {\bf b}=({\bf b}',1-z).
\end{cases}
\]
Then by Lemma \ref{lem:eqs_CutPsi}, for any , ${\bf a}\in\{0,z\}^{n}$
and ${\bf b}\in\{1,1-z\}^{m}$, we have
\[
{\rm Cut}_{({\bf a};{\bf b})}(\psi)\mid_{z\to0}=p_{{\bf a}}\cdot(\psi\mid_{z\to0})\cdot q_{{\bf b}}.
\]
Thus
\begin{align*}
{\rm ev}_{0}^{\mathfrak{m}}(\psi) & =\sum_{n,m=0}^{\infty}\sum_{\substack{{\bf a}=(a_{1},\dots,a_{n})\in\{0,z\}^{n}\\
{\bf b}=(b_{1},\dots,b_{m})\in\{1,1-z\}^{m}
}
}\big(1\otimes1\otimes{\rm ev}_{0}^{\mathfrak{m}}(e_{a_{1}}\cdots e_{a_{n}}){\rm ev}_{0}^{\mathfrak{m}}(e_{b_{1}}\cdots e_{b_{m}})\big)\cdot\tilde{L}(p_{{\bf a}}\cdot(\psi\mid_{z\to0})\cdot q_{{\bf b}})\\
 & =\Theta_{L}\cdot\tilde{L}\left(\psi\mid_{z\to0}\right)\cdot\Theta_{R}
\end{align*}
where
\begin{align*}
\Theta_{L} & =\sum_{n=0}^{\infty}\sum_{{\bf a}=(a_{1},\dots,a_{n})\in\{0,z\}^{n}}(1\otimes1\otimes{\rm ev}_{0}^{\mathfrak{m}}(e_{a_{1}}\cdots e_{a_{n}}))\cdot p_{{\bf a}},\\
\Theta_{R} & =\sum_{m=0}^{\infty}\sum_{{\bf b}=(b_{1},\dots,b_{m})\in\{1,1-z\}^{m}}(1\otimes1\otimes{\rm ev}_{0}^{\mathfrak{m}}(e_{b_{1}}\cdots e_{b_{m}}))\cdot q_{{\bf b}}.
\end{align*}
Here, by Lemma \ref{lem:ev0(u_dsh_v)=00003D0} we have
\[
{\rm ev}_{0}^{\mathfrak{m}}(\psi)=0.
\]
Since the constant terms of
\[
\Theta_{L},\Theta_{R}\in\mathcal{H}[[u,v]]\quad\quad(\subset\mathscr{X}\hat{\otimes}\mathscr{X}\hat{\otimes}\mathcal{H})
\]
are both $1$, $\Theta_{L}$ and $\Theta_{R}$ are invertible. Hence,
we find
\[
\tilde{L}\left(\psi\mid_{z\to0}\right)=0.\qedhere
\]
\end{proof}

\subsection{Proof of the second form of the block shuffle identity}

The purpose of this section is to prove the second form of the block
shuffle identity (Theorem \ref{thm:BlockShuffleIdentity2}). The second
form of the block shuffle identity is a statement about the block
regularized value, which is defined via generating function, and thus
it is convenient to work on the dual side. We denote by ${\rm dec}:\mathfrak{X}\to\mathfrak{X}\otimes\mathfrak{X}$
the deconcatenation coproduct. Then $(\mathfrak{X},\mshuffle,{\rm dec})$
is a commutative Hopf algebra (see \cite{Kei_22_qs}). Let $\mathcal{R}=\mathcal{H}\langle\langle X_{1},X_{2},\dots\rangle\rangle$
be the dual $\mathcal{H}$-module of $\mathfrak{X}\otimes_{\mathbb{Q}}\mathcal{H}$.
Then $(\mathcal{R},\cdot,\mdelta)$ is a co-commutative Hopf $\mathcal{H}$-algebra
where $\cdot$ is concatenation and $\mdelta:\mathcal{R}\to\mathcal{R}\otimes_{\mathcal{H}}\mathcal{R}$
is defined by $\mdelta(uv)=\mdelta(u)\mdelta(v)$ and
\[
\mdelta X_{k}=X_{k}\otimes1+1\otimes X_{k}-\sum_{\substack{i,j,k'\geq1\\
i+j+k'=k
}
}(X_{i}\otimes X_{j})\cdot\mdelta(X_{k'}).
\]
For an index $\Bbbk=(k_{1},\dots,k_{d})\in\mathbb{Z}_{\geq1}^{d}$,
we put $x_{\Bbbk}=x_{k_{1}}\cdots x_{k_{d}}$ and $X_{\Bbbk}=X_{k_{1}}\cdots X_{k_{d}}$.
By definition,
\[
\mdelta\Phi^{\shuffle}=\sum_{\Bbbk,\Bbbk'}L_{B}^{\shuffle}(x_{\Bbbk}\mshuffle x_{\Bbbk'})X_{\Bbbk}\otimes X_{\Bbbk'}
\]
and
\[
\mdelta\Phi=\sum_{\Bbbk,\Bbbk'}L_{B}(x_{\Bbbk}\mshuffle x_{\Bbbk'})X_{\Bbbk}\otimes X_{\Bbbk'}
\]
where $\Bbbk$ and $\Bbbk'$ runs all indices with different parities.
\begin{lem}
\label{lem:LB=00003D0_for_ev1_case}We have
\[
L_{B}(u\mshuffle v)=0
\]
for $u\in\mathfrak{X}_{{\rm ev}}'$ and $v\in\mathfrak{X}_{{\rm od}}$.
\end{lem}

\begin{proof}
By the first form of the block shuffle identity (Theorem \ref{thm:BlockShuffleIdentity}),
we have
\[
\mdelta(\Phi^{\shuffle})\in\sum_{k=0}^{\infty}X_{1}^{k}\otimes\mathcal{R}+\sum_{k=0}^{\infty}\mathcal{R}\otimes X_{1}^{k}.
\]
Since $\mdelta(X_{1})=X_{1}\otimes1+1\otimes X_{1}$, 
\[
\mdelta(\Phi)=\mdelta(\Gamma_{1}(X_{1})^{-2})\cdot\mdelta(\Phi^{\shuffle})\cdot\mdelta(\Gamma_{1}(-X_{1})^{-2})\in\sum_{k=0}^{\infty}X_{1}^{k}\otimes\mathcal{R}+\sum_{k=0}^{\infty}\mathcal{R}\otimes X_{1}^{k}
\]
which implies $L_{B}(u\mshuffle v)=0$ for all $u\in\mathfrak{X}_{{\rm ev}}'$
and $v\in\mathfrak{X}_{{\rm od}}$.
\end{proof}
We define a derivation ${\rm D}$ on $\mathcal{A}(\{0,1\})$ (with
respect to the concatenation product) by
\[
{\rm D}(e_{a})=e_{a}e_{a}-e_{a}e_{0}-e_{1}e_{a}=-e_{1}e_{0}\quad(a\in\{0,1\}).
\]

\begin{lem}
\label{lem:Ltilde_Dw}For $w\in\mathcal{A}(\{0,1\})$, we have
\begin{align*}
\tilde{L}({\rm D}(w)) & =\sum_{N\geq2}(-1)^{N}\zeta^{\mathfrak{m}}(N)\tilde{L}({\rm Cut}_{(\{0\}^{N-1};\emptyset)}(w))+\sum_{N\geq2}\zeta^{\mathfrak{m}}(N)\tilde{L}({\rm Cut}_{(\emptyset;\{1\}^{N-1})}(w)).
\end{align*}
\end{lem}

\begin{proof}
Let $*$ be a binary operation on $\mathcal{A}(\{0,1\})$ defined
in \cite{HS2020_algdif}. First, note that 
\[
{\rm D}(w)=w*e_{1}-w\shuffle e_{1}
\]
for $w\in\mathcal{A}(\{0,1\})$ since
\begin{align*}
{\rm D}(e_{a_{1}}\cdots e_{a_{k}}) & =\sum_{i=1}^{k}e_{a_{1}}\cdots e_{a_{i}}e_{a_{i}}\cdots e_{a_{k}}-\sum_{i=1}^{k}e_{a_{1}}\cdots e_{a_{i}}e_{0}e_{a_{i+1}}\cdots e_{a_{k}}-\sum_{i=1}^{k}e_{a_{1}}\cdots e_{a_{i-1}}e_{1}e_{a_{i}}\cdots e_{a_{k}}\\
e_{a_{1}}\cdots e_{a_{k}}*e_{1} & =\sum_{i=1}^{k}e_{a_{1}}\cdots e_{a_{i}}e_{a_{i}}\cdots e_{a_{k}}-\sum_{i=1}^{k}e_{a_{1}}\cdots e_{a_{i}}e_{0}e_{a_{i+1}}\cdots e_{a_{k}}+e_{a_{1}}\cdots e_{a_{k}}e_{1},\\
e_{a_{1}}\cdots e_{a_{k}}\shuffle e_{1} & =\sum_{i=1}^{k}e_{a_{1}}\cdots e_{a_{i-1}}e_{1}e_{a_{i}}\cdots e_{a_{k}}+e_{a_{1}}\cdots e_{a_{k}}e_{1}.
\end{align*}
Note that
\[
\tilde{L}({\rm D}(w))=0
\]
for $w\in\mathbb{Q}\oplus e_{1}\mathcal{A}(\{0,1\})e_{0}$ by the
regularized double shuffle relation (for motivic multiple zeta values).
Let ${\rm reg}_{0}:\mathcal{A}(\{0,1\})\to\mathcal{A}^{1}(\{0,1\})\coloneqq\mathbb{Q}\oplus e_{1}\mathcal{A}(\{0,1\})$
be the unique $\shuffle$-homomorphism satisfying ${\rm reg}_{0}(e_{0})=0$
and ${\rm reg}_{0}(w)=w$ for $w\in\mathcal{A}^{1}(\{0,1\})$. By
comparing the coefficients of $A^{n}$ of the formula in \cite[Lemma 19]{HMS_regohno},
we have
\[
{\rm reg}_{0}(e_{0}^{n}(w_{1}*w_{2}))=\sum_{i+j=n}{\rm reg}_{0}(e_{0}^{i}w_{1})*{\rm reg}_{0}(e_{0}^{j}w_{2})
\]
for $w_{1},w_{2}\in\mathcal{A}^{1}(\{0,1\})$. Fix $u\in\mathcal{A}^{1}(\{0,1\})$
and $n\geq0$. Letting $w_{1}=u$ and $w_{2}=e_{1}$ in the equality
above, we obtain
\begin{align*}
{\rm reg}_{0}(e_{0}^{n}(u*e_{1})) & =\sum_{i+j=n}{\rm reg}_{0}(e_{0}^{i}u)*{\rm reg}_{0}(e_{0}^{j}e_{1})=\sum_{i+j=n}(-1)^{j}{\rm reg}_{0}(e_{0}^{i}u)*e_{1}e_{0}^{j}.
\end{align*}
Therefore,
\begin{align*}
\tilde{L}({\rm D}(e_{0}^{n}u)) & =\tilde{L}(e_{0}^{n}u*e_{1}-e_{0}^{n}u\shuffle e_{1})\\
 & =\tilde{L}(e_{0}^{n}u*e_{1})\\
 & =\tilde{L}(e_{0}^{n}(u*e_{1}))\\
 & =\sum_{i+j=n}(-1)^{j}\tilde{L}({\rm reg}_{0}(e_{0}^{i}u)*e_{1}e_{0}^{j}).
\end{align*}
Here, by the regularized double shuffle relation,
\begin{align*}
\sum_{\substack{i+j=n\\
j\geq1
}
}(-1)^{j}\tilde{L}({\rm reg}_{0}(e_{0}^{i}u)*e_{1}e_{0}^{j}) & =\sum_{\substack{i+j=n\\
j\geq1
}
}(-1)^{j}\tilde{L}({\rm reg}_{0}(e_{0}^{i}u))\tilde{L}(e_{1}e_{0}^{j})\\
 & =\sum_{\substack{i+j=n\\
j\geq1
}
}(-1)^{j+1}\tilde{L}(e_{0}^{i}u)\zeta^{\mathfrak{m}}(j+1)\\
 & =\sum_{N\geq2}(-1)^{N}\zeta^{\mathfrak{m}}(N)\tilde{L}({\rm Cut}_{(\{0\}^{N-1};\emptyset)}(e_{0}^{n}u)).
\end{align*}
On the other hand
\begin{align*}
\sum_{\substack{i+j=n\\
j=0
}
}(-1)^{j}\tilde{L}({\rm reg}_{0}(e_{0}^{i}u)*e_{1}e_{0}^{j}) & =\tilde{L}({\rm reg}_{0}(e_{0}^{n}u)*e_{1})=\tilde{L}({\rm D}({\rm reg}_{0}(e_{0}^{n}u))).
\end{align*}
Thus
\[
\tilde{L}({\rm D}(e_{0}^{n}u))=\sum_{N\geq2}(-1)^{N}\zeta^{\mathfrak{m}}(N)\tilde{L}({\rm Cut}_{(\{0\}^{N-1};\emptyset)}(e_{0}^{n}u))+\tilde{L}({\rm D}({\rm reg}_{0}(e_{0}^{n}u)))
\]
for any $u\in\mathcal{A}^{1}(\{0,1\})$ and $n\geq0$, or equivalently
\begin{equation}
\tilde{L}({\rm D}(w))=\sum_{N\geq2}(-1)^{N}\zeta^{\mathfrak{m}}(N)\tilde{L}({\rm Cut}_{(\{0\}^{N-1};\emptyset)}(w))+\tilde{L}({\rm D}({\rm reg}_{0}(w)))\label{eq:Ltilde_D}
\end{equation}
for any $w\in\mathcal{A}(\{0,1\})$. Let $\tau$ be the anti-automorphism
of $\mathcal{A}(\{0,1\})$ defined by $\tau(e_{0})=-e_{1}$ and $\tau(e_{1})=-e_{0}$.
By definition $\tau\circ{\rm D}(w)=-{\rm D}\circ\tau(w)$. Since

\[
\tilde{L}({\rm D}(\tau(w)))=-\tilde{L}(\tau(D(w)))=-\tilde{L}(D(w)),
\]
and
\[
\tilde{L}({\rm Cut}_{(\{0\}^{N-1};\emptyset)}(\tau(w)))=(-1)^{N-1}\tilde{L}\left(\tau\left({\rm Cut}_{(\emptyset;\{1\}^{N-1})}(w)\right)\right),
\]
we have
\[
-\tilde{L}({\rm D}(w))=-\sum_{N\geq2}\zeta^{\mathfrak{m}}(N)\tilde{L}({\rm Cut}_{(\emptyset;\{1\}^{N-1})}(w))+\tilde{L}({\rm D}({\rm reg}_{0}(\tau(w)))),
\]
and especially,
\begin{equation}
\tilde{L}({\rm D}(w))=\sum_{N\geq2}\zeta^{\mathfrak{m}}(N)\tilde{L}({\rm Cut}_{(\emptyset;\{1\}^{N-1})}(w))\label{eq:Ltilde_D_special}
\end{equation}
for $w\in\mathcal{A}^{1}$. By replacing $w$ with ${\rm reg}_{0}(w)$
in (\ref{eq:Ltilde_D_special}), the last term of the right-hand side
of (\ref{eq:Ltilde_D}) now becomes
\begin{align*}
\tilde{L}({\rm D}({\rm reg}_{0}(w))) & =\sum_{N\geq2}\zeta^{\mathfrak{m}}(N)\tilde{L}({\rm Cut}_{(\emptyset;\{1\}^{N-1})}({\rm reg}_{0}(w)))\\
 & =\sum_{N\geq2}\zeta^{\mathfrak{m}}(N)\tilde{L}({\rm reg}_{0}({\rm Cut}_{(\emptyset;\{1\}^{N-1})}(w)))\\
 & =\sum_{N\geq2}\zeta^{\mathfrak{m}}(N)\tilde{L}({\rm Cut}_{(\emptyset;\{1\}^{N-1})}(w)),
\end{align*}
and hence
\begin{align*}
\tilde{L}({\rm D}(w)) & =\sum_{N\geq2}(-1)^{N}\zeta^{\mathfrak{m}}(N)\tilde{L}({\rm Cut}_{(\{0\}^{N-1};\emptyset)}(w))+\tilde{L}({\rm D}({\rm reg}_{0}(w)))\\
 & =\sum_{N\geq2}(-1)^{N}\zeta^{\mathfrak{m}}(N)\tilde{L}({\rm Cut}_{(\{0\}^{N-1};\emptyset)}(w))+\sum_{N\geq2}\zeta^{\mathfrak{m}}(N)\tilde{L}({\rm Cut}_{(\emptyset;\{1\}^{N-1})}({\rm reg}_{0}(w))).\qedhere
\end{align*}
\end{proof}
\begin{lem}
\label{lem:LB=00003D0_for_x1_case}For $u\in\mathfrak{X}_{{\rm od}}$,
we have
\[
L_{B}(x_{1}\mshuffle u)=0.
\]
\end{lem}

\begin{proof}
For $a_{0},\dots,a_{k+1}\in\{0,1\}$ with $a_{0}=0$ and $a_{k+1}=1$,
\begin{align*}
x_{1}\dshuffle e_{a_{0}}\cdots e_{a_{k+1}} & =2\sum_{i=0}^{k+1}e_{a_{0}}\cdots e_{a_{i-1}}e_{a_{i}}e_{a_{i}}e_{a_{i+1}}\cdots e_{a_{k+1}}-\sum_{i=0}^{k}e_{a_{0}}\cdots e_{a_{i}}(e_{0}+e_{1})e_{a_{i+1}}\cdots e_{a_{k+1}}\\
 & =2\sum_{i=1}^{k}e_{a_{0}}\cdots e_{a_{i-1}}e_{a_{i}}e_{a_{i}}e_{a_{i+1}}\cdots e_{a_{k+1}}+2e_{a_{0}}e_{a_{0}}e_{a_{1}}\cdots e_{a_{k+1}}+2e_{a_{0}}\cdots e_{a_{k}}e_{a_{k+1}}e_{a_{k+1}}\\
 & \quad\quad-e_{0}(e_{a_{1}}\cdots e_{a_{k}}\shuffle(e_{0}+e_{1}))e_{1}\\
 & =2\sum_{i=1}^{k}e_{a_{0}}\cdots e_{a_{i-1}}\left({\rm D}(e_{a_{i}})+e_{a_{i}}e_{0}+e_{1}e_{a_{i}}\right)e_{a_{i+1}}\cdots e_{a_{k+1}}\\
 & \quad\quad+2e_{a_{0}}e_{a_{0}}e_{a_{1}}\cdots e_{a_{k+1}}+2e_{a_{0}}\cdots e_{a_{k}}e_{a_{k+1}}e_{a_{k+1}}-e_{0}(e_{a_{1}}\cdots e_{a_{k}}\shuffle(e_{0}+e_{1}))e_{1}\\
 & =2\sum_{i=1}^{k}e_{a_{0}}\cdots e_{a_{i-1}}{\rm D}(e_{a_{i}})e_{a_{i+1}}\cdots e_{a_{k+1}}+e_{0}(e_{a_{1}}\cdots e_{a_{k}}\shuffle(-e_{0}-e_{1}+2e_{0}+2e_{1}))e_{1}\\
 & =e_{0}\left(2{\rm D}(e_{a_{1}}\cdots e_{a_{k}})+e_{a_{1}}\cdots e_{a_{k}}\shuffle(e_{0}+e_{1})\right)e_{1}.
\end{align*}
In other words, we have
\[
x_{1}\dshuffle e_{0}ue_{1}=e_{0}\left(2{\rm D}(u)+u\shuffle(e_{0}+e_{1})\right)e_{1}
\]
for $u\in\mathcal{A}(\{0,1\})$. Thus, by Lemma \ref{lem:Ltilde_Dw},
\begin{align*}
\tilde{L}(\llbracket x_{1}\dshuffle e_{0}ue_{1}\rrbracket) & =2\tilde{L}({\rm D}(u))\\
 & =2\sum_{N\geq2}(-1)^{N}\zeta^{\mathfrak{m}}(N)\tilde{L}({\rm Cut}_{(\{0\}^{N-1};\emptyset)}(u))+2\sum_{N\geq2}\zeta^{\mathfrak{m}}(N)\tilde{L}({\rm Cut}_{(\emptyset;\{1\}^{N-1})}(u))
\end{align*}
for $u\in\mathcal{A}(\{0,1\})$, or equivalently,
\begin{equation}
L_{B}^{\shuffle}(x_{1}\mshuffle w)=2\sum_{N\geq2}(-1)^{N}\zeta^{\mathfrak{m}}(N)L_{B}^{\shuffle}({\rm cut}_{(x_{1}^{N-1};1)}(w))+2\sum_{N\geq2}\zeta^{\mathfrak{m}}(N)\tilde{L}({\rm cut}_{(1;x_{1}^{N-1})}(w))\label{eq:Lbsh(x1_w)}
\end{equation}
where ${\rm cut}_{(u;v)}:\mathfrak{X}\to\mathfrak{X}$ is a linear
map defined by
\[
{\rm cut}_{(u;v)}(w)\coloneqq\begin{cases}
w' & {\rm if}\ w=uw'v\\
0 & {\rm otherwise}
\end{cases}
\]
for monomials $u,v,w\in\mathfrak{X}$. Define $f:\mathcal{R}\to\mathcal{R}$
by the composition
\[
\mathcal{R}\xrightarrow{\mdelta}\mathcal{R}\widehat{\otimes}\mathcal{R}\xrightarrow{u_{1}\otimes u_{2}\mapsto\langle u_{1},x_{1}\rangle\cdot u_{2}}\mathcal{R}.
\]
Then by (\ref{eq:Lbsh(x1_w)}), we have
\begin{align*}
f(\Phi^{\shuffle}) & =\sum_{\Bbbk:\,{\rm odd}}L_{B}^{\shuffle}(x_{1}\mshuffle x_{\Bbbk})\cdot X_{\Bbbk}=2\psi_{1}(X_{1})\cdot\Phi^{\shuffle}-\Phi^{\shuffle}\cdot2\psi_{1}(-X_{1})
\end{align*}
where
\[
\psi_{1}(t)\coloneqq\frac{d}{dt}\log\Gamma_{1}(t)=-\sum_{n=2}^{\infty}\zeta^{\mathfrak{m}}(n)(-t)^{n-1}.
\]
Note that $f$ is an derivation on $R$, i.e.,
\[
f(w_{1}\cdot w_{2})=f(w_{1})w_{2}+w_{1}f(w_{2})\quad(w_{1},w_{2}\in\mathcal{R})
\]
and in particular
\[
f(p(X_{1}))=\frac{d}{dX_{1}}p(X_{1})\quad(p\in\mathcal{H}[[t]]).
\]
Thus for $\Phi_{p,q}\coloneqq p(X_{1})\Phi^{\shuffle}q(X_{1})$ with
general $p(t),q(t)\in1+t\mathcal{H}[[t]]$, we have
\begin{align*}
f(\Phi_{p,q}) & =\left(\frac{p'(X_{1})}{p(X_{1})}+2\psi_{1}(X_{1})\right)\cdot\Phi_{p,q}+\Phi_{p,q}\cdot\left(\frac{q'(X_{1})}{q(X_{1})}-2\psi_{1}(-X_{1})\right)\\
 & =\frac{d}{dX_{1}}\log\left(p(X_{1})\Gamma_{1}(X_{1})^{2}\right)\cdot\Phi_{p,q}+\Phi_{p,q}\cdot\frac{d}{dX_{1}}\log\left(q(X_{1})\Gamma_{1}(-X_{1})^{2}\right).
\end{align*}
Therefore, by the choice $p(t)=\Gamma_{1}(t)^{-2}$, $q(t)=\Gamma_{1}(-t)^{-2}$,
\[
f(\Phi)=0,
\]
which is equivalent to the claim of the lemma.
\end{proof}
With the lemma above, we can finally complete the proof of the second
form of the block shuffle identity (Theorem \ref{thm:BlockShuffleIdentity2})
as follows:
\begin{proof}[Proof of Theorem \ref{thm:BlockShuffleIdentity2}]
By Lemma \ref{lem:LB=00003D0_for_ev1_case}, it is enough to show
\[
L_{B}(x_{1}^{n}\mshuffle v)=0
\]
for $n\geq1$ and $v\in\mathfrak{X}_{{\rm od}}$. Since
\[
x_{1}^{n}-\frac{1}{n}x_{1}\mshuffle x_{1}^{n-1}\in\mathfrak{X}_{{\rm ev}}',
\]
by Lemma \ref{lem:LB=00003D0_for_ev1_case} and \ref{lem:LB=00003D0_for_x1_case},
we have
\[
L_{B}(x_{1}^{n}\mshuffle v)=\frac{1}{n}L_{B}(x_{1}\mshuffle x_{1}^{n-1}\mshuffle v)=0.\qedhere
\]
\end{proof}
\begin{rem}
In terms of the generating series, Theorem \ref{thm:BlockShuffleIdentity2}
is equivalent to the Lie-like property $\mdelta\Phi=\Phi\otimes1+1\otimes\Phi$.
\end{rem}

\subsection{Proof of the full version of Charlton's generalized cyclic insertion
conjecture}

In this section, we derive the full version of Charlton's generalized
cyclic insertion conjecture \cite[Conjecture 6.3]{Ch21} from Theorem
\ref{thm:BlockShuffleIdentity2}.
\begin{prop}
\label{prop:Charlton_full_ver}Suppose that $d>0$ and $l_{1},\ldots,l_{d}\geq1$
such that $\sum_{i=1}^{d}(l_{i}-1)$ is odd, and let
\[
A_{r}=\begin{cases}
\frac{2\mu^{r}}{(r+2)!} & r\text{: even},\\
0 & r\text{: odd},
\end{cases}
\]
where $\mu\in\mathcal{H}$ is the motivic $2\pi i$. Then, 
\[
\sum_{\substack{r=0}
}^{d}A_{r}\sum_{\sigma\in\mathbb{Z}/d\mathbb{Z}}\Big(\prod_{j=1}^{r}\delta_{l_{\sigma+j},1}\Big)L_{B}^{\shuffle}(x_{l_{\sigma+r+1}}\cdots x_{l_{\sigma+d}})=\begin{cases}
L_{B}^{\shuffle}(x_{l_{1}+\cdots+l_{d}}) & d:\,{\rm odd}\\
0 & d:\,{\rm even}.
\end{cases}
\]
Here, the subscripts of $l$'s are regarded as elements of $\mathbb{Z}/d\mathbb{Z}$.
\end{prop}

\begin{proof}
By Theorem \ref{thm:BlockShuffleIdentity2} and Proposition \ref{prop:cyclic_in_msh},
we have
\[
\sum_{\sigma\in\mathbb{Z}/d\mathbb{Z}}L_{B}(x_{l_{\sigma+1}}\cdots x_{l_{\sigma+d}})\equiv\begin{cases}
L_{B}(x_{l_{1}+\cdots+l_{d}}) & d:{\rm odd}\\
0 & d:{\rm even},
\end{cases}
\]
and since $L_{B}^{\shuffle}(x_{l})=L_{B}(x_{l})$ for $l\geq1$, it
suffices to show
\[
\sum_{\sigma\in\mathbb{Z}/d\mathbb{Z}}L_{B}(x_{l_{\sigma+1}}\cdots x_{l_{\sigma+d}})=\sum_{\substack{r=0}
}^{d}A_{r}\sum_{\sigma\in\mathbb{Z}/d\mathbb{Z}}\Big(\prod_{j=1}^{r}\delta_{l_{\sigma+j},1}\Big)L_{B}^{\shuffle}(x_{l_{\sigma+r+1}}\cdots x_{l_{\sigma+d}}).
\]
We define $c_{k},c_{k}'$ by $\Gamma_{1}(t)^{-2}=\sum_{k=0}^{\infty}c_{k}t^{k}$
and $\Gamma_{1}(-t)^{-2}=\sum_{k=0}^{\infty}c_{k}'t^{k}$. Then, by
definition, 
\[
L_{B}(x_{l_{1}}\cdots x_{l_{d}})=\sum_{\substack{0\leq p,q\\
p+q\leq d
}
}\Big(c_{p}\prod_{j=1}^{p}\delta_{l_{j},1}\Big)L_{B}^{\shuffle}(x_{l_{p+1}}\cdots x_{l_{d-q}})\Big(c'_{q}\prod_{j=d-q+1}^{d}\delta_{l_{j},1}\Big)
\]
Thus,
\begin{align*}
\sum_{\sigma\in\mathbb{Z}/d\mathbb{Z}}L_{B}(x_{l_{\sigma+1}}\cdots x_{l_{\sigma+d}}) & =\sum_{\sigma\in\mathbb{Z}/d\mathbb{Z}}\sum_{\substack{0\leq p,q\\
p+q\leq d
}
}c_{p}c'_{q}\Big(\prod_{-q<j\leq p}\delta_{l_{\sigma+j},1}\Big)L_{B}^{\shuffle}(x_{l_{\sigma+p+1}}\cdots x_{l_{\sigma+d-q}})\\
 & =\sum_{\sigma'\in\mathbb{Z}/d\mathbb{Z}}\sum_{\substack{0\leq p\leq r\leq d}
}c_{p}c'_{r-p}\Big(\prod_{0<j\leq r}\delta_{l_{\sigma'+j},1}\Big)L_{B}^{\shuffle}(x_{l_{\sigma'+r+1}}\cdots x_{l_{\sigma'+d}})\quad\left(\text{by }\sigma'=\sigma-q,\,r=p+q\right)
\end{align*}
Since $\sum_{\substack{0\leq p\leq r}
}c_{p}c_{r-p}'$ is the coefficient of $t^{r}$ of 
\[
\Gamma_{1}(t)^{-2}\Gamma_{1}(-t)^{-2}=\sum_{\substack{r=0}
}^{\infty}A_{r}t^{r},
\]
it follows that
\begin{align*}
\sum_{\sigma\in\mathbb{Z}/d\mathbb{Z}}L_{B}(x_{l_{\sigma+1}}\cdots x_{l_{\sigma+d}}) & =\sum_{\sigma'\in\mathbb{Z}/d\mathbb{Z}}\sum_{\substack{0\leq r\leq d}
}A_{r}\cdot\Big(\prod_{0<j\leq r}\delta_{l_{\sigma'+j},1}\Big)L_{B}^{\shuffle}(x_{l_{\sigma'+r+1}}\cdots x_{l_{\sigma'+d}}).\qedhere
\end{align*}
\end{proof}

\subsection{Block shuffle relation for refined symmetric multiple zeta values}

In this section, we discuss a closed path analog of the block shuffle
identity. Let $0'$ and $1'$ be the standard tangential basepoints
$\overrightarrow{1}_{0}$ and $-\overrightarrow{1}_{1}$, respectively.
Define a $\mathbb{Q}$-linear map
\[
L^{\circ,\shuffle}:\mathbb{Q}e_{0}\oplus e_{0}\mathcal{A}(\{0,1\})e_{0}\to\mathcal{H}
\]
by
\[
L^{\circ,\shuffle}(e_{0})=0\text{ and }L^{\circ,\shuffle}(e_{0}e_{a_{1}}\cdots e_{a_{k}}e_{0})=I_{\gamma}^{\mathfrak{m}}(0';a_{1},\dots,a_{k};0')-\delta_{k,0}
\]
where $\gamma$ is a closed path from $0'$ to $0'$ that runs along
the real axis to $1-\varepsilon$ with $0<\varepsilon<1$ and encircles
$1$ once counter-clockwisely. Now, define $L_{B}^{\circ,\shuffle}:\mathfrak{X}_{{\rm ev}}^{+}\to\mathcal{H}$
by $L^{\circ,\shuffle}\circ B$ and let
\[
\Phi^{\circ,\shuffle}\coloneqq\sum_{\Bbbk\neq\emptyset:\,{\rm even}}L_{B}^{\circ,\shuffle}(x_{\Bbbk})X_{\Bbbk}\quad\in\mathcal{R}.
\]
As before, we define $\Phi^{\circ}$ and the block regularization
$L_{B}^{\circ}:\mathfrak{X}_{{\rm ev}}^{+}\to\mathcal{H}$ by 
\[
\sum_{\Bbbk\neq\emptyset:\,{\rm even}}L_{B}^{\circ}(x_{\Bbbk})X_{\Bbbk}\coloneqq\Phi^{\circ}\coloneqq\Gamma_{1}(X_{1})^{-2}\Phi^{\circ,\shuffle}\Gamma_{1}(-X_{1})^{-2}.
\]
Based on a numerical investigation, we conjecture the following:
\begin{conjecture}
\label{conj:block_shuffle_for_symmetric_MZV}For $u,v\in\mathfrak{X}_{{\rm ev}}'\times\mathfrak{X}_{{\rm ev}}^{+}\cup\mathfrak{X}_{{\rm od}}\times\mathfrak{X}_{{\rm od}}$,
\[
L_{B}^{\circ}(u\mshuffle v)\overset{?}{=}\begin{cases}
0 & (u,v)\in\mathfrak{X}_{{\rm ev}}'\times\mathfrak{X}_{{\rm ev}}^{+}\\
\mu^{2}L_{B}(u)L_{B}(v) & (u,v)\in\mathfrak{X}_{{\rm od}}\times\mathfrak{X}_{{\rm od}}.
\end{cases}
\]
\end{conjecture}

As a supporting fact of this conjecture, here we prove that at least
its ``real part'' is true. Let $\tau$ denote the complex conjugation
on $\mathcal{H}$, $\mathcal{H}^{\tau,+}$ the real part $\mathcal{H}^{\tau,+}\coloneqq\left\{ x\in\mathcal{H}\bigr|\tau(x)=x\right\} $,
and $\mathfrak{r}:\mathcal{H}\rightarrow\mathcal{H}^{\tau,+}$ the
map defined by $\mathfrak{r}(x)\coloneqq x+\tau(x)$. By definition
$\tau(\mu)=-\mu$. We extend the definition of $\mathfrak{r}$ to
\[
\mathcal{R}\rightarrow\mathcal{R}^{\tau,+}\coloneqq\mathcal{H}^{\tau,+}\langle\langle X_{1},X_{2},\dots\rangle\rangle
\]
by $\mathfrak{r}(\sum_{\Bbbk}c_{\Bbbk}X_{\Bbbk})=\sum_{\Bbbk}\mathfrak{r}(c_{\Bbbk})X_{\Bbbk}$.
\begin{thm}
\label{thm:block_shuffle_for_symmetric_MZV}For $u,v\in\mathfrak{X}_{{\rm ev}}^{+}\times\mathfrak{X}_{{\rm ev}}^{+}\cup\mathfrak{X}_{{\rm od}}\times\mathfrak{X}_{{\rm od}}$,
we have
\[
\mathfrak{r}\bigl(L_{B}^{\circ}(u\mshuffle v)\bigr)=\begin{cases}
0 & u,v\in\mathfrak{X}_{{\rm ev}}^{+}\\
2\mu^{2}L_{B}(u)L_{B}(v) & u,v\in\mathfrak{X}_{{\rm od}}.
\end{cases}
\]
\end{thm}

\begin{proof}
We will prove the theorem in the equivalent form
\[
\mdelta(\mathfrak{r}\bigl(\Phi^{\circ}\bigr))=1\otimes\mathfrak{r}\bigl(\Phi^{\circ}\bigr)+2\mu^{2}\Phi\otimes\Phi+\mathfrak{r}\bigl(\Phi^{\circ}\bigr)\otimes1.
\]
First, notice that $\gamma=\mathrm{dch}\circ C_{1}\circ\mathrm{dch}^{-1}$
where $C_{1}$ is a closed path from $1'$ to $1'$ that encircles
$1$ once counter-clockwisely without encircling $0$, and $\mathrm{dch}(t)=t$
is a standard straight path from $0'$ to $1'$. Hence, by the path
composition formula,
\begin{align*}
I_{\gamma}^{\mathfrak{m}}(0';a_{1},\dots,a_{k};0') & =\sum_{0\leq i\leq j\leq k}I_{\mathrm{dch}}^{\mathfrak{m}}(0';a_{1},\dots,a_{i};1')I_{C_{1}}^{\mathfrak{m}}(1';a_{i+1},\dots,a_{j};1')I_{\mathrm{dch}^{-1}}^{\mathfrak{m}}(1';a_{j+1},\dots,a_{k};0').
\end{align*}
Since 
\[
I_{C_{1}}^{\mathfrak{m}}(1';a_{i+1},\dots,a_{j};1')=\begin{cases}
\frac{\mu^{j-i}}{(j-i)!} & \text{if }a_{i+1}=\cdots=a_{j}=1,\\
0 & \text{otherwise,}
\end{cases}
\]
and 
\[
\sum_{\substack{0\leq i\leq k}
}I_{\mathrm{dch}}^{\mathfrak{m}}(0';a_{1},a_{2},\dots,a_{i};1')I_{\mathrm{dch}^{-1}}^{\mathfrak{m}}(1';a_{i+1},\dots,a_{k};0')=I_{\mathrm{dch}\circ\mathrm{dch}^{-1}}^{\mathfrak{m}}(0';a_{1},\dots,a_{k};0')=\delta_{k,0},
\]
we deduce that
\begin{align*}
\mathfrak{r}\bigl(L^{\circ,\shuffle}(e_{0}e_{a_{1}}\cdots e_{a_{k}}e_{0})\bigr) & =2\sum_{\substack{0\leq i\leq j\leq k\\
j-i\geq2,\,\mathrm{even}\\
a_{i+1}=\cdots=a_{j}=1
}
}I_{\mathrm{dch}}^{\mathfrak{m}}(0';a_{1},\dots,a_{i};1')\frac{\mu^{j-i}}{(j-i)!}I_{\mathrm{dch}^{-1}}^{\mathfrak{m}}(1';a_{j+1},\dots,a_{k};0').
\end{align*}
Rewriting this equation in terms of $\Phi^{\circ,\shuffle}$, 
\begin{align*}
\mathfrak{r}\bigl(\Phi^{\circ,\shuffle}\bigr) & =2\sum_{\substack{\Bbbk,\Bbbk':\,{\rm odd}\\
r\geq0
}
}L_{B}^{\shuffle}(x_{\Bbbk})\frac{\mu^{2r+2}}{(2r+2)!}L_{B}^{\shuffle}(x_{\Bbbk'})X_{\Bbbk}X_{1}^{2r}X_{\Bbbk'}=\Phi^{\shuffle}\cdot\Bigg(\sum_{r=0}^{\infty}\frac{2\mu^{2r+2}}{(2r+2)!}X_{1}^{2r}\Bigg)\cdot\Phi^{\shuffle}.
\end{align*}
By noting
\[
\sum_{r=0}^{\infty}\frac{2\mu^{2r+2}}{(2r+2)!}X_{1}^{2r}=\left(\frac{e^{\mu X_{1}/2}-e^{-\mu X_{1}/2}}{X_{1}}\right)^{2}=\mu^{2}\Gamma_{1}(-X_{1})^{-2}\Gamma_{1}(X_{1})^{-2},
\]
we find that
\begin{align*}
\mathfrak{r}\bigl(\Phi^{\circ,\shuffle}\bigr) & =\mu^{2}\Phi^{\shuffle}\cdot\Gamma_{1}(-X_{1})^{-2}\Gamma_{1}(X_{1})^{-2}\cdot\Phi^{\shuffle},
\end{align*}
or equivalently, 
\begin{align*}
\mathfrak{r}\bigl(\Phi^{\circ}\bigr) & =\mu^{2}\Phi^{2}.
\end{align*}
Hence, by Theorem \ref{thm:BlockShuffleIdentity2}, it readily follows
that
\[
\mdelta(\mathfrak{r}\bigl(\Phi^{\circ}\bigr))=\mu^{2}\mdelta(\Phi)^{2}=\mu^{2}(1\otimes\Phi+\Phi\otimes1)^{2}=1\otimes\mathfrak{r}\bigl(\Phi^{\circ}\bigr)+2\mu^{2}\Phi\otimes\Phi+\mathfrak{r}\bigl(\Phi^{\circ}\bigr)\otimes1.\qedhere
\]
\end{proof}

\section{Several remarks on the block shuffle relation\label{sec:remarks}}

In this final section, we would like to make some remarks on the block
shuffle relation to discuss its significance. In the first part, we
will show that the block shuffle relation is an ultimate generalization
of Charlton's conjectures in a certain sense. Then in the second part,
we will discuss the space of the relations between $I_{{\rm bl}}^{\mathfrak{m}}$'s
of block degree ($=$ the number of entries minus one) at most two,
make some observations based on numerical experimentations, and state
some expectations about higher block degree case. 

\subsection{A certain maximality of the block shuffle relation}

Let us refer to a multiple zeta value $I_{{\rm bl}}^{\mathfrak{m}}(l_{1},\dots,l_{d})$
with $l_{1},\dots,l_{d}\geq2$ and $\#\{l_{i}:{\rm even}\}=1$ as
\emph{quasi-Hoffman}. One important feature of the block shuffle relation
is that it gives all linear relations among quasi-Hoffman multiple
zeta values. We have already mentioned this fact in Remark \ref{rem:Hoffman},
but we restate it as a theorem in light of its importance.
\begin{thm}
\label{thm:characterizing_BS}Define $\mathfrak{X}'\subset\mathfrak{X}$
as a direct sum
\[
\bigoplus_{d=1}^{\infty}\bigoplus_{\substack{l_{1},\ldots,l_{d}\geq2\\
\#\{l_{i}:{\rm even}\}=1
}
}\mathbb{Q}\cdot x_{l_{1}}\cdots x_{l_{d}}.
\]
Then, 
\[
\ker L_{B}\cap\mathfrak{X}'=\bigl(\mathfrak{X}^{+}\mshuffle\mathfrak{X}^{+}\bigr)\cap\mathfrak{X}'.
\]
In other words, all $\mathbb{Q}$-linear relations among quasi-Hoffman
multiple zeta values are $\mathbb{Q}$-linear combinations of the
block shuffle relations.
\end{thm}

This theorem implies that the block shuffle relations is an ``ultimate''
generalization of Conjectures \ref{Conj:Ch1} and \ref{Conj:Ch2}
in the following sense: Fix a positive integer $d$ and let $\mathfrak{Y}_{d}$
be the set of formal $\mathbb{Q}$-linear sums of the sequences $(S_{1},\dots,S_{r})$
where $S_{1},\dots,S_{r}$ are disjoint subsets of $\{1,\dots,d\}$
such that $S_{1}\sqcup\cdots\sqcup S_{r}=\{1,\dots,d\}$ and $\#S_{j}$
are odd. Define a $\mathbb{Q}$-linear map 
\[
M_{d}:\mathfrak{Y}_{d}\to\prod_{l_{1},\dots,l_{d}\geq2}\mathcal{H}
\]
by 
\[
M_{d}((S_{1},\dots,S_{r}))=\left(I_{{\rm bl}}^{\mathfrak{m}}(\sum_{j\in S_{1}}l_{j},\dots,\sum_{j\in S_{r}}l_{j})\right)_{l_{1},\dots,l_{d}\geq2},
\]
where we regard $I_{{\rm bl}}^{\mathfrak{m}}(b_{1},\dots,b_{k})$
as $0$ if $b_{1}+\cdots+b_{k}-k$ is even. For example, Conjectures
\ref{Conj:Ch1} and \ref{Conj:Ch2} say that
\[
\sum_{i=0}^{d-1}(\{i+1\},\dots,\{d\},\{1\},\dots,\{i\})-\begin{cases}
(\{1,\dots,d\}) & d:\,{\rm odd}\\
0 & d:\,{\rm even}
\end{cases}
\]
and
\[
\sum_{\sigma\in\mathfrak{S}_{n+1}}{\rm sgn}(\sigma)(\{\sigma(1)\},\{n+2\},\{\sigma(2)\},\{n+3\},\dots,\{\sigma(n)\},\{2n+1\},\{\sigma(n+1)\})
\]
are elements of $\ker M_{d}$ and $\ker M_{2n+1}$, respectively.
\begin{thm}
\label{thm:exhaust_formal}All elements of $\ker M_{d}$ come from
the block shuffle relation, i.e., 
\[
\ker M_{d}=\sum_{0<i<n}\sum_{\substack{S_{1}\sqcup\cdots\sqcup S_{n}=\{1,\dots,d\}\\
\#S_{j}:\,\mathrm{odd}
}
}\mathbb{Q}\cdot(S_{1},\dots,S_{i})\hat{\shuffle}(S_{i+1},\dots,S_{n}),
\]
where $\hat{\shuffle}$ is defined by
\[
(A_{1},\dots,A_{l})\hat{\shuffle}(B_{1},\dots,B_{m})=\sum_{n=0}^{l+m}(-1)^{(l+m-n)/2}\sum_{\substack{(f,g)\in\mathscr{D}_{n}}
}(C_{1}(f,g),\dots,C_{n}(f,g)),
\]
where $\mathscr{D}_{n}$ is the set that appeared in the definition
of $\mshuffle$ and $C_{i}(f,g)\subset\{1,\dots,d\}$ is the union
of $A_{j}$ for $j\in f^{-1}(i)$ and $B_{j}$ for $j\in g^{-1}(i)$.
\end{thm}

\begin{proof}
Let $l_{1}\geq2$ be an even number, and $l_{2},\dots,l_{d}\geq3$
odd numbers such that all subsets of $\{l_{1},\dots,l_{d}\}$ of odd
cardinality have distinct sums. Then, for $x\in\ker M_{d}$, the $(l_{1},\dots,l_{d})$-component
of $M_{d}(x)$ gives a relation among quasi-Hoffman multiple zeta
values, hence is a linear combination of the block shuffle relation
by Theorem \ref{thm:characterizing_BS}. This proves $\ker M_{d}\subset(\text{right-hand side})$.
\end{proof}
\begin{rem}
Let $\mathfrak{\widetilde{Y}}_{d}$ be the set of formal $\mathbb{Q}$-linear
sums of the sequence $(S_{1},\dots,S_{r})$ where $S_{1},\dots,S_{r}$
are disjoint subsets of $\{1,\dots,d\}$ such that $S_{1}\sqcup\cdots\sqcup S_{r}=\{1,\dots,d\}$
without parity conditions on $\#S_{j}$. Let us denote by $\widetilde{M}_{d}$
the natural extension of $M_{d}$ to $\mathfrak{\widetilde{Y}}_{d}$.
Then, a stronger statement
\[
\ker\widetilde{M}_{d}=\sum_{0<i<n}\sum_{\substack{S_{1}\sqcup\cdots\sqcup S_{n}=\{1,\dots,d\}}
}\mathbb{Q}\cdot(S_{1},\dots,S_{i})\hat{\shuffle}(S_{i+1},\dots,S_{n})
\]
seems to hold. However, we do not have a proof of this conjecture
at the moment.
\end{rem}

\subsection{The space of the relations among multiples zeta values of a fixed
block degree}

Let us see other aspects of the block shuffle identity. One of the
powerful ways to study motivic multiple zeta values is to use their
filtration structures. In particular, one of the reasonable research
policies is to build understanding in order from the lowest filtration.
Let $\mathcal{Z}^{\mathfrak{m}}$ be the $\mathbb{Q}$-linear space
spanned by all motivic multiple zeta values. There are two important
filtration structures on $\mathcal{Z}^{\mathfrak{m}}$: the depth
filtration and the block filtration. The latter is a filtration introduced
by Brown \cite{BrownLetter} based on Charlton's block notation defined
as
\[
B_{n}\mathcal{Z}^{\mathfrak{m}}=\left\langle I_{{\rm bl}}^{\mathfrak{m}}(l_{0},\dots,l_{m})\mid m\leq n\right\rangle _{\mathbb{Q}}.
\]
These two filtrations have very different flavors. For example, Brown
\cite{BrownLetter} showed that the block filtration coincides with
the coradical filtration \cite{BrownLetter} and thus has a relatively
simple dimension formula, whereas the depth filtration has more complicated
dimension formula conjectured by Broadhurst-Kreimer in \cite{BK97}.
This is one of the advantages of considering the block filtration.
On the other hand, the regularized double shuffle relation works well
with the depth filtration and provides a basic tool to study the depth
filtration, but does not seem to work well with the block filtration.\footnote{For example, the double shuffle relation for the product of depth
$m$ and depth $n$ motivic multiple zeta values produces relations
among motivic multiple zeta values of depth at most $m+n$, but the
double relation for the product of block degree $m$ and $n$ motivic
multiple zeta values does not have such a property (even for $m=n=0$
case, both the shuffle and harmonic products contain terms of arbitrarily
high block degrees).} Thus, in order to study the block filtration structure, we need to
search for other families of relations than the double shuffle relations. 

In the following, we will look at concrete examples of the relations
among motivic multiple zeta values in small block degrees. For simplicity,
let us specify ourselves to the admissible indices here. The block
degree zero case is trivial: for each even integer $w\geq0$, there
exists only one motivic multiple zeta value $I_{{\rm bl}}^{\mathfrak{m}}(w+2)$
of block degree zero and weight $w$, hence there are no relations.
Next, in block degree one case, there exist $w-1$ motivic multiple
zeta values $I_{{\rm bl}}^{\mathfrak{m}}(2,w),I_{{\rm bl}}^{\mathfrak{m}}(3,w-1),\dots,I_{{\rm bl}}^{\mathfrak{m}}(w,2)$
in each odd weight $w\geq1$. By Brown's theorem, $\{I(2a+1,w-2a+1)\mid1\leq a\leq(w-1)/2\}$
forms their basis, and so all the relations are exhausted by the duality
relations
\begin{equation}
I_{{\rm bl}}^{\mathfrak{m}}(a,b)+I_{{\rm bl}}^{\mathfrak{m}}(b,a)=0\quad\quad(a+b=w+2,\,a,b\geq2).\label{eq:degree_one_rel}
\end{equation}
Since (\ref{eq:degree_one_rel}) is a special case of the block shuffle
relation, we can say that all the $\mathbb{Q}$-linear relations among
motivic multiple zeta values of degree one come from the block shuffle
relations.

Now, let us consider the degree two case. For each even integer $w\geq2$,
there are ${w-1 \choose 2}$ indices of block degree two, which are
classified into the four patterns $I_{{\rm bl}}^{\mathfrak{m}}(\mathrm{even},\mathrm{odd},\mathrm{odd})$,
$I_{{\rm bl}}^{\mathfrak{m}}(\mathrm{odd},\mathrm{even},\mathrm{odd})$,
$I_{{\rm bl}}^{\mathfrak{m}}(\mathrm{odd},\mathrm{odd},\mathrm{even})$
and $I_{{\rm bl}}^{\mathfrak{m}}(\mathrm{even},\mathrm{even},\mathrm{even})$
according to the parity of their entries. The block shuffle relation
gives linear relations among the first three patterns and $I_{{\rm bl}}^{\mathfrak{m}}(w+2)$,
or those among the last pattern and $I_{{\rm bl}}^{\mathfrak{m}}(w+2)$.
By Theorem \ref{thm:characterizing_BS}, all the linear relation among
the first three patterns and $I_{{\rm bl}}^{\mathfrak{m}}(w+2)$ are
exhausted by the block shuffle relations except for those containing
the indices of the form $(a,1,b)$ with odd $a+b$. Moreover, such
exceptional values can be expressed as
\begin{equation}
I_{{\rm bl}}^{\mathfrak{m}}(a,1,b)=-I_{{\rm bl}}^{\mathfrak{m}}(a+b+1)-\sum_{j=2}^{a-1}I_{{\rm bl}}^{\mathfrak{m}}(j,a+1-j,b)-\sum_{j=2}^{b-1}I_{{\rm bl}}^{\mathfrak{m}}(a,b+1-j,j)\label{eq:Ibl(a,1,b)}
\end{equation}
by Hoffman's relation\footnote{It is also a consequence of the non-admissible case of the block shuffle
relation 
\[
I_{{\rm bl}}^{\mathfrak{m}}(1,a,b)+I_{{\rm bl}}^{\mathfrak{m}}(a,1,b)+I_{{\rm bl}}^{\mathfrak{m}}(a,b,1)-I_{{\rm bl}}^{\mathfrak{m}}(a+b+1)=0
\]
and the definitions of $I_{{\rm bl}}^{\mathfrak{m}}(1,a,b)$ and $I_{{\rm bl}}^{\mathfrak{m}}(a,b,1)$.} \cite[Theorem 5.1]{HoffmanConj}. Thus, all the linear relations
among $I_{{\rm bl}}^{\mathfrak{m}}(\mathrm{even},\mathrm{odd},\mathrm{odd})$,
$I_{{\rm bl}}^{\mathfrak{m}}(\mathrm{odd},\mathrm{even},\mathrm{odd})$,
$I_{{\rm bl}}^{\mathfrak{m}}(\mathrm{odd},\mathrm{odd},\mathrm{even})$
and $I_{{\rm bl}}^{\mathfrak{m}}(w+2)$ are exhausted by the block
shuffle relations and (\ref{eq:Ibl(a,1,b)}). How about the linear
relations among $I_{{\rm bl}}^{\mathfrak{m}}(\mathrm{even},\mathrm{even},\mathrm{even})$
and $I_{{\rm bl}}^{\mathfrak{m}}(w+2)$? Numerical experiments suggest
that all the linear relations among $I_{{\rm bl}}^{\mathfrak{m}}(\mathrm{even},\mathrm{even},\mathrm{even})$
and $I_{{\rm bl}}^{\mathfrak{m}}(w+2)$ are exhausted by the block
shuffle relations.

Finally, let us consider the ``mixed'' case, i.e., all the linear
relations among $I_{{\rm bl}}^{\mathfrak{m}}(\mathrm{even},\mathrm{odd},\mathrm{odd})$,
$I_{{\rm bl}}^{\mathfrak{m}}(\mathrm{odd},\mathrm{even},\mathrm{odd})$,
$I_{{\rm bl}}^{\mathfrak{m}}(\mathrm{odd},\mathrm{odd},\mathrm{even})$,
$I_{{\rm bl}}^{\mathfrak{m}}(\mathrm{even},\mathrm{even},\mathrm{even})$,
and $I_{{\rm bl}}^{\mathfrak{m}}(w+2)$. The block shuffle relation
and the relation (\ref{eq:Ibl(a,1,b)}) reduces any $I_{{\rm bl}}^{\mathfrak{m}}(a,b,c)$
with $a+b+c=w+2$ into a linear combination of
\begin{equation}
\{I_{{\rm bl}}^{\mathfrak{m}}(a,b,c)\mid(a,b,c)\in\mathbb{I}_{\mathrm{ooe}}^{(w)}\cup\mathbb{I}_{\mathrm{eee}}^{(w)}\}\cup\{I_{{\rm bl}}^{\mathfrak{m}}(w+2)\}\label{eq:generators}
\end{equation}
where 
\begin{align*}
\mathbb{I}_{\mathrm{ooe}}^{(w)} & \coloneqq\left\{ (a,b,c)\in\mathbb{Z}_{>1}^{3}\mid(a,b,c)\equiv(1,1,0)\bmod{2},\,a+b+c=w+2\right\} ,\\
\mathbb{I}_{\mathrm{eee}}^{(w)} & \coloneqq\left\{ (a,b,c)\in\mathbb{Z}_{>1}^{3}\mid(a,b,c)\equiv(0,0,0)\bmod{2},\,a\leq b,\,a<c,\,a+b+c=w+2\right\} ,
\end{align*}
and vica versa does not give any linear relation among the elements
of (\ref{eq:generators}). Thus, finding linear relations independent
of the block shuffle relations and the relations (\ref{eq:Ibl(a,1,b)})
is equivalent to finding non-trivial linear relations among (\ref{eq:generators}).
By Brown's theorem, 
\begin{equation}
\{I_{{\rm bl}}^{\mathfrak{m}}(a,b,c)\mid(a,b,c)\in\mathbb{I}_{\mathrm{ooe}}^{(w)}\}\cup\{I_{{\rm bl}}^{\mathfrak{m}}(w+2)\}\label{eq:Hoffman_basis}
\end{equation}
($=$ the Hoffman basis with at most two $3$'s) forms a basis of
the $\mathbb{Q}$-linear space spanned by (\ref{eq:generators}),
and thus we see that there are exactly $\#\mathbb{I}_{\mathrm{eee}}^{(w)}$
dimensions of such $\mathbb{Q}$-linear relations. Now, what exactly
are they? Are there simple relations like the block shuffle relations
or the relations (\ref{eq:Ibl(a,1,b)})? One particular feature of
those relations is the smallness of the coefficients, so let us investigate
the relations from the perspective of the size of their coefficients.
We can search for the relations with small coefficients by computer
in the following way. First, by computing Brown's infinitesimal coaction
and solving a system of linear equations, we can explicitly calculate
a basis of the $\mathbb{Z}$-module of all the linear relations (congruences)
among
\[
\{I_{{\rm bl}}^{\mathfrak{m}}(a,b,c)\mid(a,b,c)\in\mathbb{I}_{\mathrm{ooe}}^{(w)}\cup\mathbb{I}_{\mathrm{eee}}^{(w)}\}
\]
modulo $\mathbb{Q}I_{{\rm bl}}^{\mathfrak{m}}(w+2)$. Then, by applying
the LLL-algorithm to the obtained basis, we can find a new basis consisting
of (approximately) shortest vectors\footnote{Given an integer relation $\sum c_{a,b,c}I_{{\rm bl}}^{\mathfrak{m}}(a,b,c)\equiv0$
(mod $I_{{\rm bl}}^{\mathfrak{m}}(w+2)$), we identify it with the
vector $(c_{a,b,c})$ and define its norm as the square root of the
sum of the squares of $c_{a,b,c}$'s. The adjective ``short'' then
means having a small norm.}. For instance, for weights $40$, $50$, $60$ and $70$, the norms
of the congruences obtained by the algorithm explained above are as
follows\footnote{The outcome of the computation may slightly depend on the implementation
of the LLL-algorithm and the setting of parameters.}:
\begin{itemize}
\item For $w=40$, the norms of the first and second shortest vectors in
the obtained basis are $1.28\dots\times10^{4}$ and $1.17\dots\times10^{7}$.
The norms of all the other $61$ vectors in the obtained basis lie
between $2.08\dots\times10^{13}$ and $2.98\dots\times10^{13}$.
\item For $w=50$, the norms of the first and second shortest vectors in
the obtained basis are $2.55\dots\times10^{4}$ and $3.73\dots\times10^{7}$.
The norms of all the other $98$ vectors in the obtained basis lie
between $2.32\dots\times10^{17}$ and $3.47\dots\times10^{17}$.
\item For $w=60$, the norms of the first and second shortest vectors in
the obtained basis are $4.45\dots\times10^{4}$ and $9.53\dots\times10^{7}$.
The norms of all the other $143$ vectors in the obtained basis lie
between $3.15\dots\times10^{21}$ and $5.52\dots\times10^{21}$.
\item For $w=70$, the norms of the first and second shortest vectors in
the obtained basis are $7.12\dots\times10^{4}$ and $2.09\dots\times10^{8}$.
The norms of all the other $196$ vectors in the obtained basis lie
between $7.34\dots\times10^{25}$ and $1.13\dots\times10^{26}$.
\end{itemize}
The data above suggest the existence of a two-dimensional family of
relations with conspicuously small norms in each weight. By careful
observation of the coefficients of those relations, we can find their
exact forms as follows (For the simplicity of the expression, we also
used other multiple zeta values than those in (\ref{eq:generators}))\footnote{Of course, we can also rewrite the relations as those among (\ref{eq:generators})
using the block shuffle relation and the relation (\ref{eq:Ibl(a,1,b)})
(and also the definition of $I_{{\rm bl}}^{\mathfrak{m}}(1,3,w-2)$
for the first equation), but the descriptions get more complicated.
For example, the first formula (\ref{eq:exceptional_relation}) is
equivalent to 
\begin{align*}
p_{w+2}\cdot I_{{\rm bl}}^{\mathfrak{m}}(w+2)\overset{?}{=} & 3\sum_{\substack{(a,b,c)\in\mathbb{I}_{\mathrm{ooe}}^{(w)}}
}\left(p_{a,1,b}-2\delta_{a,3}-3\delta_{b,3}\right)\cdot I_{{\rm bl}}^{\mathfrak{m}}(a,b,c)\\
 & -3\sum_{\substack{\substack{(a,b,c)\in\mathbb{I}_{\mathrm{eee}}^{(w)}}
}
}p_{a,b,c}\cdot I_{{\rm bl}}^{\mathfrak{m}}(a,b,c)-3\sum_{\substack{\substack{(a,a,c)\in\mathbb{I}_{\mathrm{eee}}^{(w)}}
}
}\left(p_{c,a,a}-3\delta_{a,2}\right)\cdot I_{{\rm bl}}^{\mathfrak{m}}(a,a,c),
\end{align*}
where $p_{n}\in\mathbb{Z}_{>0}$ is given by
\[
p_{n}\coloneqq\frac{41}{216}n^{4}-\frac{91}{108}n^{3}+\frac{17}{12}n^{2}+\begin{cases}
-\frac{11}{3}n+9 & n\equiv0\bmod{6}\\
-\frac{107}{27}n+\frac{161}{27} & n\equiv2\bmod{6}\\
-\frac{139}{27}n+\frac{169}{27} & n\equiv4\bmod{6}
\end{cases}
\]
for $n>2$.}:
\begin{conjecture}
\label{conj:except_rel}Let $p_{a,b,c}$ and $q_{a,b,c}$ be polynomials
in $a,b,c$ defined by 
\[
p_{a,b,c}=(a-b)(a+b-6c+3)
\]
and
\begin{align*}
q_{a,b,c} & =\frac{1}{12}(a-b)\Big(17(a^{3}+b^{3})+197ab(a+b)+12(a^{2}+b^{2}-16ab)c-222(a+b)c^{2}-168c^{3}\\
 & \quad-362(a^{2}+b^{2})-1238ab+798(a+b)c+1674c^{2}+1075(a+b)-3774c+590\Big),
\end{align*}
which has symmetries $p_{a,b,c}+p_{b,a,c}=p_{a,b,c}+p_{c,a,b}+p_{b,c,a}=0$
and $q_{a,b,c}+q_{b,a,c}=q_{a,b,c}+q_{c,a,b}+q_{b,c,a}=0$.
\begin{enumerate}[label=\roman{enumi})]
\item For an even integer $w\geq0$,
\begin{align}
 & 3\sum_{\substack{\substack{(a,b,c)\equiv(1,1,0)\bmod{2}\\
a+b+c=w+2,\,b>1
}
}
}(p_{a,1,b}+\delta_{a,1}\delta_{b,3})\cdot I_{{\rm bl}}^{\mathfrak{m}}(a,b,c)-\sum_{\substack{\substack{(a,b,c)\equiv(0,0,0)\bmod{2}\\
a+b+c=w+2
}
}
}p_{a,b,c}\cdot I_{{\rm bl}}^{\mathfrak{m}}(a,b,c)\nonumber \\
 & \overset{?}{=}\frac{w^{2}(w-2)(5w-11)}{24}I_{{\rm bl}}^{\mathfrak{m}}(w+2).\label{eq:exceptional_relation}
\end{align}
\item For an even integer $w\geq4$,
\begin{align}
 & 3\sum_{\substack{(a,b,c)\equiv(1,1,0)\bmod{2}\\
a+b+c=w+2,a,b>1
}
}r{}_{a,b,c}I_{{\rm bl}}^{\mathfrak{m}}(a,b,c)-\sum_{\substack{(a,b,c)\equiv(0,0,0)\bmod{2}\\
a+b+c=w+2
}
}r_{a,b,c}'I_{{\rm bl}}^{\mathfrak{m}}(a,b,c)\nonumber \\
 & \overset{?}{=}\frac{-220w^{6}+2751w^{5}-10375w^{4}+16620w^{3}-35620w^{2}-7536w+285120}{1440}I_{{\rm bl}}^{\mathfrak{m}}(w+2)\label{eq:exceptional_relation2}
\end{align}
where $r_{a,b,c}$ and $r_{a,b,c}'$ are integers defined by
\begin{align*}
r_{a,b,c} & \coloneqq q_{a,1,b}+(a-1)(c-1)(a^{2}+19ab-3b^{2}+ac-6bc-26a-17b+4c+30)\\
 & \quad+26\delta_{a,5}+39\delta_{b,5}-21\delta_{a,3}\delta_{b,3}+2(b+c-b^{2}+2bc-c^{2})\delta_{a,3}+3c(3-c)\delta_{b,3}
\end{align*}
and
\[
r_{a,b,c}'\coloneqq q_{a,b,c}+18\delta_{a,2}\delta_{b,2}+63\delta_{a,2}\delta_{b,4}+117\delta_{a,4}\delta_{b,2}.
\]
\end{enumerate}
\end{conjecture}

As the weight $w$ increases, the norms of the block shuffle relations,
the relations (\ref{eq:Ibl(a,1,b)}), and the conjectural relations
(\ref{eq:exceptional_relation}), (\ref{eq:exceptional_relation2})
increase only at a rate of polynomial order of $w$, whereas the data
seem to suggest that the norm of any other ($=$ linearly independent)
relation of block degree $2$ increases at a higher speed than any
polynomial function of $w$.

Let us estimate the dimension of the relations of each kind. As $w\to\infty$,
we have approximately $\frac{1}{2}w^{2}$ generators out of which
$\#\mathbb{I}_{\mathrm{ooe}}\fallingdotseq\frac{1}{8}w^{2}$ are independent,
by the fact that (\ref{eq:Hoffman_basis}) forms a basis. Thus, the
dimension of all the relations among $\{I_{{\rm bl}}^{\mathfrak{m}}(a,b,c)\mid a+b+c=w+2\}\cup\{I_{{\rm bl}}^{\mathfrak{m}}(w+2)\}$
is approximately $\frac{3}{8}w^{2}$. Now:
\begin{itemize}
\item The block shuffle relations have approximately $\frac{1}{3}w^{2}$
dimensions (which make up $\frac{8}{9}$ of the dimension of all relations).
\item The dimension of the exceptional relations (\ref{eq:Ibl(a,1,b)}),
(\ref{eq:exceptional_relation}), (\ref{eq:exceptional_relation2})
is about $w$ in total.
\item Hence, the remaining mysterious relations have approximately $\frac{1}{24}w^{2}$
dimensions. It seems that such relations have ``very large'' norms
as we have explained above.
\end{itemize}
Thus, we can say that roughly $\frac{8}{9}$ (resp. $0$, $\frac{1}{9}$)
of the relations in block degree two comes from the block shuffle
relations (resp. the exceptional relations, the remaining relations)
as $w$ grows large.

More generally, let us consider the higher block degree case. Fix
a block degree $n\geq3$. Then for a large $w$, we have approximately
$g_{n}(w)\coloneqq\frac{w^{n}}{n!}$ generators and $\frac{1}{2^{n}}g_{n}(w)$
basis. Thus, the dimension of the space of all relations is approximately
$(1-\frac{1}{2^{n}})g_{n}(w)$, out of which the block shuffle relations
span around $\frac{n}{n+1}g_{n}(w)$ dimensions. Thus, the space of
the remaining relations has approximately $(\frac{1}{n+1}-\frac{1}{2^{n}})g_{n}(w)$
dimensions. Although it is merely an optimistic guess, we expect that
the situation is similar to the block degree two case so that the
block shuffle relations make up a majority of the relations having
``small'' coefficients. More precisely, given a weight $w$ and a
real number $\alpha>0$, let $V^{(w)}$ be the $\mathbb{Z}$-module
spanned by all the integer relations of weight $w$ and block degree
at most $n$, $V_{1}^{(w)}$ and $V_{2}^{(w)}$ its $\mathbb{Z}$-submodules
spanned by the block shuffle relations, and all the relations whose
norms are bounded by $w^{\alpha}$, respectively. Then, for sufficiently
large $w$, $V_{1}^{(w)}\subset V_{2}^{(w)}$, and our expectation
is that the elements of $V_{2}^{(w)}$ which are not in $V_{1}^{(w)}$
are ``relatively rare'' i.e., the ratio $\dim(V_{1}^{(w)}):\dim(V_{2}^{(w)}/V_{1}^{(w)}):\dim(V^{(w)}/V_{2}^{(w)})$
converges to $\frac{n}{n+1}:0:\frac{1}{n+1}-\frac{1}{2^{n}}$ as $w\rightarrow\infty$
for any $\alpha$. 

\subsection*{Acknowledgements}

This work was supported by JSPS KAKENHI Grant Numbers JP18J00982,
JP18K13392, JP19J00835, JP20K14293, JP22K03244, and MOST Grant Number
111-2115-M-002-003-MY3.

\bibliographystyle{plain}
\bibliography{Block_shuffle}

\end{document}